\documentclass[journal,article,submit,moreauthors,pdftex,10pt,a4paper]{mdpi-arxiv}

\usepackage{etoolbox}
\newtoggle{insertremarks}
\toggletrue{insertremarks}

\usepackage{amscd}
\usepackage{amssymb}
\usepackage{amsmath}
\usepackage{mathrsfs}
\usepackage{amsfonts}
\usepackage{amssymb}

\DeclareMathAlphabet\gothic{U}{euf}{m}{n}
\firstpage{1}
\articlenumber{x}
\doinum{10.3390/------}
\pubvolume{xx}
\pubyear{2018}
\copyrightyear{2018}
\history{Received: date; Accepted: date; Published: date}

\iftoggle{insertremarks}{
   \usepackage{fancyhdr}
	\pagestyle{fancy}
	\fancyhf{} 
	\rhead{\thepage\ of \pageref*{LastPage}}
	\lfoot{}

	\fancypagestyle{plain}
	{
	   \fancyhf{}
	   \fancyfoot{}
	}
}{%
}

\Title{Fourier Transform on the Homogeneous Space of 3D Positions and Orientations for Exact Solutions to
Linear Parabolic and (Hypo-)Elliptic PDEs}

\newcommand{\ProjS}{\mathbb{P}_p^{sym}}
\newcommand{\mslim}{\operatorname{ms-lim}\nolimits}

\newcommand\eqdef{\mathrel{\overset{\makebox[0pt]{\mbox{\normalfont\tiny\sffamily def}}}{=}}}
\newcommand{\hG}{\hat{G}}
\newcommand{\quot}{\mathbb{R}^3\rtimes S^2}
\newcommand{\eps}{\varepsilon}
\newcommand{\RR}{\mathbb{R}}

\newcommand{\R}{\mathbb{R}}
\newcommand{\mF}{\mathcal{F}}
\newcommand{\bR}{\mathbf{R}}
\newcommand{\bu}{\mathbf{u}}
\newcommand{\be}{\mathbf{e}}
\newcommand{\desda}{\Leftrightarrow}

\newcommand{\sumlm}{\sum \limits_{l=0}^{\infty} \sum \limits_{m=-l}^{l}}

\newcommand{\Philmom}{\Phi^{l,m}_{\bomega}}

\newcommand{\lambdalmr}{\lambda^{l,m}_{r}}

\newcommand{\bomega}{{\boldsymbol{\omega}}}
\newcommand{\rmd}{\mathrm{d}}

\newcommand{\tK}{\tilde{K}}

\newcommand{\bx}{\mathbf{x}}

\newcommand{\bn}{\mathbf{n}}
\newcommand{\ba}{\mathbf{a}}

\newcommand{\bN}{\mathbf{N}}
\newcommand{\bX}{\mathbf{X}}
\newcommand{\bG}{\mathbf{G}}
\newcommand{\bP}{\mathbf{P}}

\newcommand{\bW}{\mathbf{W}}
\newcommand{\bI}{\mathbf{I}}
\newcommand{\bLam}{\boldsymbol{\Lambda}}

\Author{Remco Duits\orcidA{} and Erik J. Bekkers and Alexey Mashtakov\orcidC{} }

\AuthorNames{Remco Duits and Erik J. Bekkers and Alexey Mashtakov}

\address{%
Department of Mathematics and Computer Science (CASA) Eindhoven University of Technology}
\corres{Correspondence: R.Duits@tue.nl; Tel.: +31-40-247-2859.}

\firstnote{Current address: Dep. of Mathematics and Computer Science (CASA) Eindhoven University of Technology}
\abstract{
Fokker-Planck PDEs (incl. diffusions) for stable L\'{e}vy processes (incl. Wiener processes) on the joint space of positions and orientations play a major role in mechanics, robotics, image analysis, directional statistics and probability theory. Exact analytic designs and solutions are known in the 2D case, where they have been obtained using Fourier transform on $SE(2)$. Here we extend these approaches to 3D using Fourier transform on the Lie group $SE(3)$ of rigid body motions. 
More precisely, we define the homogeneous space of 3D positions and orientations $\mathbb{R}^{3}\rtimes S^{2}:=SE(3)/(\{\mathbf{0}\} \times SO(2))$ as the quotient in $SE(3)$. In our construction, two group elements are equivalent if they are equal up to a rotation around the reference axis. On this quotient we design a specific Fourier transform. We apply this Fourier transform to derive new exact solutions to Fokker-Planck PDEs of $\alpha$-stable L\'{e}vy processes on $\quot$. This reduces classical analysis computations and provides an explicit algebraic spectral decomposition of the solutions. We compare the exact probability kernel for $\alpha = 1$ (the diffusion kernel) to the kernel for $\alpha=\frac12$ (the Poisson kernel). We set up SDEs for the  L\'{e}vy processes on the quotient and derive corresponding Monte-Carlo methods. We verify that the exact probability kernels arise as the limit of the Monte-Carlo approximations. }

\keyword{Fourier Transform; Rigid Body Motions; Partial Differential Equations; L\'{e}vy Processes; Lie Groups; Homogeneous Spaces; Stochastic Differential Equations}

\begin{document}

\section{Introduction}
The Fourier transform has had a tremendous impact on various fields of mathematics including
analysis, algebra and probability theory. It has a broad range of applied fields such as
signal and image processing, quantum mechanics, classical mechanics, robotics and system theory.
Thanks to Jean-Baptiste Joseph Fourier (1768-1830), who published
his pioneering work ``Th\'{e}ory analytique de la chaleur'' in 1822,
the effective technique of using a Fourier transform to solve linear PDE-systems (with appropriate boundary conditions) for heat transfer evolutions on compact subsets $\Omega$ of $\R^d$ was born.
The Fourier series representations of the solutions helped to understand the physics of heat transfer. Due to the linearity of the evolution operator that maps the possibly discontinuous square integrable initial condition to the square integrable solution at a fixed time $t>0$, one can apply a spectral decomposition which shows how each eigenfunction is dampened over time. Thanks to contributions of Johann Peter Gustav Lejeune Dirichlet (1805-1859),
completeness of the Fourier basis could then be formalized for several boundary conditions. Indeed separation of variables (also known as `the Fourier method') directly provides a Sturm-Liouville problem \cite{Zettl} and an orthonormal basis of eigenfunctions for $\mathbb{L}_{2}(\Omega)$, which is complete due to compactness of the associated self-adjoint kernel operator. When dilating the subset $\Omega$ to the full space $\R^d$, the
discrete set of eigenvalues start to fill $\R$ and the discrete spectrum approximates a continuous spectrum (see e.g.\cite{kato_operators_1976}). Then a diffusion system on $\R^d$ can be solved via a unitary Fourier transform on $\mathbb{L}_{2}(\R^d)$, cf.~\!\cite{Rudin}. 

Nowadays, in fields such as mechanics/robotics \cite{chirikjian_engineering_2000,chirikjian_stochastic_2011,Saccon,Nijmijer}, mathematical physics/harmonic analysis \cite{Alibook}, machine learning \cite{BekkersMICCAI,bekkers_template_2018,cohen2018intertwiners,Wellink,Mallat} 
and image analysis \cite{duits_image_2006,citti_cortical_2006,DuitsACHA,prandigauthierbook,Janssen2018,boscain_anthropomorphic_2012}
it is worthwhile to extend the spatial domain of functions on $M=\R^d$
(or $M=\mathbb{Z}^{d}$) to groups $G= M \rtimes T$ that are the semi-direct product of an Abelian
group $M$ 
and another matrix group $T$.  This requires a generalization of the Fourier transforms on the Lie group $(\mathbb{R}^{d},+)$ towards the groups $G= \R^d \rtimes T$.
Then the Fourier transform gives rise to an invertible decomposition of a square integrable function into irreducible representations. This is a powerful mechanism in view of the Schur's lemma \cite{Schur,Dieudonne} and spectral decompositions \cite{folland_course_1994,agrachev_intrinsic_2009}. However, it typically involves regularity constraints \cite[ch:3.6]{fuhr_abstract_2005}, \cite[ch:6.6]{folland_course_1994} on the structure of the dual orbits
in order that Mackey's imprimitivity theory \cite{mackey_imprimitivity_1949} can be applied to characterize \emph{all} unitary irreducible representations (UIRs) of $G$. This sets the Fourier transform on the Lie group $G$, \cite{fuhr_abstract_2005,folland_course_1994,sugiura_unitary_1990}.
Here, we omit technicalities on regularity constraints on the dual orbits and the fact that $G$ may not be of type I
(i.e. the quasi-dual group of $G$ may not be equal to the dual group of $G$, cf.~\!\cite[thm.7.6, 7.7]{folland_course_1994}, \cite{Dixmier}, \cite[ch:3]{fuhr_abstract_2005}), as this does not play a role in our case of interest.

We are concerned with the case $M=\R^3$ and $T=SO(3)$ where $G=SE(3)= M \rtimes SO(3)$ is the Lie group of 3D rigid body motions. It is a (type I) Lie group with an explicit Fourier transform $\mathcal{F}_{G}$ where the irreducible representations are determined by regular dual orbits (which are spheres in the Fourier domain indexed by their radius $p > 0$) and an integer index $s \in \mathbb{Z}$, cf.~\cite{sugiura_unitary_1990,chirikjian_engineering_2000}.

In this article we follow the idea of Joseph Fourier: we apply the Fourier transform $\mathcal{F}_{G}$ on the rigid body motion group $G=SE(3)$ to solve hypo-elliptic and elliptic heat flow evolutions respectively on the Lie group $G$. More precisely, we design a Fourier transform $\mathcal{F}_{G/H}$ on the homogeneous space of positions and orientations $G/H$ with $H\equiv \{\mathbf{0}\} \times SO(2)$ to solve hypo-elliptic and elliptic heat flow evolutions on the homogeneous space $G/H$. We also simultaneously solve related PDEs (beyond the diffusion case) as we will explain below.

The idea of applying Fourier transforms to solve linear (hypo-elliptic) PDEs on non-commutative groups of the type $\R^{d} \rtimes T$ is common and has been studied by many researchers. For example, tangible probability kernels for heat transfer (and fundamental solutions) on the Heisenberg group were derived by Gaveau \cite{Gaveau1977}. They can be derived by application \cite[ch:4.1.1]{agrachev_intrinsic_2009} of the Fourier transform on the Heisenberg group \cite[ch:1]{folland_course_1994}.
This also applies to probability kernels for hypo-elliptic diffusions on $SE(2)=\R^{2} \rtimes SO(2)$,
where 3 different types (a Fourier series, a rapidly decaying series, and a single expression) of explicit solutions to probability kernels for (convection)-diffusions are derived in previous works by Duits et al. \cite{duits_explicit_2008-1,duits_line_2009,DuitsCASA2005,duits_left-invariant_2010-1}. For a concise review see \cite[ch:5.1]{zhang_numerical_2016}. Here, the two fundamental models for contour perception by respectively
Mumford \cite{mumford_elastica_1994-1}, Petitot~\cite{petitot_neurogeometry_2003} and Citti \& Sarti \cite{citti_cortical_2006} formed great sources of inspiration to study the hypo-elliptic diffusion problem on $SE(2)$.

The hypo-elliptic diffusion kernel formula in terms
of a Fourier series representation was generalized to the much more wide setting of unimodular Lie groups by Agrachev, Boscain, Gauthier and Rossi \cite{agrachev_intrinsic_2009-1}.
This approach was then pursued by Portegies \& Duits (with additional classical analysis techniques) to achieve explicit exact solutions to (hypo)-elliptic (convection)-diffusions on the particular $SE(3)$ case, see  \cite{portegies_new_2017}.

Here, we structure the results in \cite{portegies_new_2017} and we present new, simpler formulas (relying on an algebraic approach rather than a classical analysis approach).
To this end we first set up a specific Fourier transform on the homogeneous space of positions and orientations in Theorem~\ref{corr:1}. Then we use it to derive explicit spectral decompositions of the evolution operator in Theorem~\ref{th:decomposition}, from which we deduce explicit new kernel expressions in Theorem~\ref{th:three}.
Finally, we generalize the exact solutions to other PDE systems beyond the diffusion case:
We will simultaneously solve the Forward-Kolmogorov PDEs for $\alpha$-stable L\'{e}vy processes on the homogeneous space of positions and orientations. Next we address their relevance in the fields of image analysis, robotics and probability theory.

In image analysis, left-invariant diffusion PDEs on $SE(3)$
have been widely used for crossing-preserving diffusion and enhancement of fibers in diffusion-weighted MRI images of brain white matter \cite{portegies_improving_2015,momayyez-siahkal_3d_2009,Skibbe,MEESTERS2017,reisert_fiber_2011,prckovska_extrapolating_2010}, or for crossing-preserving enhancements of 3D vasculature in medical images \cite{Janssen2018}.
They extend classical works on multi-scale image representations \cite{Iiji59a,Koenderink,ter_haar_romeny_front-end_2003,weickert_anisotropic_1998} to Lie groups \cite{duits_scale_2007}.

In robotics, they play a role via the central limit theorem \cite{Benoist} in work-space generation of robot arms \cite[ch.12]{chirikjian_engineering_2000} and they appear indirectly in Kalman-filtering on $SE(3)$ for tracking  \cite{Barbaresco} and camera motion estimation \cite{Berger}.

In probability theory, diffusion systems on Lie groups describe Brownian motions \cite{Oksendal,Hsu}
and they appear as limits in central limit theorem on Lie groups \cite{Benoist}.

Both in probability theory \cite{Feller} and in image analysis \cite{Fels2003,Duits2003a,duits_axioms_2004,Pedersen},
the spectral decomposition of the evolution operator also allows to simultaneously deal with important variants of the diffusion evolution. These variants of the heat-evolution are obtained by taking fractional powers $-(-\Delta)^{\alpha}$, cf.~\!\cite{yosida_functional_1980}, of the minus Laplacian operator $\Delta=\textrm{div} \circ \textrm{grad}$ that
generates the heat flow (due to Fick's law and the Gauss divergence theorem), where $\alpha \in (0,1]$.

This generalization allows for heavy tailed distributions of $\alpha$-stable L\'{e}vy processes, that arise in a fundamental generalization \cite{Feller} of the central limit theorem \emph{where one drops the finite variance condition}. Here we note that recently an extension of the central limit on linear groups (such as $SE(3)$) has been achieved for finite 2nd order moments \cite{Benoist}.
In engineering applications where (iterative group-) convolutions are applied \cite[ch.12\&13]{chirikjian_engineering_2000}, \cite{BekkersMICCAI,Wellink,winkels20183d,worrall2018cubenet,weiler20183d,CittiX,Mallat,oyallon2013generic}
the `kernel width' represents the spread of information or the scale of observing the signal.
In case the applications allow for an underlying probabilistic model with finite variances, variance
is indeed a good measure for `kernel width'. But often this is not the case. Probability kernels for stochastic L\'{e}vy processes (used in directional statistics \cite{MardiaJuppBook1999}, stock market modeling \cite{Wu},
natural image statistics \cite{Pedersen}), modeling of point-spread functions in acquired images (e.g. in spectroscopy \cite{BelkicandBelkic2010CRCPress}),
do require distributions with heavier tails than diffusion kernels. Therefore, `full width at half maximum' is a more generally applicable measure for kernel width than variance, as it applies
to all $\alpha$-stable L\'{e}vy processes. The probability distributions for $\alpha<1$ encode a longer range of interaction via their heavy tails and still allow for unlimitedly sharp kernels. 
\begin{figure}
\centerline{
\includegraphics[width=\iftoggle{insertremarks}{0.76}{0.76}\hsize]{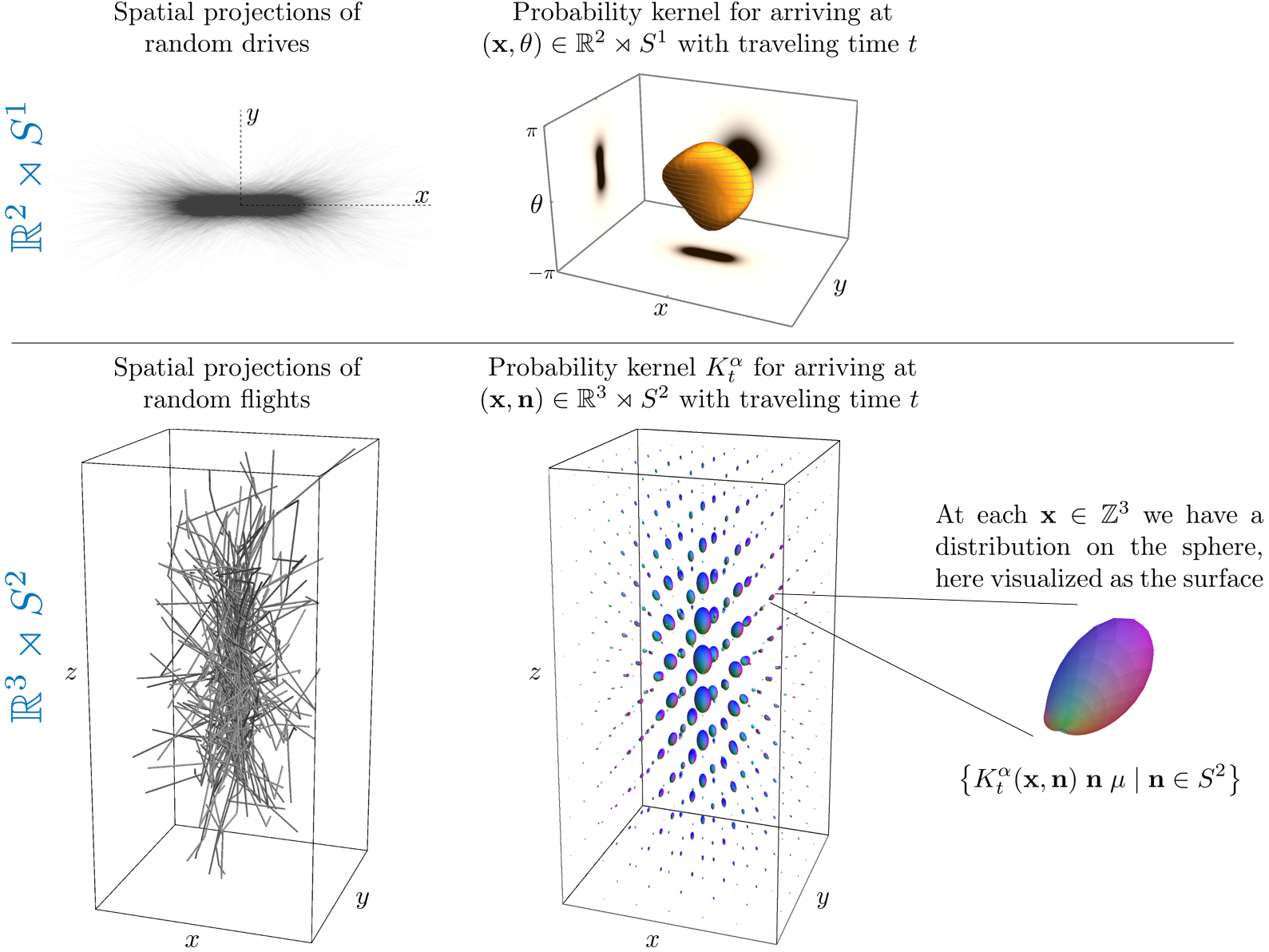}
}
\caption{Various visualization of the diffusion process ($\alpha=1$) on $\mathbb{R}^{d} \rtimes S^{d-1}$, for $d=2$ and $d=3$. Top: random walks (or rather `drunk man's drives') and an iso-contour of the limiting diffusion kernel, for the case $d=2$ studied in previous works (see e.g. \cite{zhang_numerical_2016,duits_left-invariant_2010-2,citti_cortical_2006-1}).
Bottom: random walks (or rather `drunk man's flights') and a visualization of the limiting
distribution for the case $d=3$. This limiting distribution is a hypo-elliptic diffusion kernel $(\mathbf{x},\mathbf{n}) \mapsto K^{\alpha=1}_{t}(\mathbf{x},\mathbf{n})$ that we study in this article.
We visualize the kernel $K^{\alpha=1}_{t}$ by a spatial grid of surfaces, where all surfaces are scaled by the same $\mu>0$.
\label{fig:sample-paths}}
\end{figure}
\subsection{Structure of the Article }

The structure of the article is as follows. In the first part of the introduction we briefly discussed the history of the Fourier transform, and its generalization to other groups that are the semi-direct product of the translation group and another matrix group, where we provided an overview of related works. Then we specified our domain of interest --- the Fourier transform on the homogeneous space $G/H$ of positions and orientations which is a Lie group quotient of the Lie group $G=SE(3)$ with a subgroup $H$ isomorphic to $\{\mathbf{0}\} \times SO(2)$.
Then we addressed its application of solving PDE systems on $G/H$, motivated from applications in image analysis, robotics and probability theory.

There are 4 remaining subsections of the introduction. In Subsection~\ref{ch:intro-1} we provide basic facts on the homogeneous space $G/H$ of positions and orientations and we provide preliminaries for introducing a Fourier transform on $G/H$. In Subsection~\ref{ch:intro-2} we formulate the PDEs of interest on $G/H$ that we will solve. In Subsection~\ref{ch:intro-3} we formulate the corresponding PDEs on the group $G$. In Subsection~\ref{rem:3} we provide a roadmap on the spectral decomposition of the PDE evolutions.

In Section~\ref{ch:symmetry}, based on previous works, we collect the necessary prior information about the PDEs of interest and the corresponding kernels. We also describe how to extend the
 case $\alpha=1$ (the diffusion case) to the general 
 case $\alpha \in (0,1]$.

In Section~\ref{ch:FTSE3}, we describe the Fourier transform on the Lie group $SE(3)$, where we rely on UIRs of $SE(3)$. In particular, by relating the UIRs to the dual orbits of $SO(3)$ and by using a decomposition w.r.t. an orthonormal basis of modified spherical harmonics, we recall an explicit formula for the inverse Fourier transform. 

In Section~\ref{ch:FTquot}, we present a Fourier transform $\mathcal{F}_{G/H}$ on the quotient $G/H = \quot$. Our construction requires an additional constraint --- an input function must be bi-invariant w.r.t. subgroup $H$, as explained in Remark 2. This extra symmetry constraint is satisfied by the PDE kernels of interest. We prove a theorem, where we present: 1) a matrix representation for the Fourier transform on the quotient, 2) an explicit inversion formula and 3) a Plancherel formula.

In Section~\ref{ch:FTPDE}, we apply our Fourier transform on the quotient to solve the PDEs of interest. The solution is given by convolution of the initial condition with the specific kernels (which are the probability kernels of $\alpha$-stable L\'{e}vy process). We find the exact formulas for the kernels in the frequency domain relying on a spectral decomposition of the evolution operator (involved in the PDEs). We show that this result can be obtained either via conjugation of the evolution operator with our Fourier transform on $\mathbb{R}^{3}\rtimes S^2$ or (less efficiently) via conjugation of the evolution operator with the Fourier transform acting only on the spatial part $\R^3$. Then we present a numerical scheme to approximate the kernels via Monte-Carlo simulation and we provide a comparison of the exact solutions and their approximations. Finally, in Section~\ref{ch:coclusion}, we summarize our results and discuss their applications. In the appendices, we address the probability theory and stochastic differential equations (SDEs) regarding L\'{e}vy processes on $\mathbb{R}^{3} \rtimes S^2$. 
\\[6pt]
The main contributions of this article are:
\begin{itemize}
\item construction of $\mathcal{F}_{\quot}$ --- the Fourier transform on the quotient $\quot$, in formula~(\ref{FThom});
\item matrix representation for $\mathcal{F}_{\quot}$, explicit inversion and Plancherel formulas, in Theorem~\ref{corr:1};
\item explicit spectral decompositions of PDE evolutions for $\alpha$-stable L\'{e}vy process on $\quot$, both in the Fourier domains of $\quot$ and $\R^3$, in Theorem~\ref{th:decomposition}.
The new spectral decomposition in the Fourier domain of $\quot$ is simpler and involves ordinary spherical harmonics.
\item the exact formulas for the probability kernels of $\alpha$-stable L\'{e}vy processes on $\quot$, in Theorem~\ref{th:three}. This also includes new formulas for the heat kernels (the case $\alpha=1$), that are more efficient than the heat kernels presented in previous work \cite{portegies_new_2017}.
\item simple formulation and verifications (Monte-Carlo simulations) of discrete random walks for $\alpha$-stable L\'{e}vy processes on $\quot$ in Proposition~\ref{prop:RW}. The corresponding SDEs are in Appendix~\ref{ch:PT}.
\end{itemize}
\subsection{Introduction to the Fourier Transform on the Homogeneous Space of Positions and Orientations}\label{ch:intro-1}

Let $G=SE(3)$ denote the Lie group of rigid body motions, equipped with group product:
\begin{equation} \label{product}
g_1 g_2 =(\mathbf{x}_1, \mathbf{R}_{1})(\mathbf{x}_{2},\mathbf{R}_{2})= (\mathbf{R}_{1}\mathbf{x}_{2}+ \mathbf{x}_{1}, \mathbf{R}_{1} \mathbf{R}_{2}), \quad \textrm{ with } g_k=(\mathbf{x}_k,\mathbf{R}_{k})\in G, \; k=1,2.
\end{equation}
Here $\mathbf{x}_k \in \RR^3$ and $\mathbf{R}_{k} \in SO(3)$.
Note that $SE(3)\!=\!\R^{3}\rtimes\!SO(3)$ is a semi-direct product of $\R^3$ and $SO(3)$.

\begin{Definition}
Let $B(\mathcal{H})$ denote the vector space of bounded linear operators
on some Hilbert space $\mathcal{H}$. \\
Within the space $B(\mathcal{H})$ we
denote the subspace of bounded linear trace-class operators by
\[
B_2(\mathcal{H})= \left\{A: \mathcal{H} \rightarrow \mathcal{H} \; | \; A \text{ linear and }
|\!|\!|A|\!|\!|^2:=\operatorname{trace}(A^*A)< \infty \right\}.
\]
\end{Definition}
\begin{Definition}
Consider a mapping
$\sigma: G \to B(\mathcal{H}_{\sigma})$, where $\mathcal{H}_{\sigma}$ denotes the Hilbert space on which each $\sigma_{g}$ acts. Then $\sigma$ is a Unitary Irreducible Representation (UIR) of $G$ if
\begin{enumerate}
\item
 $\sigma: G \to B(\mathcal{H}_{\sigma})$ is a homomorphism,
 \item
$\sigma_g^{-1}=\sigma_{g}^*$ for all $g\in G$,
\item
there does not exist a closed subspace $V$ of $\mathcal{H}_{\sigma}$ other than $\{0,\mathcal{H}_{\sigma}\}$ such that $\sigma_g V \subset V$.
\end{enumerate}
\end{Definition}
We denote by $\hG$ the dual group of $G$. Its elements are equivalence classes of UIR's, where one identifies elements via
$\sigma_1 \sim \sigma_2 \Leftrightarrow \textrm{ there exists a unitary linear operator $\upsilon$, s.t. } \sigma_1 = \upsilon \circ \sigma_2 \circ \upsilon^{-1}.$
\begin{Definition}
The Fourier transform $\mathcal{F}_{G}(f)=((\mathcal{F}_{G}f)(\sigma))_{\sigma \in \hat{G}}$ of a square-integrable, measurable and bounded function $f$ on $G$ is a measurable field of bounded operators indexed by unitary irreducible representations (UIR's) $\sigma$.
Now $\hat{G}$ can be equipped with a canonical Plancherel measure $\nu$ and the Fourier transform $\mathcal{F}_{G}$ admits\footnote{Lie group $SE(3)$ is unimodular since the left-invariant Haar measure and the right-invariant Haar measures coincide.
Furthermore, it is of type I (i.e. dual group and quasi dual group coincide) and thereby admits a Plancherel theorem \cite{fuhr_abstract_2005,folland_course_1994}.} an extension unitary operator from $\mathbb{L}_{2}(G)$ to the direct-integral space $\int_{\hat{G}}^{\oplus} B_{2}(\mathcal{H}_{\sigma}) {\rm d}\nu(\sigma)$.
This unitary extension \cite[4.25]{folland_course_1994}
 (also known as  `Plancherel transform' \cite[thm.3.3.1]{fuhr_abstract_2005})
is given by
 \begin{equation} \label{FT}
 \begin{array}{l}
\mathcal{F}_{G}(f) = \int \limits^{\oplus}_{\hat{G}} \hat{f}(\sigma)\; {\rm d}\nu(\sigma), \; \textrm{with }\\
\hat{f}(\sigma) = \left(\mathcal{F}_{G}f\right)(\sigma) = \int \limits_G f(g)\; \sigma_{g^{-1}} \rmd g \; \in B_2(\mathcal{H}_\sigma), \; \textrm{ for all } \sigma \in \hat{G},
\end{array}
 \end{equation}
for all $f \in \mathbb{L}_1(G) \cap \mathbb{L}_2(G)$.
\end{Definition}
\noindent
The Plancherel theorem states
that $\|\mathcal{F}_{G}(f)\|^{2}_{\mathbb{L}_{2}(\hat{G})}= \int_{\hat{G}}|\!|\!|\mathcal{F}_{G}(f)(\sigma)|\!|\!|^2 {\rm d}\nu(\sigma) =\int_G |f(g)|^2 {\rm d}g=\|f\|^2_{\mathbb{L}_{2}(G)}$ for all $f \in \mathbb{L}_{2}(G)$,
and we have the inversion formula $f = \mathcal{F}_{G}^{-1} \mathcal{F}_{G}f= \mathcal{F}_{G}^* \mathcal{F}_{G}f$.
For details see \cite{folland_course_1994,fuhr_abstract_2005} and for detailed explicit computations see \cite{chirikjian_engineering_2000}.

In this article, we will constrain and modify the Fourier transform $\mathcal{F}_{G}$ on $G=SE(3)$ such that we obtain a suitable Fourier transform $\mathcal{F}_{G/H}$ defined on a homogeneous space\footnote{Although the semi-direct product notation $\R^{3} \rtimes S^2$ is formally not correct as $S^2$ is not a Lie group, it is convenient:
It reminds that $G/H$ denotes the homogeneous space of positions and orientations.}
\begin{equation} \label{quot}
\R^{3} \rtimes S^2 :=
G/H \textrm{ with subgroup }H=\{\mathbf{0}\} \times \operatorname{Stab}_{SO(3)}(\mathbf{a})
\end{equation}
of left cosets, where
$\operatorname{Stab}_{SO(3)}(\mathbf{a})=
\{\mathbf{R} \in SO(3)\;|\; \mathbf{R}\mathbf{a}=\mathbf{a}\}$ denotes the subgroup of $SO(3)$ that stabilizes an a priori reference axis $\mathbf{a} \in S^{2}$, say $\mathbf{a}=\mathbf{e}_z=(0,0,1)^T$. In the remainder of this article we set this choice $\ba = \be_z$.
\begin{Remark}\label{rem:terminology} (notation and terminology) \\
Elements in (\ref{quot}) denote equivalence classes of rigid body motions $g=(\mathbf{x}, \mathbf{R}_{\mathbf{n}}) \in SE(3)$ that map $(\mathbf{0},\ba)$ to $(\mathbf{x},\mathbf{n})$:
$$[g] =: (\bx,\bn) \in \quot \quad \Leftrightarrow \quad g \odot (\mathbf{0},\ba) = (\bx,\bn),$$
under the (transitive) action
\begin{equation} \label{actiononset}
g \odot (\mathbf{x}',\mathbf{n}')=
(\mathbf{R}\mathbf{x}' +\mathbf{x}, \mathbf{R} \mathbf{n}'), \qquad \textrm{ for all }g=(\mathbf{x},\mathbf{R}) \in SE(3), \; (\mathbf{x}',\mathbf{n}') \in \R^{3}\rtimes S^2.
\end{equation}
Therefore, we simply denote the equivalence classes $[g]$ by $(\mathbf{x},\mathbf{n})$. This is similar to the conventional writing $\mathbf{n} \in S^{2}=SO(3)/SO(2)$. Throughout this manuscript we refer to $G/H$ as `the homogeneous space of positions and orientations' and henceforth $\bR_{\bn}$ denotes \underline{any} rotation that maps the reference axis $\ba$ into $\bn$.
\end{Remark}
The precise definition of the Fourier transform $\mathcal{F}_{G/H}$ on the homogeneous space $G/H$ will be presented in Section~\ref{ch:FTquot}. It relies on the decomposition into unitary irreducible representations in (\ref{FT}), but we must take both a domain and a range restriction into account. This will be explained later in Section~\ref{ch:FTquot}. Next we address an a priori domain constraint that is rather convenient than necessary.
\begin{Remark}
We will
constrain the Fourier transform $\mathcal{F}_{G/H}$ to
\begin{equation} \label{extraconstraint}
\mathbb{L}_{2}^{sym}(G/H) := \left\{f \in \mathbb{L}_{2}(G/H) \;|\; \forall_{\mathbf{R} \in \operatorname{Stab}_{SO(3)}(\mathbf{a})} :\, f(\mathbf{x},\mathbf{n})=
f(\mathbf{R}\mathbf{x},\mathbf{R}\mathbf{n})
\right\}.
\end{equation}
\end{Remark}
This constraint is convenient in view of the PDEs of interest (and the symmetries of their kernels) that we will formulate in the next subsection, and that will solve via Fourier's method in Section~\ref{ch:FTPDE}.
\subsection{Introduction to the PDEs of Interest on the Quotient $\quot$}\label{ch:intro-2}
\noindent
Our main objective is to use the Fourier transform $\mathcal{F}_{G/H}$ to solve the following PDEs on $\R^3 \rtimes S^{2}$:
\begin{equation}\label{eq:differentialQ}
\boxed{\left\{\begin{aligned}
\frac{\partial}{\partial t} W_{\alpha}(\mathbf{x},\mathbf{n},t) &=
Q_{\alpha} W_{\alpha}(\mathbf{x},\mathbf{n},t),  \\
W_{\alpha}(\mathbf{x},\mathbf{n},0) &= U(\mathbf{x},\mathbf{n}),
\end{aligned} \right. }
\end{equation}
where  $(\mathbf{x},\mathbf{n}) \in \quot$, $t \geq 0$, $\alpha \in (0,1]$ and the generator \begin{equation} \label{Qalphafund}
\boxed{Q_{\alpha}:=-(-Q)^{\alpha}}
\end{equation}
is expressed via
\[\boxed{
Q=D_{11} \|\mathbf{n} \times \nabla_{\R^3}\|^2 + D_{33}(\mathbf{n} \cdot \nabla_{\R^3})^2 + D_{44}\Delta^{S^2}_{\bn},
}
\]
with $D_{33}>D_{11}\geq 0, D_{44}>0$, and with $\Delta^{S^2}_{\bn}$ the Laplace-Beltrami operator on $S^2 \!=\! \left\{\bn \in \R^3 \,\big|\, \|\bn\| = 1\right\}$.

Note, that the generator $Q$ is a self-adjoint unbounded operator with domain
\[
\mathcal{D}(Q):=\mathbb{H}_{2}(\R^3) \otimes \mathbb{H}_{2}(S^2),
\]
where $\mathbb{H}_2$ denotes the Sobolev space $\mathbb{W}_2^2$.

The semigroup for $\alpha=1$ is a strongly continuous semigroup on $\mathbb{L}_{2}(\R^3 \rtimes S^2)$ with a closed generator, and by taking the fractional power of the generator one obtains another strongly continuous semigroup, as defined and explained in a more general setting in the work by Yosida \cite[ch:11]{yosida_functional_1980}. The fractional power is formally defined by
\begin{equation} \label{def}
Q_{\alpha}W=-(-Q)^{\alpha} W:= \frac{\sin \alpha \pi}{\pi} \int_{0}^{\infty} \lambda^{\alpha-1} (Q-\lambda I)^{-1}(-Q W) {\rm d}\lambda \textrm{ for all }W \in \mathcal{D}(Q).
\end{equation}
\iftoggle{insertremarks}{
One has $(-Q)^{\alpha} \circ (-Q)^{\beta}W= (-Q)^{\alpha+\beta}W$ for all $W \in \mathcal{D}(Q^{\alpha+\beta})$ and all $\alpha,\beta>0$ with $\alpha+\beta <1$ and
\begin{equation} \label{domain}
\mathcal{D}(Q_{\alpha}):=\mathbb{H}_{2\alpha}(\R^3) \otimes \mathbb{H}_{2\alpha}(S^2) \supset \mathcal{D}(Q),
\end{equation}
where $f \in \mathbb{H}_{2\alpha}(\R^3)$ implies $\int\limits_{\R^3} |\hat{f}(\omega)|^2 \left(1 +|\omega|^{2 \alpha}\right) {\rm d} \omega < \infty$, and $f\!\in\!\mathbb{H}_{2\alpha}(S^2)$ implies $\left|\Delta^{S^2}_{\bn}\right|^{\alpha} f\!\in\! \mathbb{L}_2(S^2)$.}{}

In Subsection~\ref{rem:3} below, we will show that the common technical representation~(\ref{def}) is not really needed for our setting.
In fact, it will be very easy to account for $\alpha \in (0,1]$ in the solutions; by a spectral decomposition we only need to take fractional powers of certain eigenvalues in the Fourier domain. For the moment the reader may focus on the case $\alpha=1$, where system (\ref{eq:differentialQ}) becomes an ordinary elliptic diffusion system which is hypo-elliptic (in the sense of H\"{o}rmander \cite{hormander_hypoelliptic_1967}) if $D_{11}=0$.

The PDEs~\ref{eq:differentialQ} have our interest as they are Forward-Kolmogorov equations for $\alpha$-stable L\'{e}vy processes on $G/H$.
See Appendix~\ref{ch:PT} for a precise formulation of discrete and continuous stochastic processes. This generalizes previous works on such basic processes \cite{Feller,duits_axioms_2004} with applications in financial mathematics \cite{Misiorek2012} and computer vision \cite{felsberg_monogenic_2004,Pedersen,Kanters2007,schmidt_morphological_2016}, from Lie group $\R^3$ to the Lie group quotient $\R^{3} \rtimes S^{2}$.

See Figure~\ref{fig:sample-paths} for a visualization of sample paths from the discrete stochastic processes explained in Appendix~\ref{ch:PT}. Roughly speaking, they represent `drunk man's flights' rather than `drunk man's walks'.

\subsection{Reformulation of the PDE on the Lie Group $SE(3)$}\label{ch:intro-3}

Now we reformulate and extend our PDEs~(\ref{eq:differentialQ}) to the Lie group $G=SE(3)$ of rigid body motions, equipped with group product (\ref{product}).
This will help us to better recognize symmetries, as we will show in Subsection~\ref{sec:Sym}.
To this end, the PDEs are best expressed in a basis of left-invariant vector fields
$\{g \mapsto  \left.\mathcal{A}_{i}\right|_{g}\}_{i=1}^6$ on $G$. Such left-invariant vector fields are obtained by push forward from the left-multiplication $L_{g_1} g_2 := g_1 g_2$ as
\[
\left.\mathcal{A}_{i}\right|_{g}= (L_{g})_{*} A_i \in T_{g}(G),
\]
where $A_i:=\left.\mathcal{A}_{i}\right|_{e}$ form an orthonormal basis for the Lie algebra $T_{e}(G)$. We choose such a basis typically such that the first 3 are spatial generators $A_{1}=\partial_x, A_{2}=\partial_y, A_{3}=\partial_{z}=\mathbf{a} \cdot \nabla_{\R^3}$ and the remaining 3 are rotation generators, in such a way that $A_{6}$ is the generator of a counter-clockwise rotation around the reference axis $\mathbf{a}$. They are infinitesimal generators w.r.t. the right regular representation $g \mapsto \mathcal{R}_{g}$, with $\mathcal{R}_{g_1}\tilde{U}(g_2)=\tilde{U}(g_2 g_1)$ for all $g_1,g_2 \in G$.
This means that for all $\tilde{U} \in C^{1}(G)$ and $g \in G$ one has
\begin{equation} \label{li}
\mathcal{A}_{i}\tilde{U}(g)= ({\rm d}\mathcal{R}(A_i)(\tilde{U}))(g):=\lim \limits_{t \downarrow 0} \frac{\left(\;(\mathcal{R}_{e^{t A_i}}-I)\tilde{U}\right)(g)}{t} =
\lim \limits_{t \downarrow 0} \frac{\tilde{U}(g\, e^{t A_i})-\tilde{U}(g)}{t},
\end{equation}
where $A \mapsto e^A$ denotes the exponent that maps Lie algebra element
$A \in T_{e}(G)$ to the corresponding Lie group element.
The explicit formulas for the left-invariant vector fields in Euler-angles (requiring 2 charts)
can be found in Appendix~\ref{ch:LI}, or in \cite{chirikjian_engineering_2000,Duits2011}.

Now we can re-express the PDEs (\ref{eq:differentialQ}) on the group $G=SE(3)$ as follows:
\begin{equation} \label{PDEgroup}
\boxed{\left\{\begin{aligned}
\frac{\partial}{\partial t} \tilde{W}_{\alpha}(g,t) &=
\tilde{Q}_{\alpha} \tilde{W}_{\alpha}(g,t), & g \in G, t\geq 0 \\
\tilde{W}_{\alpha}(g,0) &= \tilde{U}(g), & g \in G,
\end{aligned} \right. }
\end{equation}
where the generator 
\begin{equation}\label{eq:Qtildealpha}
\boxed{
\tilde{Q}_{\alpha}:=-(-\tilde{Q})^{\alpha}
}
\end{equation}
is again a fractional power ($\alpha \in (0,1]$)
of the diffusion generator $\tilde{Q}$ given by
\begin{equation}
\boxed{
\tilde{Q}=D_{11} (\mathcal{A}_{1}^2+ \mathcal{A}_{2}^{2}) + D_{33}\,\mathcal{A}_{3}^2 + D_{44}(\mathcal{A}_{4}^2+ \mathcal{A}_{5}^{2}),
}
\end{equation}
where $\mathcal{A}_{i}^2=\mathcal{A}_{i} \circ \mathcal{A}_{i}$ for all $i\in \{1,\ldots,5\}$.
The initial condition in (\ref{PDEgroup}) is given by
\[
\tilde{U}(g)= \tilde{U}(\mathbf{x},\mathbf{R})= U(\mathbf{x},\mathbf{R}\mathbf{a}).
\]
Similar to the previous works \cite{portegies_new_2017,duits_morphological_2012}
one has
\begin{equation}\label{eq:16b2}
\begin{array}{l}
\tilde{W}_{\alpha}(\mathbf{x},\mathbf{R},t)=W_{\alpha}(\mathbf{x},\mathbf{R}\mathbf{a},t),
\end{array}
\end{equation}
that holds for all $t \geq 0, \; (\mathbf{x},\mathbf{R}) \in SE(3)$. 
\begin{Remark}
The above formula~(\ref{eq:16b2}) relates the earlier PDE formulation (\ref{eq:differentialQ}) on the quotient $G/H$ to the PDE formulation (\ref{PDEgroup}) on the group $G$.
It holds since we have the relations
\[
\begin{array}{l}
\mathcal{A}_{6}\tilde{W}_{\alpha}(\mathbf{x},\mathbf{R},t)=0,\\
(\mathcal{A}_{5}^2 +\mathcal{A}_{4}^2)\tilde{W}_{\alpha}(\mathbf{x},\mathbf{R},t)=\Delta^{S^2}
W_{\alpha}(\mathbf{x},\mathbf{R}\mathbf{a},t), \\[6pt]
\mathcal{A}_{3}\tilde{W}_{\alpha}(\mathbf{x},\mathbf{R}_{\mathbf{n}},t)= \mathbf{n} \cdot \nabla_{\R^3} W_{\alpha}(\mathbf{x},\mathbf{n},t), \\
\left(\mathcal{A}_{1}^{2}+\mathcal{A}_{2}^2\right) \tilde{W}_{\alpha}(\mathbf{x},\mathbf{R},t)= \left(\Delta^{\mathbb{R}^{3}}- \mathcal{A}_{3}^2\right)\tilde{W}_{\alpha}(\mathbf{x},\mathbf{R},t)= \|\mathbf{n} \times \nabla_{\R^3}\|^2 \;W_{\alpha}(\mathbf{x},\mathbf{R}\mathbf{a},t)
\end{array}
\]
so that the generator of the PDE (\ref{PDEgroup}) on $G$ and the generator of the PDE (\ref{eq:differentialQ}) on $G/H$ indeed stay related via
\[
\tilde{Q}_{\alpha}\tilde{W}_{\alpha}(\mathbf{x},\mathbf{R},t)=Q_{\alpha}W_{\alpha}(\mathbf{x},\mathbf{R}\mathbf{a},t) \textrm{ for all }t\geq 0.
\]
\end{Remark}
\iftoggle{insertremarks}{
\begin{Remark}
For $\alpha=\frac{1}{2}$, one has
$\frac{\partial^2}{\partial t^2} +
\tilde{Q}= \left( \frac{\partial}{\partial t} - \sqrt{-\tilde{Q}}\right) \left( \frac{\partial}{\partial t} + \sqrt{-\tilde{Q}}\right)$ and by the Fourier transform on $G$ it then follows that
the PDE system (\ref{PDEgroup}) can be replaced by the following Poisson system on $G \times \mathbb{R}^+$ :
\[
\left\{
\begin{array}{rcll}
\left(\frac{\partial^2}{\partial t^2} +
\tilde{Q}\right) \tilde{W}_{\frac{1}{2}}(g,t) &=&0 &\  g \in G, t\geq 0, \ \textrm{ with } \forall_{t\geq 0}:\tilde{W}_{\frac{1}{2}}(\cdot,t) \in \mathbb{L}_{2}(G)\\
\tilde{W}_{\frac{1}{2}}(g,0) &=& \tilde{U}(g), & \ g \in G,
\end{array} \right.
\]
that no longer involves a fractional power.
\end{Remark}
}{
}
\subsection{\label{rem:3} A Preview on the Spectral Decomposition of the PDE Evolution Operator and the Inclusion of $\alpha$}

Let \iftoggle{insertremarks}{$U \in \mathcal{D}(Q_{\alpha})$}{$U$} be in the domain of the generator $Q_{\alpha}$ given by (\ref{Qalphafund}), of our evolution (\ref{eq:differentialQ}). \iftoggle{insertremarks}{}{For a formal definition of this domain we refer to~\cite[Eq. 9]{arxivFT}.} Let its spatial Fourier transform be given by
\begin{equation}\label{eq:spatFT0}
\overline{U}(\bomega,\mathbf{n})=\left[\mathcal{F}_{\R^{3}} U(\cdot, \mathbf{n})\right](\bomega):= \frac{1}{(2\pi)^{\frac{3}{2}}}\int \limits_{\R^3} U(\mathbf{x}, \bn) \; e^{- i
\bomega \cdot \bx} \;{\rm d}\bx.
\end{equation}
To the operator $Q_{\alpha}$ we associate the corresponding operator $-(-\mathcal{B})^{\alpha}$ in the spatial Fourier domain by
\begin{equation} \label{correspondence}
-(-\mathcal{B})^{\alpha}=\left(\mathcal{F}_{\R^3} \otimes 1_{\mathbb{L}_{2}(S^2)}\right) \circ Q_{\alpha} \circ \left(\mathcal{F}_{\R^3}^{-1} \otimes 1_{\mathbb{H}_{2\alpha}(S^2)}\right).
\end{equation}
Then direct computations show us that
\begin{equation} \label{BB}
-(-\mathcal{B})^{\alpha} \overline{U}(\boldsymbol{\omega},\mathbf{n})=
\left[-(-\mathcal{B}_{\boldsymbol{\omega}})^{\alpha} \overline{U}(\boldsymbol{\omega},\cdot)\right](\mathbf{n}), \textrm{ for all } \mathbf{n} \in S^2,
\end{equation}
where for each fixed $\boldsymbol{\omega} \in \R^3$, the operator $-(-\mathcal{B}_{\boldsymbol{\omega}})^{\alpha}: \mathbb{H}_{2\alpha}(S^2) \to \mathbb{L}_{2}(S^2)$ is given by
\begin{equation} \label{Bexp}
-(-\mathcal{B}_{\boldsymbol{\omega}})^{\alpha}=- \left( -D_{44}\Delta^{S^2}_{\bn} + D_{11} \|\boldsymbol{\omega} \times \mathbf{n}\|^2 + D_{33} (\boldsymbol{\omega}\cdot\mathbf{n})^2 \right)^{\alpha}.
\end{equation}
In this article, we shall employ Fourier transform techniques 
to derive a complete orthonormal basis (ONB) of eigenfunctions
\begin{equation} \label{ef}
\left\{\Phi_{\boldsymbol{\omega}}^{l,m} \; |\; l \in \mathbb{N}_0, m \in \mathbb{Z} \textrm{ with }|m| \leq l\right\},
\end{equation}
in $\mathbb{L}_{2}(S^2)$
for the operator $-(-\mathcal{B}_{\boldsymbol{\omega}}):=-(-\mathcal{B}_{\boldsymbol{\omega}})^{\alpha=1}$. Then clearly, this basis is also
an ONB of eigenfunctions for $-(-\mathcal{B}_{\boldsymbol{\omega}})^{\alpha}$, as we only need to take the fractional power of the eigenvalues. Indeed once eigenfunctions (\ref{ef}) and the eigenvalues
\begin{equation} \label{efunc2}
\mathcal{B}_{\boldsymbol{\omega}}\Phi_{\boldsymbol{\omega}}^{l,m} = \lambda^{l,m}_{r} \, \Phi_{\boldsymbol{\omega}}^{l,m}, \textrm{ with }r=\|\boldsymbol{\omega}\|,
\end{equation}
are known, the exact solution of (\ref{eq:differentialQ}) is given by (shift-twist) convolution on $\R^3 \rtimes S^2$
as defined below
\begin{equation} \label{decompose}
\begin{array}{l}
W_{\alpha}(\bx,\bn,t) = (K_t^{\alpha} \ast U)(\bx,\bn) := \int \limits_{S^2} \int \limits_{\R^{3}}
K_{t}^{\alpha}(\bR_{\bn'}^T(\bx-\bx'),\bR_{\bn'}^T\bn)\, U(\bx',\bn') \; {\rm d}\bx'{\rm d}\mu_{S^2}(\bn')
\\
=\int \limits_{\mathbb{R}^3} \sumlm \, \left\langle \; \overline{U}(\bomega,\cdot) \, , \, \Philmom (\cdot)\; \right\rangle_{\mathbb{L}_2(S^2)} \Philmom(\bn) \; e^{-(-\lambdalmr)^{\alpha} t} e^{i \bx \cdot \bomega} \rmd \bomega,
\\[9pt]
\textrm{with the probability kernel
given by } \\[9pt]
\qquad
K_{t}^{\alpha}(\bx,\bn)= \left[\mathcal{F}^{-1}_{\R^3}\left(\overline{K}_{t}^{\alpha}(\cdot,\bn)\right)\right](\bx), \\
\qquad \; \; \textrm{ with } \overline{K}_t^{\alpha}(\bomega,\mathbf{n})=\sumlm \overline{\Phi^{l,m}_{\bomega}(\mathbf{a})}\,\Phi^{l,m}_{\bomega}(\mathbf{n})\,e^{-(-\lambdalmr)^{\alpha} t}.
\end{array}
\end{equation}
Here the inner product in $\mathbb{L}_{2}(S^2)$ is given by
\begin{equation}\label{eq:inerprodL2S2}
\left\langle y_1(\cdot) , y_2(\cdot)\right\rangle_{\mathbb{L}_{2}(S^2)} := \int \limits_{S^2} y_1(\mathbf{n})\, \overline{y_2(\mathbf{n})}\, {\rm d}\mu_{S^2}(\mathbf{n}).
\end{equation}
where $\mu_{S^2}$ is the usual Lebesgue measure on the sphere $S^2$.
\begin{Remark}
The eigenvalues $ \lambda^{l,m}_{r}$ only depend on $r=\|\bomega\|$ due to the symmetry $\Phi_{\bR \bomega}^{l,m}(\bR \bn)=\Phi_{\bomega}^{l,m}(\bn)$ that one directly recognizes from (\ref{Bexp}) and (\ref{efunc2}).
\end{Remark}
\begin{Remark}
The kernels $K_{t}^{\alpha}$ are the probability density kernels of stable L\'{e}vy processes on $\R^{3}\rtimes S^2$, see Appendix~\ref{app:A1}. Therefore, akin to the $\R^n$-case \cite{Feller,Pedersen}, we refer to them as the $\alpha$-stable L\'{e}vy kernels on $\mathbb{R}^{3} \rtimes S^2$.
\end{Remark}
\section{Symmetries of the PDEs of Interest \label{ch:symmetry}}
Next we employ PDE-formulation~(\ref{PDEgroup}) on the group $G = SE(3)$ to
summarize the symmetries for the probability kernels $K_{t}^{\alpha}: \R^{3} \rtimes S^2 \to \R^+$. For
details see \cite{portegies_new_2015} and \cite{portegies_new_2017}.
\subsection{PDE symmetries}\label{sec:Sym}

Consider PDE system~(\ref{PDEgroup}) on the group $G = SE(3)$. Due to left-invariance (or rather left-covariance) of the PDE, linearity of the map $\tilde{U}(\cdot) \mapsto \tilde{W}_{\alpha}(\cdot,t)$, and the Dunford-Pettis theorem~\cite{Bukhvalov}, the solutions are obtained by group convolution with a kernel
$\tilde{K}_{t}^{\alpha} \in \mathbb{L}_{1}(G)$:
\begin{equation} \label{solconv}
\tilde{W}_{\alpha}(g,t)= \left(\tilde{K}_{t}^{\alpha} * \tilde{U}\right)(g):= \int \limits_{G} \tilde{K}_{t}^{\alpha}(h^{-1}g)\, \tilde{U}(h)\; {\rm d} h,
\end{equation}
where we took the convention that the probability kernel acts from the left\iftoggle{insertremarks}{
\footnote{We take the convention $\left(f_1 * f_2\right)(g) = \int\limits_G f_1(h^{-1} g) f_2(h) {\rm d} h$ since we want the kernel to act from the left and since we want $\mathcal{F}_{G}\left(f_1 * f_2\right) =\mathcal{F}_{G}\left(f_1\right) \circ \mathcal{F}_{G}\left(f_2\right)$. Often one also encounters the convention where $f_1$ and $f_2$ are switched, e.g.~\cite{chirikjian_engineering_2000,MIngLiao,fuhr_abstract_2005}.}}{}.
In the special case $U=\delta_e$ with unity element $e=(\mathbf{0},\mathbf{I})$ we get
$\tilde{W}_{\alpha}(g,t)=\tilde{K}_{t}^{\alpha}(g)$.

Thanks to the fundamental relation (\ref{eq:16b2}) that holds in general, we have in particular that
\begin{equation}\label{eq:16b}
\begin{array}{l}
\forall_{t \geq 0} \; \forall_{(\mathbf{x},\mathbf{R}) \in G}\; :\;\tilde{K}_{t}^{\alpha}(\mathbf{x},\mathbf{R})=K_{t}^{\alpha}(\mathbf{x},\mathbf{R}\mathbf{a}).
\end{array}
\end{equation}
Furthermore, the PDE system given by (\ref{PDEgroup}) is invariant under $\mathcal{A}_{i} \mapsto - \mathcal{A}_{i}$, and,
since inversion on the Lie algebra corresponds to inversion on the group,
the kernels must satisfy
\begin{equation} \label{sym2}
\forall_{t \geq 0} \;\forall_{g \in G}\; :\;\tilde{K}_{t}^{\alpha}(g)= \tilde{K}_{t}^{\alpha}(g^{-1}),
\end{equation}
and for the corresponding kernel on the quotient this means
\begin{equation}\label{sym22}
\forall_{t \geq 0} \;\forall_{(\mathbf{x},\mathbf{n}) \in G/H}\; :\; K_{t}^{\alpha}(\mathbf{x},\mathbf{n})=
K_t^{\alpha}(-\mathbf{R}_{\mathbf{n}}^T \mathbf{x}, \mathbf{R}_{\mathbf{n}}^T\mathbf{a}).
\end{equation}

Finally, we see invariance of the PDE w.r.t. right actions of the subgroup $H$. This is due to the isotropy
of the generator $\tilde{Q}_{\alpha}$ in the tangent subbundles $\textrm{span}\{\mathcal{A}_1,\mathcal{A}_{2}\}$
and $\textrm{span}\{\mathcal{A}_4,\mathcal{A}_{5}\}$. This due to (\ref{Za}) in Appendix~\ref{ch:LI}.
Note that invariance of the kernel w.r.t. right action of the subgroup $H$ and invariance of the kernel w.r.t. inversion (\ref{sym2}) also implies invariance of the kernel w.r.t. left-actions of the subgroup $H$, since $(g^{-1}(h')^{-1})^{-1}=h'g$ for all $h' \in H$ and $g\in G$.
Therefore, we have
\begin{equation} \label{nice}
\begin{array}{lll}
\forall_{t \geq 0} \;
\forall_{g \in G}\forall_{h,h' \in H}:\, &\tilde{K}_{t}^{\alpha}(\,g \,h\,) &= \tilde{K}_{t}^{\alpha}(g) = \tilde{K}_{t}^{\alpha}(h'g)
, \\
\forall_{t \geq 0} \;\forall_{(\mathbf{x},\mathbf{n}) \in G/H}\forall_{\overline{\alpha} \in [0,2\pi)}:\, &K_{t}^{\alpha}(\mathbf{x},\mathbf{n}) &=
K_t^{\alpha}(\mathbf{R}_{\mathbf{a},\overline{\alpha}}\mathbf{x}, \mathbf{R}_{\mathbf{a},\overline{\alpha}} \mathbf{n}).
\end{array}
\end{equation}
\begin{Remark}\label{rem:notation} (notations, see also the list of abbreviations at the end of the article)\\
To avoid confusion between the Euler angle $\overline{\alpha}$ and the $\alpha$ indexing the $\alpha$-stable L\'{e}vy distribution, we put an overline for this specific angle. Henceforth  $\mathbf{R}_{\mathbf{v},\psi}$ denotes a counter-clockwise rotation over axis $\mathbf{v}$ with angle $\psi$.
This applies in particular to the case where the axis is the reference axis $\mathbf{v}=\mathbf{a}=(0,0,1)^T$ and
$\psi=\overline{\alpha}$.
 Recall that
$\mathbf{R}_{\mathbf{n}}$ (without an angle in the subscript) denotes any 3D rotation that maps reference axis $\mathbf{a}$ onto $\mathbf{n}$.\\
We write the symbol \, $\hat{\left. \cdot \right.}$ \, above a function to indicate its Fourier transform on\! $G$ and \!$G/H$; we use the  \\
symbol $\overline{\left. \cdot  \right.}$ for strictly spatial Fourier transform; the symbol \, $\tilde{\left. \cdot  \right.}$ \, above a function/operator to indicate that it is defined on the group $G$ and the function/operator without symbols when it is defined on the quotient $G/H$.
\end{Remark}
\subsection{Obtaining the kernels with $D_{11}>0$ from the kernels with $D_{11}=0$ \label{ch:D11}}

In \cite[cor.2.5]{portegies_new_2017} it was deduced 
that for $\alpha=1$ the elliptic diffusion kernel ($D_{11}>0$) directly follows from the hypo-elliptic diffusion kernel ($D_{11}=0$) in the spatial Fourier domain via
\[
\overline{K}_{t}^{1, \textrm{elliptic}}\left(\bomega, \mathbf{n}\right)= e^{-r^2 D_{11} t}\overline{K}_{t}^{1, \textrm{hypo-elliptic}}\left(\sqrt{\frac{D_{33}-D_{11}}{D_{33}}} \bomega, \mathbf{n}\right), \quad \textrm{ with } \ r = \|\bomega\|, \; \; 0 \leq D_{11}<D_{33}.
\]
For the general case $\alpha \in (0,1]$ the transformation from the case $D_{11}=0$  to the case $D_{11}>0$ is achieved by
replacing $-(-\lambda^{l,m}_r)^{\alpha} \mapsto - (-\lambda^{l,m}_r +r^{2} D_{11})^{\alpha} \textrm{ and }$
\mbox{$r \mapsto r \sqrt{\frac{D_{33}-D_{11}}{D_{33}}}$} in formula (\ref{decompose}) for the kernel. Henceforth, we set $D_{11}=0$.

\section{The Fourier Transform on $SE(3)$}\label{ch:FTSE3}

%
%
The group $G=SE(3)$ is a unimodular Lie group (of type I) with (left- and right-invariant) Haar measure ${\rm d}g = {\rm d}\mathbf{x} {\rm d} \mu_{SO(3)}(\mathbf{R})$ being the product
  of the Lebesgue measure on $\mathbb{R}^{3}$ and the Haar-Measure $\mu_{SO(3)}$ on $SO(3)$.
 Then for all $f \in \mathbb{L}_1(G) \cap \mathbb{L}_2(G)$ the Fourier transform $\mathcal{F}_{G}f$ is given by (\ref{FT}).
For more details see \cite{fuhr_abstract_2005,sugiura_unitary_1990,folland_course_1994}.
One has the inversion formula:
 \begin{equation}
f(g) =  (\mF^{-1}_G \mF_G f)(g) = \int_{\hat{G}} \text{trace}\left\{\left(\mF_G f\right)(\sigma)\; \sigma_g \right\} \rmd \nu(\sigma) = \int_{\hat{G}} \text{trace}\left\{\hat{f}(\sigma)\; \sigma_g \right\} \rmd \nu(\sigma).
 \end{equation}
 In our Lie group case of $SE(3)$ we identify all unitary irreducible representations $\sigma^{p,s}$ having non-zero dual measure with the pair $(p,s) \in \mathbb{R}^{+} \times \mathbb{Z}$. This identification is commonly applied, see e.g.~\cite{chirikjian_engineering_2000}. Using the method~\cite[Thm.~\!2.1]{sugiura_unitary_1990}, \cite{mackey_imprimitivity_1949} of induced representations, all unitary irreducible representations (UIR's) of $G$, up to equivalence, with non-zero Plancherel measure are given by:

 \begin{equation} \label{inter}
 \boxed{
 \begin{array}{l}
\sigma = \sigma^{p,s}: SE(3) \to B(\mathbb{L}_{2}(p \, S^2)), \qquad p > 0, \; s \in \mathbb{Z}, \\[7pt] 
 \left(\sigma^{p,s}_{(\bx,\bR)} \phi\right)(\bu) = e^{- i\, \bu \cdot \bx}\; \phi\left(\bR^{-1}\bu\right) \; \Delta_s \left(\bR^{-1}_{\frac{\bu}{p}} \bR \bR_{\frac{\bR^{-1}\bu}{p}}\right), \ \ \ \bu \in p S^{2}, \ \phi \in \mathbb{L}_{2}(p S^2),
 \end{array}
 }
 \end{equation}
where $pS^2$ denotes a 2D sphere of radius $p = \|\bu\|$; $\Delta_s$ is a unitary irreducible representation of $SO(2)$ (or rather of the stabilizing subgroup $\textrm{Stab}_{SO(3)}(\ba) \subset SO(3)$ isomorphic to $SO(2)$) producing a scalar.

In (\ref{inter}), $\bR_{\frac{\bu}{p}}$ denotes a rotation that maps $\ba$ onto $\frac{\bu}{p}$.
So direct computation
\[
\bR^{-1}_{\frac{\bu}{p}} \bR \bR_{\frac{\bR^{-1}\bu}{p}}\ba= \bR^{-1}_{\frac{\bu}{p}} \bR \bR^{-1}\left(\frac{\bu}{p}\right)= \ba
\]
shows us that it is a rotation around the $z$-axis (recall $\ba = \be_z$), say about angle $\overline{\alpha}$. This yields character $\Delta_s \left(\bR^{-1}_{\frac{\bu}{p}} \bR \bR_{\frac{\bR^{-1}\bu}{p}}\right)=e^{- i s \overline{\alpha}}$, for details see \cite[ch.10.6]{chirikjian_engineering_2000}. Thus, we can rewrite~(\ref{inter}) as
$$\left(\sigma^{p,s}_{(\bx,\bR)} \phi\right)(\bu) = e^{- i\, \left(\bu \cdot \bx + s\overline{\alpha}\right)}\; \phi(\bR^{-1}\bu), \ \quad \textrm{ where } (\bx,\bR) \in G, \ \bu \in p S^{2}, \ \phi \in \mathbb{L}_{2}(p S^2) .$$

Mackey's theory~\cite{mackey_imprimitivity_1949} relates the UIR $\sigma^{p,s}$ to the dual orbits $p S^2$ of $SO(3)$. Thereby, the dual measure $\nu$ can be identified with a measure on the family of dual orbits of $SO(3)$ given by $\{p S^2 \, | \, p > 0\}$, and
$$
 \left(\mF^{-1}_G \hat{f}\right)(g) =  \int\limits_{\hat{G}} \operatorname{trace}\left\{\hat{f}(\sigma^{p,s})\; \sigma^{p,s}_g \right\}\; {\rm d} \nu(\sigma^{p,s}) = \int\limits_{\R^{+}} \operatorname{trace}\left\{\hat{f}(\sigma^{p,s})\; \sigma^{p,s}_g \right\} p^2 {\rm d} p,
$$
for all $p>0$, $s \in \mathbb{Z}$. For details see~\cite[ch. 3.6.]{fuhr_abstract_2005}.

The matrix elements of $\hat{f}=\mF_G f$ w.r.t. an orthonormal basis of modified spherical harmonics $\{Y^{l,m}_s(p^{-1}\cdot)\}$, with $|m|, |s| \leq l$, see~\!\cite[ch.9.8]{chirikjian_engineering_2000}, for $\mathbb{L}_{2}(p S^2)$ are given by
\begin{equation} \label{matrixone}
\hat{f}^{p,s}_{l,m,l',m'} := \int \limits_G f(g) \; \left\langle \; \sigma_{g^{-1}}^{p,s} Y^{l',m'}_s(p^{-1}\,\cdot)\,, \, Y^{l,m}_s(p^{-1}\,\cdot)\right\rangle_{\mathbb{L}_{2}(p S^2)}\; \, \rmd g,
\end{equation}
where the $\mathbb{L}_{2}$ inner product is given by
$\left\langle\, y_1(\cdot)\, , \,y_2(\cdot)\, \right\rangle_{\mathbb{L}_{2}(p S^2)} := \left\langle\, y_1(p\cdot)\, , \,y_2(p\cdot)\, \right\rangle_{\mathbb{L}_{2}(S^2)}$, recall (\ref{eq:inerprodL2S2}).

For an explicit formula for the modified spherical harmonics $Y^{l,m}_s$ see \cite{chirikjian_engineering_2000}, where they are denoted by $h_{m,s}^{l}$. The precise technical analytic expansion of the modified spherical harmonics is not important for this article.
The only properties of $Y^{l,m}_s$ that we will need are gathered in the next proposition.
\begin{Proposition} \label{prop}
The modified spherical harmonics $Y^{l,m}_s$ have the following properties:
\begin{equation*}
\begin{array}{l}
\textrm{1)\ for }s= 0 \textrm{ or }m=0 \textrm{ they coincide with standard spherical harmonics $Y^{l,m}$, cf.~\!\cite[eq.4.32]{griffiths} : } \\[7pt]
\ \ Y^{l,m}_{s=0}= Y^{l,m} \textrm{ and }Y^{l,0}_s=(-1)^s \, Y^{l,s}, \textrm{ where }Y^{l,m}(\mathbf{n}(\beta,\gamma))= \frac{\epsilon_m}{\sqrt{2\pi}}\, P^{m}_{l}(\cos \beta) \, e^{i m \gamma}, \\[5pt]
\textrm{ with } \mathbf{n}(\beta,\gamma)=(\cos \gamma \sin \beta, \sin \gamma \sin \beta, \cos \beta)^T,   \textrm{ with spherical angles }   \beta \in [0,\pi], \gamma \in [0,2\pi), \\[5pt]
\textrm{ with $P^{m}_l$ the normalized associated Legendre polynomial\iftoggle{insertremarks}{\footnotemark}{}
and  $\epsilon_m=(-1)^{\frac{1}{2}\left(m+|m|\right)}$;} \\[5pt]
\textrm{2)\ they have a specific rotation transformation property in view of (\ref{inter}): } \\[3pt]
\ \ \sigma^{p,s}_{(\mathbf{0}, \mathbf{R})}  Y^{l,m}_s= \sum \limits_{m'=-l}^{l} \mathcal{D}^l_{m'm}(\mathbf{R}) \;Y^{l,m'}_{s},
\textrm{ where $\mathcal{D}^{l}_{m'm} (\cdot)$ denotes the Wigner D-function~\cite{Wigner};} \\[5pt]
\textrm{3)\ for each $s \in \mathbb{Z}$ fixed they form a complete orthonormal basis for }\mathbb{L}_{2}(S^2): \\[7pt]
\  \left\langle Y^{l,m}_s(\cdot), Y^{l',m'}_{s}(\cdot)\right\rangle_{\mathbb{L}_{2}(S^2)}=\delta^{l, l'} \delta^{m, m'} \textrm{ for all } m, m' \in \mathbb{Z}, \,  l, l' \in \mathbb{N}_0,  \; \textrm{ with } |m|\leq l, |m'|\leq l', \; l,l' \geq |s|.
\end{array}
\end{equation*}
\end{Proposition}
\footnotetext{
We take the convention $P^m_l(x)=N^{l,m} (1- x^2)^{\frac{|m|}{2}} \frac{d^{|m|}P_l(x)}{dx^{|m|}}$ with $N^{l,m}=
\sqrt{\frac{(2l+1)\, (l-|m|)!}{2(l+|m|)!}}$,
$P_{l}(x)=\frac{1}{2^l l!}\frac{d^l}{dx^l}(x^2-1)^l$ for $|x|\leq 1$.}
For details and relation between different Euler angle conventions, see~\cite[ch:9.4.1]{chirikjian_engineering_2000}.
In our convention of ZYZ-Euler angles, see Appendix~\ref{ch:LI}, one has
\begin{equation}\label{Wigner}
\mathcal{D}^{l}_{m'm}(
\mathbf{R}_{\mathbf{e}_{z},\overline{\alpha}}
\mathbf{R}_{\mathbf{e}_{y},\beta}
\mathbf{R}_{\mathbf{e}_{z},\gamma}
)= e^{-im' \overline{\alpha}} P^{l}_{m'm}(\cos \beta)e^{-im\gamma}, \end{equation}
with $P^{l}_{m'm}$ a generalized associated Legendre polynomial given in~\cite[eq.9.21]{chirikjian_engineering_2000}.
\iftoggle{insertremarks}{
\begin{Remark}
The Wigner D-functions in Proposition~\ref{prop} are independent of $s\in\mathbb{Z}$ due to the construction of $Y^{l,m}_s$, cf.~\!\cite[ch.9,10.7]{chirikjian_engineering_2000}.
Now since $Y^{l,m}_{s=0}= Y^{l,m}$ this means that
the basis of modified spherical harmonics $\{Y^{l,m}_{s}\}_{m=-l}^l$ is chosen such that the matrix representation of
operator $\sigma_{(\mathbf{0},\mathbf{R})}^{p,s}$ relative to this basis is the same as the matrix representation of the rotation operator
$\sigma_{(\mathbf{0},\mathbf{R})}^{p,s=0}$ relative to the basis of spherical harmonics $\{Y^{l,m}\}_{m=-l}^l$.
\end{Remark}
}{

}
\noindent
Moreover, we have inversion formula (\cite[Eq.10.46]{chirikjian_engineering_2000}):
\begin{equation} \label{C5}
f(g) = \frac{1}{2 \pi^2} \sum_{s \in \mathbb{Z}} \sum_{l' = |s|}^{\infty} \sum_{l = |s|}^\infty \sum_{m'=-l'}^{l'} \sum_{m=-l}^l \int \limits_0^\infty \hat{f}^{p,s}_{l,m,l',m'}\, \left(\sigma^{p,s}_g\right)_{l',m',l,m} \; p^2 \rmd p,
\end{equation}
with matrix coefficients (independent of $f$) given by
\begin{equation} \label{matrix}
\left(\sigma^{p,s}_g\right)_{l',m',l,m}= \left\langle \, \sigma^{p,s}_g Y^{l,m}_s(p^{-1} \cdot) \, , \, Y^{l',m'}_{s}(p^{-1} \cdot) \,\right\rangle_{\mathbb{L}_{2}(p S^2)}.
\end{equation}
Note that $\sigma^{p,s}$ is a UIR so we have
\begin{equation}\label{property}
\left(\sigma^{p,s}_{g^{-1}}\right)_{l',m',l,m}= \overline{\left(\sigma^{p,s}_g\right)_{l,m,l',m'}}\,.
\end{equation}

\section{A Specific Fourier Transform on the Homogeneous Space $\R^{3}\rtimes S^2$ of Positions and Orientations \label{ch:FTquot}}
Now that we have introduced the notation of Fourier transform on the Lie group $G = SE(3)$, we will define the Fourier transform $\mathcal{F}_{G/H}$ on the homogeneous space $G/H = \quot$.
Afterwards, in the subsequent section, we will solve the Forward-Kolmogorov/Fokker-Planck PDEs~(\ref{eq:differentialQ}) via application of this transform. Or more precisely, via conjugation with Fourier transform $\mathcal{F}_{G/H}$.
\subsection{The Homogeneous Space $\R^{3}\rtimes S^2$  \label{ch:FTquot2a}}

Throughout this manuscript we shall rely on a Fourier transform on the homogeneous space
of positions and orientations that is defined by the partition of left-cosets:
$\RR^3 \rtimes S^2:=G/H$, given by (\ref{quot}).

Note that subgroup $H$ can be parameterized as follows:
\begin{equation} \label{subgroup}
H=\{h_{\overline{\alpha}}:=(\mathbf{0}, \mathbf{R}_{\mathbf{a},\overline{\alpha}}) \;|\; \overline{\alpha} \in [0,2\pi)\},
\end{equation}
where we recall, that $\bR_{\mathbf{a},\overline{\alpha}}$ denotes a (counter-clockwise) rotation around the reference axis $\mathbf{a} =\be_z$.
The reason behind this construction is that the group $SE(3)$ acts transitively on $\R^{3} \rtimes S^{2}$
by $(\bx',\bn') \mapsto g \odot (\mathbf{x}',\mathbf{n}')$ given by (\ref{actiononset}).
Recall that by the definition of the left-cosets one has
\[
H = \{\mathbf{0}\} \times SO(2), \textrm{ and }g_{1}\sim g_{2} \desda g_{1}^{-1}g_2 \in H.
\]
The latter equivalence simply means that
for $
g_{1}=(\mathbf{x}_{1},\mathbf{R}_{1})$ and $g_{2}=(\mathbf{x}_{2},\mathbf{R}_{2})$ one has
\[
g_{1}\sim g_{2} \desda \mathbf{x}_{1}=\mathbf{x}_{2} \textrm{ and }\exists_{\overline{\alpha} \in [0,2\pi)}\;:\; \mathbf{R}_{1}= \mathbf{R}_{2} \mathbf{R}_{\mathbf{a},\overline{\alpha}}.
\]
The equivalence classes $[g]=\{g' \in SE(3)\; |\; g'\sim g\}$ are often just denoted by $(\mathbf{x},\mathbf{n})$ as they consist of all rigid body motions
$g=(\mathbf{x},\mathbf{R}_{\mathbf{n}})$ that map reference point $(\mathbf{0},\mathbf{a})$ onto $(\mathbf{x},\mathbf{n}) \in \R^3 \rtimes S^2$\ :
\begin{equation}\label{eq:qoutfromgroup}
g \odot (\mathbf{0},\mathbf{a})=(\mathbf{x},\mathbf{n}),
\end{equation}
where we recall $\mathbf{R}_{\mathbf{n}}$ is \emph{any} rotation that maps $\mathbf{a} \in S^2$ onto $\mathbf{n} \in S^2$.

\subsection{Fourier Transform on $\RR^3 \rtimes S^2$ \label{ch:FTquot2}}
Now we can define the Fourier transform $\mathcal{F}_{G/H}$ on the homogeneous space $G/H$. Prior this, we specify a class of functions where this transform acts.
\begin{Definition}
Let $p>0$ be fixed and $s\in \mathbb{Z}$.
We denote
\[
\mathbb{L}_2^{sym}(p S^2) = \left\{\left. f \in \mathbb{L}_{2}(p S^2) \; \right|\; \forall_{\overline{\alpha} \in [0,2\pi)} \; \sigma^{p,s}_{h_{\overline{\alpha}}}f=f\right\}
\]
the subspace of spherical functions that have the prescribed axial symmetry, w.r.t. the subgroup $H$, recall (\ref{subgroup}).
\end{Definition}
\begin{Definition}
We denote the orthogonal projection from $\mathbb{L}_2(p S^2)$ onto the closed subspace $\mathbb{L}_2^{sym}(p S^2)$ by $\ProjS$.
\end{Definition}
\begin{Definition}\label{def:quotrep}
To the group representation $\sigma^{p,s}: SE(3) \to B(\mathbb{L}_{2}(p S^2))$ given by~\!(\ref{inter}),
we relate a ``representation'' $\overline{\sigma}^{p,s}: \R^{3} \rtimes S^{2} \to B(\mathbb{L}_{2}(p S^2))$ on $\quot$, defined by
\begin{equation}\label{coordinatefree}
\overline{\sigma}^{p,s}_{[g]}:= \frac{1}{(2\pi)^2}\int \limits_{0}^{2\pi} \int \limits_{0}^{2\pi} \sigma_{h_{\tilde{\alpha}}g h_{\overline{\alpha}}}^{p,s} \, {\rm d}\overline{\alpha} {\rm d} \tilde{\alpha} = \ProjS \circ \sigma^{p,s}_{g} \circ \ProjS.
\end{equation}
\end{Definition}
\begin{Definition}
A function $\tilde{U}: G \to \mathbb{C}$ is called axially symmetric if
\begin{equation} \label{AS}
\tilde{U}(\mathbf{x},\mathbf{R})= \tilde{U}(\mathbf{x}, \mathbf{R} \mathbf{R}_{\mathbf{a},\overline{\alpha}}) \ \ \ \textrm{ for all }\overline{\alpha} \in [0,2\pi)  \textrm{ and all }(\mathbf{x},\mathbf{R}) \in G.
\end{equation}
\end{Definition}

To each function $U:G/H \to \mathbb{C}$ we relate an axially symmetric function $\tilde{U}:G \to \mathbb{C}$ by
\begin{equation} \label{UtildeU}
\tilde{U}(\mathbf{x},\mathbf{R}):=U(\mathbf{x},\mathbf{R} \mathbf{a}).
\end{equation}
\begin{Definition}
We define the Fourier transform of function $U$ on $G/H=\R^{3} \rtimes S^{2}$ by
\begin{equation} \label{FThom}
\boxed{
\hat{U}(\overline{\sigma}^{p,s}) = \left(\mathcal{F}_{G/H}U\right)(\overline{\sigma}^{p,s}):= \ProjS \circ \mathcal{F}_{G}\tilde{U}(\sigma^{p,s})\circ \ProjS.
}
\end{equation}
\end{Definition}
Standard properties of the Fourier transform $\mathcal{F}_{G}$ on $SE(3)$ such as the Plancherel theorem and the inversion formula~\cite{chirikjian_engineering_2000,sugiura_unitary_1990} naturally carry over to
$\mathcal{F}_{G/H}$ with `simpler formulas'. This is done by a domain and range restriction via the projection operators $\ProjS$ in (\ref{FThom}).
The reason for the specific construction~(\ref{FThom}) will become clear from the next lemmas, and the `simpler formulas' for the Plancherel and inversion formulas are then summarized in a subsequent theorem, where we constrain ourselves to the case $m=m'=0$ in the formulas. The operator $\ProjS$ that is most right in~(\ref{FThom}) constrains the basis $Y_{s}^{l,m}$ to $m=0$, whereas the operator $\ProjS$ that is most left in (\ref{FThom}) constrains the basis $Y^{l',m'}_s$ to $m'=0$.
\begin{Lemma} \label{Lemma:axial} (axial symmetry)
Let $\tilde{U}: G \to \mathbb{C}$ be axially symmetric. Then
\begin{enumerate}
\item
it relates to a unique function $U: G/H \to \mathbb{C}$  via $U(\mathbf{x},\mathbf{n})=\tilde{U}(\mathbf{x},\mathbf{R}_{\mathbf{n}})$;
\item
the matrix coefficients
\[
\hat{U}^{p,s}_{l,m,l',m'}=
\left[\mathcal{F}_{G}\tilde{U}(\sigma^{p,s})\right]_{l,m,l',m'} \textrm{ of linear operator }\mathcal{F}_{G}\tilde{U}(\sigma^{p,s})
\]
relative to the modified spherical harmonic basis $\{Y_{s}^{l,m}\}$ vanish if $m\neq 0$.
\item
the matrix coefficients
\[
\hat{U}^{p,s}_{l,m,l',m'}=
\left[\mathcal{F}_{G/H}U(\overline{\sigma}^{p,s})\right]_{l,m,l',m'} \textrm{ of linear operator }\mathcal{F}_{G/H}U(\overline{\sigma}^{p,s})
\]
relative to the modified spherical harmonic basis $\{Y_{s}^{l,m}\}$ vanish if $m\neq 0$ or $m' \neq 0$;
\end{enumerate}
Conversely, if $\tilde{U}=\mathcal{F}_{G}^{-1}(\hat{U})$ and
\begin{equation} \label{conversely}
\forall_{p>0} \forall_{l \in \mathbb{N}_0}\forall_{s \in \mathbb{Z}, \textrm{ with }|s|\leq l} \; \forall_{m' \in \mathbb{Z}, \textrm{ with }|m'|\leq l} \;
\forall_{m \neq 0}\;:\; \hat{U}^{p,s}_{l,m,l',m'}=0,
\end{equation}
then $\tilde{U}$ satisfies the axial symmetry (\ref{AS}).
\end{Lemma}
\begin{proof}
Item 1:
Uniqueness of $U$ follows by the fact that the choice of $\mathbf{R}_{\mathbf{n}}$ of some rotation that maps $\mathbf{a}$ onto $\mathbf{n}$ does not matter. Indeed $U(\mathbf{x},\mathbf{n})=\tilde{U}(\mathbf{x},\mathbf{R}_{\mathbf{n}}\mathbf{R}_{\mathbf{a},\overline{\alpha}})=\tilde{U}(\mathbf{x},\mathbf{R}_{\mathbf{n}})$.

Item 2: Assumption (\ref{AS}) can be rewritten as $\tilde{U}(g)=\tilde{U}(gh_{\overline{\alpha}})$ for all $h_{\overline{\alpha}} \in H$, $g \in G$. This gives:
\begin{equation} \label{derivation}
\begin{array}{ll}
\hat{U}^{p,s}_{l,m,l',m'} &=\left\langle\;(\mathcal{F}_{G}\tilde{U})(Y^{l',m'}_s(p^{-1}\cdot))\, ,\, Y^{l,m}_s(p^{-1}\cdot)\;\right\rangle_{\mathbb{L}_{2}(pS^2)} \\[6pt] &=\int \limits_{G} \tilde{U}(g) \; \left\langle\;
\sigma_{g^{-1}}^{p,s} Y^{l',m'}_s(p^{-1}\cdot)\, ,\, Y^{l,m}_s(p^{-1}\cdot)\;\right\rangle_{\mathbb{L}_{2}(pS^2)} \; {\rm d}g \\
&= \int \limits_{G} \tilde{U}(g) \; \left\langle\;   Y^{l',m'}_s(p^{-1}\cdot)\; \, ,\, \sigma_{g}^{p,s} Y^{l,m}_s(p^{-1}\cdot) \right\rangle_{\mathbb{L}_{2}(pS^2)} \; \, {\rm d}g \\
 &= \int \limits_{G} \tilde{U}(g h_{\overline{\alpha}}) \; \left\langle\;  Y^{l',m'}_s(p^{-1}\cdot)\, ,\, \sigma_{g h_{\overline{\alpha}}}^{p,s} Y^{l,m}_s(p^{-1}\cdot) \;\right\rangle_{\mathbb{L}_{2}(pS^2)} \; {\rm d}(g h_{\overline{\alpha}}) \\
 &= \int \limits_{G} \tilde{U}(g ) \; \left\langle\; Y^{l',m'}_s(p^{-1}\cdot)\, ,\,  \sigma_{g}^{p,s} \circ \sigma_{h_{\overline{\alpha}}}^{p,s} Y^{l,m}_s(p^{-1}\cdot)\;\right\rangle_{\mathbb{L}_{2}(pS^2)} \; {\rm d}(g h_{\overline{\alpha}}) \\
 &= e^{-im \overline{\alpha}}\; \hat{U}^{p,s}_{l,m,l',m'} \textrm{ for all }\overline{\alpha} \in [0,2\pi),
\end{array}
\end{equation}
where we recall that $\sigma$ is a UIR and that the Haar-measure on $G$ is bi-invariant.
In the first step we used the 3rd property whereas in the final step we used the 2nd property
of Proposition~\ref{prop} together with
\begin{equation} \label{alphasym}
\mathcal{D}^l_{m'm}(\mathbf{R}_{\ba,\overline{\alpha}})= e^{-i m \overline{\alpha}} \delta_{m'm}\ \textrm{ so that }
\sigma_{h_{\overline{\alpha}}}^{p,s} Y^{l,m}_s(p^{-1}\cdot) = e^{-i m \overline{\alpha}} Y^{l,m}_{s}(p^{-1}\cdot).
\end{equation}
We conclude that $(1-e^{-i m \overline{\alpha}}) \hat{U}^{p,s}_{l,m,l',m'}=0$ for all $\overline{\alpha} \in [0,2\pi)$ so $m \neq 0 \Rightarrow  \hat{U}^{p,s}_{l,m,l',m'}=0$.

Item 3:
Due to the 2nd property in Proposition~\ref{prop} we have
\[
\sigma^{p,s}_{(\mathbf{0},\mathbf{R})} Y^{l,m}_s(p^{-1}\cdot)=
\sum \limits_{m'=-l}^l \mathcal{D}^{l}_{m'm}(\mathbf{R})\, \; Y^{l,m'}_{s}(p^{-1}\cdot).
\]
Thereby the projection $\ProjS$ is given by
\begin{equation}\label{eq:ProjSym}
\ProjS \left(\sum \limits_{l=0}^{\infty} \sum \limits_{m=-l}^l \alpha_{l,m} Y^{l,m}_s \right)= \sum \limits_{l=0}^{\infty}\alpha_{l,0} Y^{l,0}_s.
\end{equation}
Now the projection $\ProjS$ that is applied first in (\ref{FThom}) filters out $m=0$ as the only possible nonzero component. The second projection filters out $m'=0$ as the only possible nonzero component.

Conversely, if (\ref{conversely}) holds one has by inversion formula (\ref{C5}) that
\[
\tilde{U}(g)= \frac{1}{2\pi^2} \sum \limits_{s \in \mathbb{Z}} \sum \limits_{l=|s|}^{\infty} \sum \limits_{l'=|s|}^{\infty} \sum \limits_{m'=-l'}^{l'}\; \int\limits_0^{\infty}
\hat{U}^{p,s}_{l,0,l',m'} \; \left(\sigma^{p,s}_g\right)_{l',m',l,0} \; p^2 {\rm d}p,
\]
so then the final result follows by the identity
\begin{equation} \label{symmetr}
\left(\sigma_{g h_{\overline{\alpha}}}^{p,s}\right)_{l',m',l,0}=
\left(\sigma_{g}^{p,s}\right)_{l',m',l,0}.
\end{equation}
So it remains to show why Eq.~\!(\ref{symmetr}) holds.
It is due to $\sigma_{(\mathbf{x},\mathbf{R})}^{p,s}=
\sigma_{(\mathbf{x},\mathbf{I})}^{p,s} \circ \sigma_{(\mathbf{0},\mathbf{R})}^{p,s}$ and (\ref{alphasym}), as one has
\begin{equation} \label{sigsym1}
\begin{array}{l}
\sigma_{g h_{\overline{\alpha}}}^{p,s}=\sigma^{p,s}_{(\mathbf{x},\mathbf{R}) (\mathbf{0}, \mathbf{R}_{\ba,\overline{\alpha}} )}= \sigma^{p,s}_{(\mathbf{x},\mathbf{R} \mathbf{R}_{\ba,\overline{\alpha}})}=
\sigma_{(\mathbf{x},\mathbf{R})}^{p,s} \circ \sigma^{p,s}_{(\mathbf{0},\mathbf{R}_{\ba,\overline{\alpha}})}, \textrm{  and }
Y^{l,0}_s(p^{-1}\mathbf{R}^{-1}_{\ba,\overline{\alpha}}\cdot)=Y^{l,0}_s(p^{-1}\cdot)
\end{array}
\end{equation}
and thereby Eq.~\!(\ref{symmetr}) follows
by Eq.\!~(\ref{matrix}).
\end{proof}
\begin{Lemma}
If $\tilde{K} \in \mathbb{L}_2(G)$ is real-valued and satisfies the axial symmetry (\ref{AS}), and moreover the following holds
 \begin{equation} \label{inversionsymmetry}
 \tilde{K}(g^{-1})=
 \tilde{K}(g)
 \end{equation}
then the Fourier coefficients satisfy $\hat{K}^{p,s}_{l,m,l',m'}=
\overline{\hat{K}^{p,s}_{l',m',l,m}}$ and they vanish for $m \neq 0$ and for $m'\neq 0$.

\end{Lemma}
\begin{proof}
\iftoggle{insertremarks}
{
From Lemma~\ref{Lemma:axial} we know that if $\tilde{K}$ satisfies the axial symmetry then
$\hat{K}^{p,s}_{l,m,l',m'}$ vanishes if $m \neq 0$. Furthermore, $\hat{K}^{p,s}_{l,m,l',m'}=
\overline{\hat{K}^{p,s}_{l',m',l,m}}$ holds due to
\[
\begin{array}{ll}
\hat{K}^{p,s}_{l',m',l,m} &=\left\langle\left(\mathcal{F}_{G} \tilde{K}\right) Y_{s}^{l,m}(p^{-1}\cdot)\, ,\, Y_{s}^{l',m'}(p^{-1}\cdot)\right\rangle_{\mathbb{L}_{2}(pS^2)}=
\int \limits_{G} \tilde{K}(g) \, \left(\sigma_{g^{-1}}^{p,s}\right)_{l',m',l,m} \, {\rm d}g \\ &= \int \limits_{G} \tilde{K}(g^{-1}) \, \left(\sigma_{g}^{p,s}\right)_{l',m',l,m} \, {\rm d} (g^{-1}) = \int \limits_{G} \tilde{K}(g) \, \left(\sigma_{g}^{p,s}\right)_{l',m',l,m} \, {\rm d}g = \int \limits_{G} \tilde{K}(g) \, \overline{\left(\sigma_{g^{-1}}^{p,s}\right)_{l,m,l',m'}} \, {\rm d}g \\[10pt]
&= \overline{\hat{K}^{p,s}_{l,m,l',m'}}.
\end{array}
\]
where at the end we used that $\sigma^{p,s}$ is a unitary representation.
Then applying Lemma~\ref{Lemma:axial} to the axial symmetric function $\tilde{K}$ and to the coefficient $\overline{\hat{K}^{p,s}_{l',m',l,m}}$ 
yields the result.
}
{
The proof follows by Eq.~\!(\ref{property}) and inversion invariance of the Haar measure on $G$,  see~\cite{arxivFT}. 
}
\end{proof}
The next lemma shows that~(\ref{inversionsymmetry}) is a sufficient but not a necessary condition for the Fourier coefficients to vanish for both the cases $m'\neq 0$ and $m\neq 0$.
\begin{Lemma} \label{Lemma:crucial}
Let $\tilde{K} \in \mathbb{L}_{2}(G)$ and $K \in \mathbb{L}_{2}(G/H)$ be related by (\ref{UtildeU}). Then we have the following equivalences:
 \begin{equation}
 \begin{array}{c}
K(\mathbf{x},\mathbf{n})= K(\mathbf{R}_{\mathbf{a},\overline{\alpha}} \mathbf{x}, \mathbf{R}_{\mathbf{a},\overline{\alpha}} \mathbf{n}), \ \ \
\textrm{ for all }\overline{\alpha} \in [0,2 \pi),
(\mathbf{x},\mathbf{n}) \in G/H
  \\
 \Updownarrow \\
 \tilde{K}(g h)=
 \tilde{K}(g)= \tilde{K}(h g), \ \ \  \textrm{ for all }g \in G, \ h \in H \label{weaksymmetry} \\
 \Updownarrow \\
 \textrm{The Fourier coefficients $\hat{K}^{p,s}_{l,m,l',m'}$ vanish for $m \neq 0$ and for $m'\neq 0$. }
 \end{array}
 \end{equation}
\end{Lemma}
\begin{proof}
We show $a \Rightarrow b \Rightarrow c \Rightarrow a$ to get $a \desda b \desda c$. \\
$a\Rightarrow b$: Denoting $h=h_{\overline{\alpha}}=
(\mathbf{0},\bR_{\mathbf{a},\overline{\alpha}})$,
  $g=(\mathbf{x},\mathbf{R})$, we have
\[
\begin{array}{r l}
\forall_{\overline{\alpha},\overline{\alpha}' \in [0,2\pi)} \forall_{\mathbf{x} \in \R^3} \forall_{\mathbf{R} \in SO(3)}\;:\;
\tilde{K}(g h_{\overline{\alpha}})\!\!\! & =
\tilde{K}(\mathbf{x}, \mathbf{R} \mathbf{R}_{\mathbf{a},\overline{\alpha}})=
K(\mathbf{x}, \mathbf{R} \mathbf{R}_{\mathbf{a},\overline{\alpha}} \mathbf{a})=K(\mathbf{x}, \mathbf{R} \mathbf{a})
=
\tilde{K}(\mathbf{x}, \mathbf{R})=\tilde{K}(g)\\
& =
K(\mathbf{R}_{\mathbf{a},\overline{\alpha}} \mathbf{x},  \mathbf{R}_{\mathbf{a},\overline{\alpha}} \mathbf{R} \mathbf{a})=
\tilde{K}( \mathbf{R}_{\mathbf{a},\overline{\alpha}} \mathbf{x},  \mathbf{R}_{\mathbf{a},\overline{\alpha}} \mathbf{R})=\tilde{K}( h_{\overline{\alpha}}g).
\end{array}
\]
$b \Rightarrow c$: By Lemma~\ref{Lemma:axial} we know that the Fourier coefficients vanish for $m \neq 0$. Next we show they also vanish for $m' \neq 0$. Similar to (\ref{sigsym1}) we have
\begin{equation} \label{sigsym2}
\begin{array}{l}
\sigma_{h_{\overline{\alpha}} g }^{p,s}
=
\sigma_{(\mathbf{R}_{\ba,\overline{\alpha}}\mathbf{x},\mathbf{R}_{\ba,\overline{\alpha}}\mathbf{R})}^{p,s}=
\sigma_{(\mathbf{R}_{\ba,\overline{\alpha}}\mathbf{x},\mathbf{I})}^{p,s} \circ \sigma_{(\mathbf{0},\mathbf{R}_{\ba,\overline{\alpha}}\mathbf{R})}^{p,s}\ ,
 \end{array}
\end{equation}
which gives the following relation for the matrix-coefficients:
\begin{equation} \label{inbetween}
\begin{array}{l}
\left(\sigma^{p,s}_{g=(\mathbf{x},\mathbf{R})}\right)_{l',m',l,m} = \sum \limits_{j=-l}^l \left\langle \sigma^{p,s}_{(\mathbf{x},\mathbf{I})}Y^{l,j}_s(p^{-1}\cdot)\, , \, Y^{l',m'}_s (p^{-1}\cdot)\right\rangle_{\mathbb{L}_{2}(pS^2)} \;
\mathcal{D}^l_{jm}(\mathbf{R})\qquad \Rightarrow \\
\left(\sigma^{p,s}_{h_{\overline{\alpha}}\, g}\right)_{l',m',l,m} =\sum \limits_{j=-l}^l e^{-i(m'-j)\overline{\alpha}}\left\langle \sigma^{p,s}_{(\mathbf{x},\mathbf{I})}Y^{l,j}_s(p^{-1}\cdot) \, , \, Y^{l',m'}_s(p^{-1}\cdot)\right\rangle_{\mathbb{L}_{2}(pS^2)} \;
e^{-i\, j \overline{\alpha}}\, \mathcal{D}^l_{jm}(\mathbf{R}) \Rightarrow \\[8pt]
\left(\sigma^{p,s}_{h_{\overline{\alpha}}g}\right)_{l',m',l,m} =e^{-im' \overline{\alpha}}\; \left(\sigma^{p,s}_g\right)_{l',m',l,m}. \\
\end{array}
\end{equation}
The implication can be directly verified by Proposition~\ref{prop}, (\ref{Wigner}), (\ref{sigsym2}),
and
\[
\begin{array}{ll}
\left\langle Y^{l',m'}_s(p^{-1}\cdot)\, , \, \sigma^{p,s}_{(\mathbf{R}_{\ba,\overline{\alpha}}\mathbf{x},\mathbf{I})}Y^{l,j}_s(p^{-1}\cdot)\right\rangle_{\mathbb{L}_{2}(pS^2)} &= \int \limits_{pS^2} e^{-i p (\mathbf{x} \cdot \mathbf{R}_{\ba,\overline{\alpha}}^{T} \mathbf{u})}\;
Y^{l,j}_s(\mathbf{u})\; \overline{Y^{l',m'}_s(\mathbf{u})}\, {\rm d}_{\mu_{pS^2}}(\mathbf{u})\\
 &=
\int \limits_{pS^2} e^{-i p (\mathbf{x} \cdot  \mathbf{v})}
Y^{l,j}_s(\mathbf{R}_{\ba,\overline{\alpha}}\mathbf{v}) \; \overline{Y^{l',m'}_s(\mathbf{R}_{\ba,\overline{\alpha}}\mathbf{v})}\, {\rm d}_{\mu_{pS^2}}(\mathbf{v}).
\end{array}
\]
From (\ref{inbetween}) we deduce that:
\[
\begin{array}{ll}
\hat{K}^{p,s}_{l,m,l',m'} &= \int \limits_{G} \tilde{K}(g) \; \left\langle\; \sigma_{g}^{p,s} Y^{l,m}_s(p^{-1}\cdot) \,,\,  Y^{l',m'}_s(p^{-1}\cdot)\;\right\rangle_{\mathbb{L}_{2}(pS^2)} \; \, {\rm d}g \\
 &= \int \limits_{G} \tilde{K}(h_{\overline{\alpha}} g) \; \left\langle\; \sigma_{h_{\overline{\alpha}}g}^{p,s} Y^{l,m}_s(p^{-1}\cdot) \,,\,  Y^{l',m'}_s(p^{-1}\cdot)\;\right\rangle_{\mathbb{L}_{2}(pS^2)} \; {\rm d}(h_{\overline{\alpha}}g) \\
 &= \int \limits_{G} \tilde{K}(g ) \; \left\langle \; \sigma_{g}^{p,s}  Y^{l,m}_s(p^{-1}\cdot) \,,\, \sigma_{h_{\overline{\alpha}}^{-1}}^{p,s} Y^{l',m'}_s(p^{-1}\cdot)\;\right\rangle_{\mathbb{L}_{2}(pS^2)} \; {\rm d}g = e^{+im'\, \overline{\alpha}}\; \hat{K}^{p,s}_{l,m,l',m'},
\end{array}
\]
which holds for all $\overline{\alpha} \in [0,2\pi)$. Thereby if $m'\neq 0$ then $\hat{K}^{p,s}_{l,m,l',m'}=0$.

$c \Rightarrow a$: By inversion formula (\ref{C5}), where the only contributing terms have $m=0$ and $m'=0$, we see that $\tilde{K}(gh)=\tilde{K}(hg)=\tilde{K}(g)$ for all
$h=(\mathbf{0},\mathbf{R}_{\ba,\overline{\alpha}})$. Thereby $\tilde{K}$ is axially symmetric
and by Lemma~\ref{Lemma:axial} it relates to a unique kernel on $G/H$ via $K(\bx,\bn)=\tilde{K}(\bx,\bR_{\bn})$ and the result follows by (\ref{nice}).
\end{proof}
  Now that we characterized all functions $K \in \mathbb{L}_{2}(G/H)$ for which the Fourier coefficients $\hat{K}^{p,s}_{l,m,l',m'}$
  vanish for $m\neq 0$ and $m'\neq 0$ in the above lemma, we will considerably simplify the inversion and Plancherel formula for Fourier transform   $\mathcal{F}_{G}$ on the group $G=SE(3)$
  to the Fourier transform $\mathcal{F}_{G/H}$ on the homogeneous space $G/H=\R^{3}\rtimes S^2$ in the next theorem.  This will be important in our objective of deriving the kernels for the linear PDEs (\ref{eq:differentialQ}) that we address in the next section.
\begin{Theorem} (matrix-representation for $\mathcal{F}_{G/H}$, explicit inversion and Plancherel formula)  \label{corr:1}\\
Let $K \in \mathbb{L}_{2}^{sym}(G/H)$ and $\tilde{K} \in \mathbb{L}_{2}(G)$ be related by~(\ref{UtildeU}).
Then the matrix elements of $\mathcal{F}_{G/H}K$ are given by
\[
\begin{array}{ll}
\hat{K}^{p,s}_{l',0,l,0} &=\int \limits_{G} \tilde{K}(g) \, \left(\sigma_{g^{-1}}^{p,s}\right)_{l',0,l,0} \, {\rm d}g\ , \\
\textrm{with }\left(\sigma_{g}^{p,s}\right)_{l',0,l,0} &=\sum \limits_{j=-l}^l \; \left[l',0 \; |\; p,s \;|l,j\right](\mathbf{x})\;
\mathcal{D}^l_{j0}(\mathbf{R})\ \  \  \textrm{ for all }g=(\mathbf{x},\mathbf{R}) \in G.
\end{array}
\]
The constants
$\left[l',0 \; |\; p,s \;|l,j \right](\mathbf{x}):= \left\langle\; \sigma^{p,s}_{(\mathbf{x},\mathbf{I})}Y^{l,j}_s(p^{-1}\cdot)\, , \, Y^{l',0}_s(p^{-1}\cdot)\;\right\rangle_{\mathbb{L}_{2}(pS^2)}$
admit an analytic expression in terms of elementary functions \cite[Eq.10.34]{chirikjian_engineering_2000} and the Wigner D-functions~(\ref{Wigner}).
\\[6pt]
Furthermore, we have the following Plancherel and inversion formula:
\[
\begin{array}{ll}
\|K\|^{2}_{\mathbb{L}_{2}(G/H)} &= \|\mathcal{F}_{G/H} K\|^2= \sum \limits_{s \in \mathbb{Z}}\,\int \limits_{\R^+}  \|\!|
\left(\mathcal{F}_{G/H} K\right)(\overline{\sigma}^{p,s})\|\!|^2 \; p^2 {\rm d}p=
\int \limits_{\R^+} \sum \limits_{s = -\infty}^{\infty}
\sum \limits_{l'=|s|}^{\infty} \sum \limits_{l=|s|}^{\infty}
|\hat{K}^{p,s}_{l,0,l',0}|^2 \; p^2 {\rm d}p ,  
 \\[12pt]
K(\mathbf{x},\mathbf{n}) &= \left(\mathcal{F}_{G/H}^{-1}\mathcal{F}_{G/H}K\right)(\mathbf{x},\mathbf{n}) 
=\sum \limits_{s \in \mathbb{Z}} \ \int\limits_{\R^+} \operatorname{trace}\left\{(\mathcal{F}_{G/H} K)(\overline{\sigma}^{p,s}) \; \overline{\sigma}^{p,s}_{(\mathbf{x},\mathbf{n})}\right\}\; p^2 {\rm d}p \\[6pt]
 &=\frac{1}{2\pi^2} \sum \limits_{s \in \mathbb{Z}} \sum \limits_{l'=|s|}^{\infty} \sum \limits_{l=|s|}^{\infty} \
\int \limits_{\R^+} \hat{K}^{p,s}_{l,0,l',0}\, \left(\overline{\sigma}^{p,s}_{(\mathbf{x}, \mathbf{n})}\right)_{l',0,l,0} \; p^2 {\rm d}p,
\end{array}
\]
with matrix coefficients given by (for analytic formulas see \cite[eq.10.35]{chirikjian_engineering_2000})
\begin{equation}\label{matrixcoeff}
\begin{array}{ll}
\left(\overline{\sigma}^{p,s}_{(\mathbf{x},\mathbf{n})}\right)_{l',0,l,0} &=
\left(\sigma^{p,s}_{g}\right)_{l',0,l,0}
 =\left\langle\; \sigma^{p,s}_g Y_{s}^{l,0}(p^{-1}\cdot)\, ,\, Y_{s}^{l',0}(p^{-1}\cdot)\;\right\rangle_{\mathbb{L}_{2}(p S^2)}\\
  &=
 \left\langle\; \sigma^{p,s}_g Y^{l,s}(p^{-1}\cdot)\, ,\, Y^{l',s}(p^{-1}\cdot)\;\right\rangle_{\mathbb{L}_{2}(p S^2)} \ \  \textrm{ for }g=(\mathbf{x},\mathbf{R}_{\mathbf{n}}).
\end{array}
\end{equation}
\end{Theorem}
\begin{proof}
The above formulas are a direct consequence of Lemma~\ref{Lemma:crucial} and the Plancherel and inversion formulas (see \cite{sugiura_unitary_1990}, \cite[ch:10.8]{chirikjian_engineering_2000}) for Fourier transform on $SE(3)$. Recall that a coordinate-free definition of $\overline{\sigma}^{p,s}$ is given in (\ref{coordinatefree}). Its matrix coefficients are given by
(\ref{matrixcoeff}), where we recall the 1st item of Prop.~\ref{prop} and where we note that they are independent on the choice of $\mathbf{R}_{\mathbf{n}} \in SO(3)$ mapping $\mathbf{a}$ onto $\mathbf{n}$. 
\end{proof}
\begin{Corollary}\label{corr:new}
Let $K_1$, $K_2 \in \mathbb{L}_2^{sym}(G/H)$. Then for shift-twist convolution on $G/H = \quot$ given by
$$(K_1 * K_2)(\bx,\bn) = \int \limits_{S^2} \int \limits_{\R^{3}}
K_1(\bR_{\bn'}^T(\bx-\bx'),\bR_{\bn'}^T\bn)\, K_2(\bx',\bn') \; {\rm d}\bx'{\rm d}\mu_{S^2}(\bn')$$
we have
$\mF_{G/H}(K_1 * K_2) =(\mF_{G/H}K_1) \circ (\mF_{G/H}K_2).$
\end{Corollary}
\begin{proof}
Set $\tilde{K}_1(g) \!=\! K_1(g \odot (\mathbf{0},\ba))$. Standard Fourier theory~\cite{chirikjian_stochastic_2011} gives
$\mF_{G}(\widetilde{K_1 * K_2})\!=\!\mF_{G}(\tK_1 * \tK_2)$, so
$$
\begin{array}{l c l}
\mF_{G/H} (K_1 * K_2) & \eqdef  & \ProjS \circ \mF_{G}(\widetilde{K_1 * K_2})\circ \ProjS \\
 & = & \ProjS \circ \mF_{G}(\tK_1)\circ \mF_{G}(\tK_2)\circ \ProjS \\
 & = &  \ProjS \circ \mF_{G}(\tK_1)\circ \ProjS \circ \ProjS \circ \mF_{G}(\tK_2)\circ \ProjS\\
 & = &  (\mF_{G/H}K_1) \circ (\mF_{G/H}K_2),
\end{array}
$$
where the 1st equality is given by~(\ref{FThom}) and the 3rd equality follows by Lemma~\ref{Lemma:crucial} and~(\ref{eq:ProjSym}).
\end{proof}

\iftoggle{insertremarks}{\newpage}{}
\section{Application of the Fourier Transform on $\RR^{3} \rtimes S^{2}$ for
Explicit Solutions of the Fokker-Plank PDEs of $\alpha$-stable L\'{e}vy Processes on $\RR^{3} \rtimes S^{2}$}\label{ch:FTPDE}

Our objective is to solve PDE system (\ref{eq:differentialQ}) on the homogeneous space of positions and orientations $G/H$.
Recall that we extended this PDE system to $G$ in~(\ref{PDEgroup}). As the cases $D_{11}>0$ follow from the case $D_{11}=0$ (recall Subsection~\ref{ch:D11}), we consider the case $D_{11}=0$ in this section.
From the symmetry consideration in Section~\ref{ch:symmetry} it follows that the solution of~(\ref{PDEgroup}) is given by
$\tilde{W}_{\alpha}(g,t)= (\tilde{K}_{t}^{\alpha} * \tilde{U})(g)$ with a probability kernel $\tilde{K}_{t}^{\alpha}:G \to \R^+$,
whereas the solution of~(\ref{eq:differentialQ}) is given by
\[
W_{\alpha}(\bx,\bn,t) = (K_t^{\alpha} \ast U)(\bx,\bn) := \int \limits_{S^2} \int \limits_{\R^{3}}
K_{t}^{\alpha}(\bR_{\bn'}^T(\bx-\bx'), \bR_{\bn'}^T \bn) \; U(\bx',\bn') \; {\rm d}\bx'{\rm d}\mu_{S^2}(\bn'),
\]
where the kernels $K_{t}^{\alpha}$ are invariant with respect to left-actions of the subgroup $H$,
recall Eq.~\!(\ref{nice}).  This invariance means that the condition for application of the Fourier transform $\mathcal{F}_{G/H}$ on $\R^{3}\rtimes S^2$ is satisfied (recall Lemma~\ref{Lemma:crucial}) and we can indeed employ Theorem~\ref{corr:1} to keep all
our computations, spectral decompositions and Fourier transforms in the 5D homogeneous space $\R^{3}\rtimes S^2 = G/H$  rather than a technical and less direct approach \cite{portegies_new_2017} in the 6D group $G=SE(3)$.

\begin{Remark}
For the underlying probability theory, and sample paths of discrete random walks of the $\alpha$-Stable L\'{e}vy stochastic processes we refer to Appendix~\ref{ch:PT}. To get a general impression of how Monte Carlo simulations of such stochastic processes can be used to approximate the exact probability kernels $K_{t}^{\alpha}$, see Fig.~\!\ref{fig:sample-paths}. In essence, such a stochastic approximation is computed by binning the endpoints of the random walks.
A brief mathematical explanation will follow in Subsection~\ref{ch:montecarlo}.
\end{Remark}
For now let us ignore the probability theory details and let us first focus on deriving exact analytic solutions to~(\ref{eq:differentialQ}) and its kernel $K_{t}^{\alpha}$ via Fourier transform $\mathcal{F}_{G/H}$ on $G/H=\R^{3}\rtimes S^2$.

\subsection{Exact Kernel Representations by Spectral Decomposition in the Fourier Domain}

Let us consider evolution~(\ref{eq:differentialQ}) for $\alpha$-stable L\'{e}vy process on the quotient $G/H=\R^{3}\rtimes S^2$.
Then the mapping from the initial condition $W(\cdot,0)=U(\cdot) \in \mathbb{L}_{2}(G/H)$ to the solution $W(\cdot,t)$ at a fixed time $t\geq0$ is a bounded linear mapping. It gives rise to a strongly continuous (holomorphic) semigroup~\cite{yosida_functional_1980}.
\iftoggle{insertremarks}{While keeping in mind that the domain of the evolution generator $Q_{\alpha} = -(-Q)^{\alpha}$ is given by (\ref{domain}) w}{W}e will conveniently denote the bounded linear operator on $\mathbb{L}_{2}(G/H)$ as follows:
\begin{equation}\label{eq:etQ}
W_{\alpha}(\cdot,t)=(e^{ t Q_{\alpha}} U)(\cdot),  \qquad \textrm{ for all $t\geq 0$.} 
\end{equation}

In the next main theorem we will provide a spectral decomposition of the operator using both a direct sum and a direct integral decomposition. Note that definitions of direct integral decompositions (and the underlying measure theory) can be found in \cite[ch:3.3\&3.4]{fuhr_abstract_2005}.

\subsubsection{Eigenfunctions and Preliminaries}
\begin{figure}[t!]
   \centering
   \includegraphics[width=0.38\textwidth]{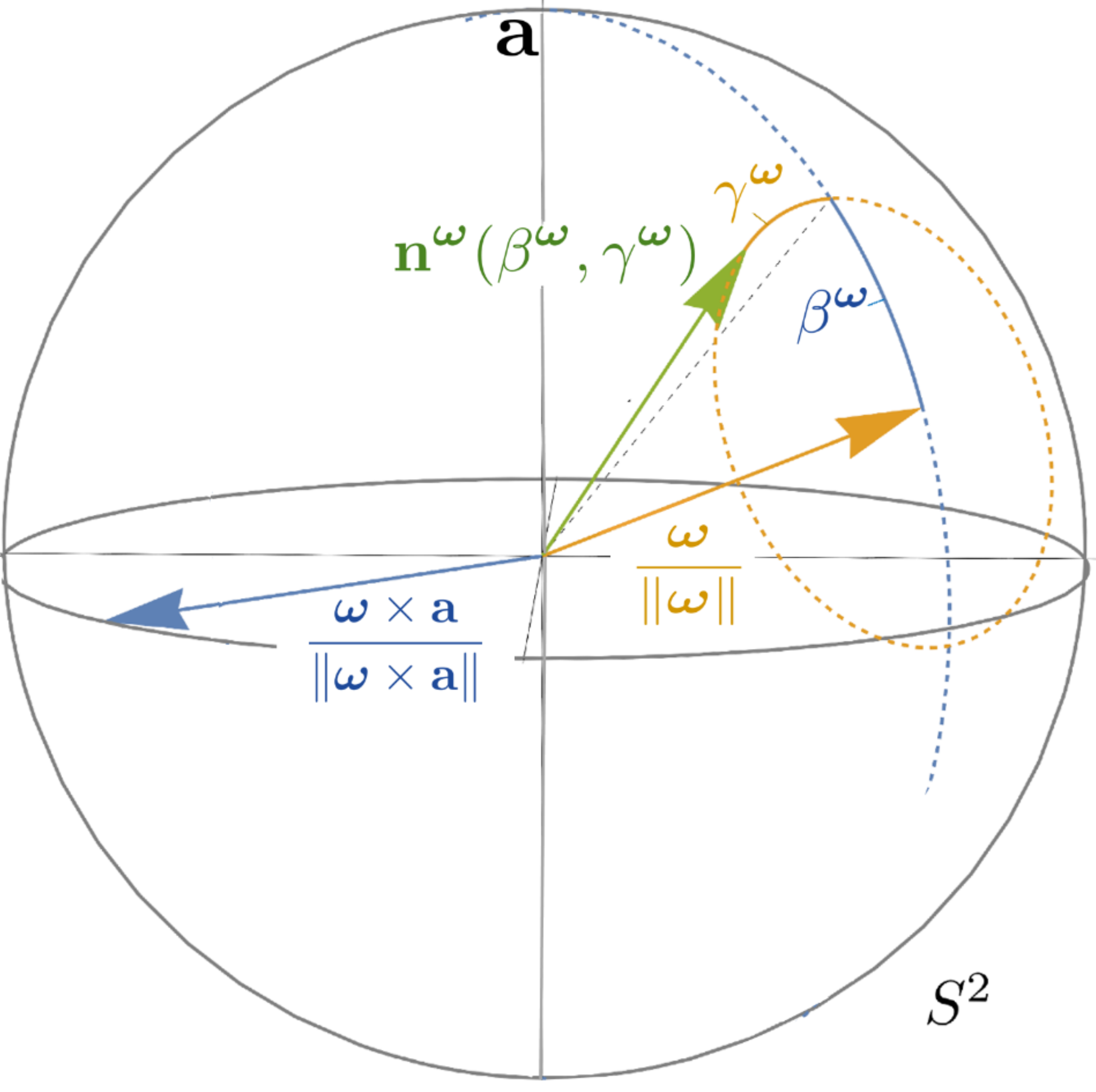}
 \caption{For $\bomega \neq \mathbf{a}$, we parameterize every orientation $\bn$ (green) by rotations around $r^{-1} \bomega$ (orange) and $\frac{\bomega \times \mathbf{a}}{||\bomega \times \mathbf{a}||}$ (blue). In other words, $\bn^{\bomega} (\beta^{\bomega},\gamma^{\bomega}) = \bR_{r^{-1} \bomega,\gamma^{\bomega}} \bR_{\frac{\bomega \times \mathbf{a}}{||\bomega \times \mathbf{a}||}, \beta^{\bomega}} (r^{-1} \bomega )$.}\label{fig:figureParametrization}
 \end{figure}

In order to formulate the main theorem
we need some preliminaries and formalities. First of all let us define
$\overline{\mathcal{F}}_{\R^3}: \mathbb{L}_{2}(\R^{3}\rtimes S^2) \to \mathbb{L}_{2}(\R^{3}\rtimes S^2)$
by
\begin{equation} \label{Ftilde}
(\overline{\mathcal{F}}_{\R^3}U)(\boldsymbol{\omega},\mathbf{n}):=
\left[ \mathcal{F}_{\R^{3}}U(\cdot,\mathbf{n})
\right](\boldsymbol{\omega}).
\end{equation}
Recall (\ref{correspondence}). Then we re-express the generator in the spatial Fourier domain:
\begin{equation} \label{reexpress}
\begin{array}{ll}
-(-\mathcal{B})^{\alpha} &=
\overline{\mathcal{F}}_{\R^3} \circ Q_{\alpha} \circ \overline{\mathcal{F}}_{\R^3}^{-1} \Rightarrow \\[6pt]
-\left(-\mathcal{B}_{\bomega} \right)^{\alpha} &=
-\left(-D_{33} \,(i \boldsymbol{\omega} \cdot \mathbf{n})^2 - D_{44}\, \Delta^{S^2}_{\bn}\right)^{\alpha}\\
  &= -\left(D_{33} \,r^2\, \left(\mathbf{a}\cdot (\mathbf{R}^T_{r^{-1}\boldsymbol{\omega}}\mathbf{n})\right)^2 - D_{44}\,\Delta^{S^2}_{\bn}\right)^{\alpha} \\[6pt]
  &= -\left(D_{33}\, r^2 \, \cos^{2}(\beta^{\boldsymbol{\omega}}) - D_{44}\, \Delta^{S^2}_{\bn}\right)^{\alpha}, \ \ \textrm{ with }r=\|\boldsymbol{\omega}\|, \alpha \in (0,1],
\end{array}
\end{equation}
where $\beta^{\boldsymbol{\omega}}$ denotes the angle between $\mathbf{n}$ and $r^{-1}\bomega$, see Fig.~\!\ref{fig:figureParametrization}.
This re-expression is the main reason for the following definitions. 

Instead of the modified spherical Harmonics $Y^{l,m}_s$ in Prop.~\ref{prop}, which are commonly used as a standard basis to represent each operator in the Fourier transform on $SE(3)$, we will use our generalized spherical harmonics, depending on a spatial frequency vector as this is in accordance with~(\ref{reexpress}).
\begin{Definition} \label{def:1}
Let $l \in \mathbb{N}_0$. Let $m \in \mathbb{Z}$ such that $|m|\leq l$.
Let $\boldsymbol{\omega} \in \R^3$ be a frequency vector. We define
\begin{equation} \label{SHw}
Y_{\boldsymbol{\omega}}^{l,m}(\mathbf{n})= Y^{l,m}(\mathbf{R}_{r^{-1}\boldsymbol{\omega}}^{T}\mathbf{n}), \qquad \textrm{ with }r=\|\boldsymbol{\omega}\|, \ \ \bn \in S^2,
\end{equation}
where we take the rotation
which maps $\mathbf{a}$ onto $r^{-1}\boldsymbol{\omega}$ whose matrix representation in the standard basis is:
\[
\bR_{r^{-1} \bomega } =  \begin{pmatrix}
\frac{(\bomega \times \mathbf{a}) \times \bomega}{||(\bomega \times \mathbf{a}) \times \bomega||} & \vline & \frac{\bomega \times \mathbf{a}}{||\bomega \times \mathbf{a}||} & \vline & r^{-1}\bomega\
\end{pmatrix}  \textrm{ for }r^{-1}\boldsymbol{\omega} \neq \mathbf{a}, \textrm{ and }\bR_{\mathbf{a}}=\bI,  \textrm{ and }\bR_{\mathbf{0}}=\bI.
\]
\end{Definition}
\noindent
Recall the standard spherical angle formula
$\mathbf{n}(\beta,\gamma)=
(\sin \beta \cos \gamma,  \sin \beta \sin \gamma,\cos \beta)^T$ from Proposition~\ref{prop}.
These are Euler-angles relative to the reference axis $\mathbf{a}=\mathbf{e}_{z}$.
For the Euler-angles relative to the (normalized) frequency $r^{-1} \bomega$ one has (see also Fig.~\ref{fig:figureParametrization}):
\begin{equation}\label{nw}
\mathbf{n}^{\bomega}(\beta^{\bomega},\gamma^{\bomega})=\bR_{r^{-1} \bomega } \mathbf{n}(\beta^{\bomega},\gamma^{\bomega}).
\end{equation}
\begin{Definition} \label{def:2}
Let $l \in \mathbb{N}_0$. Let $m \in \mathbb{Z}$ such that $|m|\leq l$.
We define the functions $\Phi_{\bomega}^{l,m} \in \mathbb{L}_{2}(S^{2})$ by
\begin{equation} \label{Philm}
\Phi^{l,m}_{\bomega}(\mathbf{n})= \sum \limits_{j=0}^{\infty} \frac{d_j^{l,m}(r)}{\|\mathbf{d}^{l,m}\left(r \right)\|} \; Y_{\bomega}^{|m|+j,m}(\mathbf{n}),
\end{equation}
where $r = \|\bomega\|$ and $\mathbf{d}^{l,m}\left(r\right):=\left(d_j^{l,m}\left(r\right)\right)_{j=0}^{\infty}$ are coefficients such that
\[
\Phi^{l,m}_{\bomega}(\mathbf{n}^{\bomega}(\beta^{\bomega},\gamma^{\bomega}))= S_{\rho}^{l,m}(\cos \beta^{\boldsymbol{\omega}}) \; \frac{e^{i m \gamma^{\boldsymbol{\omega}}}}{\sqrt{2\pi}},  \textrm{ with }\rho= r\sqrt{\frac{D_{33}}{D_{44}}},
\]
where $S_{\rho}^{l,m}(\cdot)$ denotes
the $\mathbb{L}_{2}$-normalized spheroidal wave function.
%
\begin{Remark} \label{rem:completebasis}
The spheroidal wave function arises from application of the method of separation on operator $\mathcal{B}_{\bomega}$ in (\ref{reexpress}) where
basic computations (for details see \cite{portegies_new_2017}) lead to the following singular Sturm-Liouville problem:
\begin{equation}
(L y)(x) = \frac{d}{dx}\left[ p(x) \frac{dy(x)}{dx} \right] + q(x) y(x) = -\lambda(r) \; y(x), \quad x=\cos \beta^{\boldsymbol{\omega}} \in [-1,1].
\end{equation}
with $p(x) = (1-x^2)$, $q(x) = -\rho^2 x^2 - \frac{m^2}{1-x^2}$, and again $\rho=r \sqrt{D_{33}/D_{44}}$. In this formulation, $p(x)$ vanishes at the boundary of the interval, which makes our problem a singular Sturm-Liouville problem. It is sufficient to require boundedness of the solution and its derivative at the boundary points to have nonnegative, distinct, simple eigenvalues $\lambda^{l,m}_r$ and existence of a countable, complete orthonormal basis of eigenfunctions $\{y_j\}_{j=0}^{\infty}$ \cite{margenau_mathematics_1956} for the spheroidal wave equation.

As a result 
 standard Sturm-Liouville theory (that applies the spectral decomposition theorem for compact self-adjoint operators to a kernel operator that is the right-inverse of $L$), provides us (for each $\bomega$ fixed) a complete orthonormal basis of eigenfunctions
$\{\Phi^{l,m}_{\bomega}\}$
in $\mathbb{L}_{2}(S^2)$ with eigenvalues of our (unbounded) generators:
\begin{equation} \label{efunc}
 -\left(-\mathcal{B}_{\bomega} \right)^{\alpha}\Phi^{l,m}_{\bomega}= -(-\lambda^{l,m}_r)^{\alpha} \; \Phi^{l,m}_{\bomega}, \quad \textrm{ for all } |m| \leq l.
\end{equation}
\end{Remark}
\begin{Remark}
Define $\mathcal{Y}_{l,m}(\beta,\gamma) := Y^{l,m}(\mathbf{n}(\beta,\gamma))$. Then (\ref{SHw}), (\ref{nw})
imply $Y^{l,m}_{\bomega}(\mathbf{n}^{\bomega}(\beta^{\bomega},\gamma^{\bomega}))=\mathcal{Y}_{l,m}(\beta^{\bomega},\gamma^{\bomega})$.
\end{Remark}
\iftoggle{insertremarks}{
\begin{Remark}
We have the special cases
\begin{equation} \label{sc}
\Phi^{l,m}_{\mathbf{0}}=\Phi^{l,m}_{\mathbf{a}}=Y^{l,m}
\end{equation} and if $D_{33}=0$ then $\Phi^{l,m}_{\bomega}=Y^{l,m}$. In general, $\Phi^{l,m}_{\boldsymbol{\omega}}(\cdot)$ is to be considered as a basis function in $\mathbb{L}_{2}(S^2)$ that depends both on the radius $r=\|\bomega\|$ and the direction of the frequency
$r^{-1}\bomega$ according to (\ref{Philm}). The coefficients in (\ref{Philm}) depend only on the radius and the basis functions depend on the direction.
\end{Remark}
}{}
\begin{Remark}
The matrix-representation of $ -\left(-\mathcal{B}_{\bomega} \right)^{\alpha}$
w.r.t. orthonormal basis $\left\{Y^{|m|+j,m}_{\bomega}\right\}_{j \in \mathbb{N}_0, m \in \mathbb{Z}}$ equals
\[
\bigoplus_{m \in \mathbb{Z}}
-(D_{33}r^2\mathbf{M}^{m} + D_{44}\boldsymbol{\Lambda}^{m})^{\alpha},
\]
where $\boldsymbol{\Lambda}^m:=\operatorname{diag}\{l(l+1)\}_{l=|m|}^{\infty}=\operatorname{diag}\{(|m|+j)(|m|+j+1)\}_{j=0}^{\infty}$, $r=\|\bomega\|$ and
where $\mathbf{M}^{m}$ is the tri-diagonal matrix (that can be computed analytically \cite[eq.~106]{portegies_new_2017})
given by
\begin{equation} \label{matrixM}
(\cos \beta)^2 Y^{|m|+j,m}(\mathbf{n}(\beta,\gamma))= \sum \limits_{j'=0}^{\infty} \left((\mathbf{M}^{m})^T\right)_{j,j'} Y^{|m|+j',m}(\mathbf{n}(\beta,\gamma)).
\end{equation}
As a result, we see from (\ref{Philm}), (\ref{efunc}) that the coefficients $\mathbf{d}^{l,m}\left(r\right)$ for our eigenfunctions are eigenvectors of a matrix
\begin{equation} \label{efmatrix}
-\left(D_{33}r^2\mathbf{M}^{m} + D_{44}\boldsymbol{\Lambda}^{m}\right) \mathbf{d}^{l,m}(r) = \lambda^{l,m}_r \mathbf{d}^{l,m}(r), \qquad \textrm{ for }l\geq |m|.
\end{equation}
This matrix (and its diagonalization) will play a central role for our main spectral decomposition theorem both in the spatial Fourier domain and in the Fourier domain of the homogeneous space of positions and orientations.
\end{Remark}
\end{Definition}

\subsubsection{The Explicit Spectral Decomposition of the Evolution Operators}\label{sec:Expl}

In Theorem~\ref{th:decomposition} we will present the explicit spectral decompositions both in the spatial Fourier domain and in the Fourier domain of the homogeneous space of positions and orientations.

Prior to this theorem we explain the challenges that appear when we apply $\mathcal{F}_{G/H}$ to the PDE of interest (\ref{eq:differentialQ}) on the quotient $G/H$. In order to get a grip on the evolution operator and
the corresponding kernel, we set the initial condition equal to a delta distribution at the origin, i.e. we consider
\[
U= \delta_{(\mathbf{0},\mathbf{a})} \Rightarrow W_{\alpha}(\cdot,t)=e^{t Q_{\alpha}}U=e^{-t (-Q)^{\alpha}} \delta_{(\mathbf{0},\mathbf{a})} =K_{t}^{\alpha}.
\]
In this case, the necessary condition~(\ref{weaksymmetry}) in Lemma~\ref{Lemma:crucial} for application of $\mathcal{F}_{G/H}$ is indeed satisfied, due to
the symmetry property of the kernel, given by~(\ref{nice}).
Now due to linearity
\[
\mathcal{F}_{G/H} \circ e^{t Q_{\alpha}} \circ
\mathcal{F}_{G/H}^{-1}= e^{t\, \left(\mathcal{F}_{G/H}\; \circ \; Q_{\alpha} \;\circ \;\mathcal{F}_{G/H}^{-1}\right)},
\]
we just need to study the generator in the Fourier domain.

For the moment we set $\alpha=1$ (the hypo-elliptic diffusion case) and return to the general case later on (recall Subsection~\ref{rem:3} and Subsection~\ref{ch:D11}). Then it follows that (for details see \cite[App.D]{portegies_new_2017})
\begin{equation} \label{GenFquot}
\begin{array}{l}
\left(\mathcal{F}_{G/H} \circ Q \circ \mathcal{F}_{G/H}^{-1} \hat{K}_t^{1}\right)(\overline{\sigma}^{p,s})=
\left(
-  D_{33}\, (\mathbf{a} \cdot \mathbf{u})^2 +
D_{44}\, \Delta^{p S^2}_{\mathbf{u}} \right) \hat{K}_{t}^{1}(\overline{\sigma}^{p,s}) , \\[6pt]
\textrm{with the kernel } \hat{K}_{t}^{1}:= \mathcal{F}_{G/H}K_t^{1}(\cdot).
\end{array}
\end{equation}
Here $\Delta^{p S^2}_{\mathbf{u}}$ denotes the Laplace-Beltrami operator on a sphere $pS^2\!=\!\{\bu \in \R^3 \big| \|\bu\| = p\}$ of radius~$p>0$.\\
We recall that $\mathbf{u} \in p S^{2}$ is the variable of the functions on which $\overline{\sigma}^{p,s}$ acts. Recall Eq.~\!(\ref{inter}).
So the first part in the righthand side of (\ref{GenFquot}) denotes a multiplier operator $\mathcal{M}$ given by
\[
(\mathcal{M}\phi)(\mathbf{u}) := 
  - \, (\mathbf{a} \cdot \mathbf{u})^2 \phi(\mathbf{u}), \textrm{ for all }\phi \in \mathbb{L}_{2}(pS^2), \textrm{ and almost every } \mathbf{u} \in p S^{2}.
\]
As a result we obtain the following PDE system for $\hat{K}_t^{\alpha}$ (now for general $\alpha \in (0,1]$):
\[
\boxed{
\left\{
\begin{array}{rl}
\frac{\partial}{\partial t}
\hat{K}_{t}^{\alpha}(\overline{\sigma}^{p,s}) &=
-\left(
 -D_{33}\,\mathcal{M}  -
D_{44}\, \Delta^{p S^2}_{\mathbf{u}}\right)^{\alpha}  \hat{K}_{t}^{\alpha}(\overline{\sigma}^{p,s}) \\
\hat{K}_{0}^{\alpha}(\overline{\sigma}^{p,s}) &=\textrm{1}_{\mathbb{L}_{2}(pS^2)}.
\end{array}
\right.
}
\]
\begin{Remark}
There is a striking analogy between the operators $\mathcal{F}_{G/H} \circ Q_{\alpha} \circ \mathcal{F}_{G/H}^{-1}$  and
$\overline{\mathcal{F}}_{\R^3} \circ Q_{\alpha} \circ \overline{\mathcal{F}}_{\R^3}^{-1}$ given by~(\ref{reexpress}), where the role of $r \mathbf{R}_{\bomega/r}^{T} \mathbf{n}$
corresponds to $\mathbf{u}$. This correspondence ensures that the multipliers of the multiplier operators in the generator
coincide and that the roles of $p$ and $r$ coincide:
\[
\mathbf{u}=r \mathbf{R}_{r^{-1}\bomega}^{T} \mathbf{n} \Rightarrow
(\mathbf{a} \cdot \mathbf{u})^2= r^2 (\mathbf{R}_{r^{-1}\bomega}^{T}\mathbf{a} \cdot \mathbf{n})^2 =(\bomega \cdot \mathbf{n})^2
\textrm{ \ \ and \ \ }\|\mathbf{u}\|=p=r=\|\bomega\|.
\]
\end{Remark}
\begin{Lemma}\label{lemma:mr}
Let $t\geq 0$ and $p>0$ be fixed.
The matrix-representation of operator $e^{t(
  D_{33}\,\mathcal{M}  +
D_{44}\, \Delta^{p S^2}_{\mathbf{u}})}: \mathbb{L}_{2}(pS^2) \to \mathbb{L}_{2}(pS^2)$
w.r.t. the orthonormal basis of spherical harmonics $\left\{Y^{l=|s|+j\;,\;s}(p^{-1}\cdot)\right\}_{j \in \mathbb{N}_0\ , \ s \in \mathbb{Z}}$ equals
\begin{equation}\label{eq:matrepexptQ}
\bigoplus_{s \in \mathbb{Z}}
e^{-t\, (D_{33}p^2\mathbf{M}^{s} + D_{44}\boldsymbol{\Lambda}^{s})} \
.
\end{equation}
\end{Lemma}
\begin{proof} Recall~(\ref{matrixM}) that defines matrix $\mathbf{M}^{m}$ (for analytic formulas of this tri-diagonal matrix see~\cite{portegies_new_2017}).
This may be re-written as follows:
\[
(\mathbf{a} \cdot \mathbf{n})^2 Y^{|m|+j,m}(\mathbf{n})= \sum \limits_{j'=0}^{\infty} \left((\mathbf{M}^{m})^T\right)_{j,j'} Y^{|m|+j',m}(\mathbf{n}).
\]
Now fix $s \in \mathbb{Z}$ and set $m=s$ and $\mathbf{n}= p^{-1}\mathbf{u}$ and we have:
\[
 \left\langle \left(D_{33} \mathcal{M}  +
D_{44}\, \Delta^{p S^2}_{\mathbf{u}}\right) Y^{l,s}(p^{-1}\cdot)\; , \; Y^{l',s}(p^{-1}\cdot) \; \right\rangle_{\mathbb{L}_{2}(pS^2)}=
-p^2\, D_{33} \left(\mathbf{M}^s\right)_{j',j} - D_{44} l (l+1) \delta_{jj'},
\]where again $l=|s|+j$, $l'=|s|+j'$ and $j,j' \in \mathbb{N}_0$.

Finally, we note that operator $ D_{33}\,\mathcal{M}  +
D_{44}\, \Delta^{p S^2}_{\mathbf{u}}$ is negative definite and maps each subspace $\operatorname{span}\left\{\{Y^{l,s}(p^{-1}\cdot)\}_{l=|s|}^{\infty}\right\}$ for fixed
$s \in \mathbb{Z}$ onto itself, which explains direct sum decomposition~(\ref{eq:matrepexptQ}).
\end{proof}
Next we formulate the main result, where we apply a standard identification of tensors $\mathbf{a} \otimes \mathbf{b}$ with linear maps:
\begin{equation} \label{tensormap}
\mathbf{x} \mapsto (\mathbf{a} \otimes \mathbf{b})(\mathbf{x})= \left\langle\mathbf{x}\, , \, \mathbf{b}\right\rangle\mathbf{a}.
\end{equation}
\begin{Theorem} \label{th:decomposition}
We have the following spectral decompositions for the Forward-Kolomogorov evolution operator of $\alpha$-stable Levy-processes on
the homogeneous space $G/H=\R^{3}\rtimes S^2$:
\begin{itemize} \item In the Fourier domain of the homogeneous space of positions and orientations we have:
\begin{equation} \label{FD1}
\boxed{
\begin{array}{l}
\mathcal{F}_{G/H} \circ e^{-t(-Q)^{\alpha}} \circ
\mathcal{F}_{G/H}^{-1} \\
\qquad =
\int \limits_{\R^+}^{\oplus}  \bigoplus \limits_{s \in \mathbb{Z}} \sum \limits_{l,l'=|s|}^{\infty}
\left[e^{-(D_{33}p^2\mathbf{M}^{s} + D_{44}\bLam^{s})^{\alpha}t }\right]_{l,l'} \left(Y^{l,s}(p^{-1}\cdot) \otimes Y^{l',s}(p^{-1}\cdot)\right) \; p^{2}{\rm d}p\\
 \qquad = \int \limits_{\R^+}^{\oplus}
 \bigoplus \limits_{s \in \mathbb{Z}} \sum \limits_{l=|s|}^{\infty}\; e^{-(- \lambda^{l,s}_p)^{\alpha} t}\;
 \left(\Phi^{l,s}_{p\mathbf{a}}(p^{-1}\cdot) \otimes \Phi^{l,s}_{p\mathbf{a}}(p^{-1}\cdot)\right) \; p^{2}{\rm d}p
 \end{array}
 }
\end{equation}
\item In the spatial Fourier domain we have
\begin{equation} \label{FD2}
\boxed{
\begin{array}{l}
\left(
\overline{\mathcal{F}}_{\R^3} \circ e^{-t(-Q)^{\alpha}} \circ
\overline{\mathcal{F}}_{\R^3}^{-1} \overline{U}\right)(\bomega,\cdot) =\overline{W}(\bomega,\cdot,t) \\ \qquad =
\sum \limits_{m \in \mathbb{Z}} \sum \limits_{l,l'=|m|}^{\infty}
\left[e^{-(D_{33}r^2\mathbf{M}^{m} + D_{44}\bLam^{m})^{\alpha}t }\right]_{l,l'} \left(Y^{l,m}_{\bomega} \otimes Y^{l',m}_{\bomega}\right)(\overline{U}(\bomega,\cdot)) \\ \qquad
  = 
 \sum \limits_{m \in \mathbb{Z}} \sum \limits_{l=|m|}^{\infty} e^{-(- \lambda^{l,m}_r)^{\alpha} t}\;
 \left(\Phi^{l,m}_{\boldsymbol{\omega}} \otimes \Phi^{l,m}_{\boldsymbol{\omega}}\right)(\overline{U}(\bomega,\cdot))
 \end{array}
 }
  \end{equation}
 where $\overline{W}(\bomega,\cdot,t)=\overline{\mathcal{F}}_{\R^3}W(\bomega,\cdot,t)$ and
 $\overline{U}(\bomega,\cdot)=\overline{\mathcal{F}}_{\R^3}U(\bomega,\cdot)$, recall (\ref{Ftilde}).
\end{itemize}
In both cases the normalized eigenfunctions $\Phi_{\bomega}^{l,m}$ are given by (\ref{Philm})
in Definition~\ref{def:2}. The eigenvalues $ \lambda^{l,m}_r$ are the eigenvalues of the spheroidal wave equation as explained in Remark~\ref{rem:completebasis}.
\end{Theorem}
\begin{proof}
The first identity (\ref{FD1}) follows by:
\[
\begin{array}{rl}
\mathcal{F}_{G/H} \circ e^{-t(-Q)^{\alpha}} \circ
\mathcal{F}_{G/H}^{-1} =&
e^{t \left(\mathcal{F}_{G/H} \circ -(-Q)^{\alpha} \circ
\mathcal{F}_{G/H}^{-1}\right)} \\
 \overset{\textrm{{\scriptsize  \cite[App.D]{portegies_new_2017} \& Thm.~\ref{corr:1}}}}{=}&
 \int \limits^{\oplus}_{\R^+} e^{- t \left(-D_{33} \mathcal{M} + D_{44}\Delta^{p S^2}_{\mathbf{u}}\right)^{\alpha}} p^2{\rm d}p \\
 \overset{\textrm{{\scriptsize Lemma~\ref{lemma:mr}\,\& \,Thm.~\ref{corr:1}}}}{=}&
 \int \limits_{\R^+}^{\oplus} \bigoplus \limits_{s \in \mathbb{Z}} \sum \limits_{l,l'=|s|}^{\infty}
\left[e^{-t \left(D_{33}p^2\mathbf{M}^{s} + D_{44}\bLam^{s}\right)^{\alpha}}\right]_{l,l'} \left(Y^{l,s}(p^{-1}\cdot) \otimes Y^{l',s}(p^{-1}\cdot) \right)\; p^{2}{\rm d}p \\
 \overset{(\ref{Philm})\iftoggle{insertremarks}{,(\ref{sc})}{}}{=}&
  \int \limits_{\R^+}^{\oplus} \bigoplus \limits_{s \in \mathbb{Z}} \sum \limits_{l=|s|}^{\infty} e^{-(- \lambda^{l,s}_p)^{\alpha} t}\;
 \left(\Phi^{l,s}_{p\boldsymbol{a}}(p^{-1}\cdot) \otimes \Phi^{l,s}_{p\boldsymbol{a}}(p^{-1}\cdot)\right) \; p^{2}{\rm d}p\ .
\end{array}
\]
\iftoggle{insertremarks}{}{In the last equality we use the fact that $\Phi^{l,m}_{\mathbf{a}}=Y^{l,m}$. }By applying identification~(\ref{tensormap}), one observes that formula~(\ref{FD2}) is a reformulation of~(\ref{decompose}), that has already been derived for $\alpha=1$ in previous work by the first author with J.M.Portegies~\cite[Thm.2.3 \& Eq.31]{portegies_new_2017}. The key idea behind the derivation, the expansion and the completeness of the eigenfunctions $\{\Phi^{l,m}_{\bomega}\}$ is summarized in Remark~\ref{rem:completebasis}. The general case $\alpha \in (0,1]$ then directly follows by Subsection~\ref{rem:3}.
\end{proof}
\iftoggle{insertremarks}
{
\begin{Remark}
The advantage of doing the spatial decomposition in the Fourier domain of $G/H$ as done in~(\ref{FD1}) over its spatial counterpart~(\ref{FD2}) is that it involves 2 integrations less and involves ordinary spherical harmonics (rather than generalized spherical harmonics).
\end{Remark}
}
{}
Recently, exact formulas for the (hypo-elliptic) heat-kernels on $G=SE(3)$ and on $G/H=\R^{3} \rtimes S^{2}$ (i.e. the case $\alpha=1$) have been published in 
\cite{portegies_new_2017}. In the next theorem we
\begin{itemize}
\item
extend these results to the kernels of PDE~(\ref{eq:differentialQ}), which are Forward Kolmogorov equations of $\alpha$-stable L\'{e}vy process with $\alpha \in (0,1]$,
\item
provide a structured alternative formula via the transform $\mathcal{F}_{G/H}$ characterized in Theorem~\ref{corr:1}.
\end{itemize}
\begin{Theorem}\label{th:three}
We have the following formulas for the probability kernels of $\alpha$-stable L\'{e}vy processes on $\R^{3}\rtimes S^2$:
\begin{itemize}
\item
Via conjugation with $\mathcal{F}_{\R^{3}\rtimes S^2}$:
\begin{equation} \label{formule1}
\boxed{
K_{t}^{\alpha}(\mathbf{x},\mathbf{n})= \frac{1}{(2\pi)^2} \int \limits_{0}^{\infty} \sum \limits_{s \in \mathbb{Z}}
\sum \limits_{l=|s|}^{\infty} e^{-(- \lambda^{l,s}_p)^{\alpha} t}\;
\left[\overline{\sigma}^{p,s}_{(\mathbf{x},\mathbf{n})}\right]_{l,0,l,0}\; p^{2} {\rm d}p.
}
\end{equation}
where {\small
$\left[\overline{\sigma}^{p,s}_{(\mathbf{x},\mathbf{n})}\right]_{l,0,l,0}=\left\langle \sigma_{(\mathbf{x},\mathbf{R}_{\mathbf{n}})}^{p,s} \Phi_{p \mathbf{a}}^{l,s}(p^{-1}\cdot)\, , \, \Phi_{p \mathbf{a}}^{l,s}(p^{-1}\cdot)\right\rangle_{\mathbb{L}_{2}(pS^2)}$} can be derived analytically (see~\iftoggle{insertremarks}{Remark~\ref{rem:analytic}}{\cite[Rem.~\!18]{arxivFT}}).
\item Via conjugation with $\overline{\mathcal{F}}_{\R^{3}}$:
\begin{equation} \label{formule2}
\boxed{
K_{t}^{\alpha}(\mathbf{x},\mathbf{n})= \frac{1}{(2\pi)^3} \int \limits_{\R^{3}} \left(\sum \limits_{l=0}^{\infty} \sum \limits_{m=-l}^{l}  e^{-(-\lambda^{l,m}_{\|\bomega\|})^{\alpha}t}\;
\overline{\Phi_{\bomega}^{l,m}(\mathbf{a})} \, \Phi_{\bomega}^{l,m}(\mathbf{n})
\right) \; e^{i \mathbf{x} \cdot \bomega}\; {\rm d}\bomega.
}
\end{equation}
\end{itemize}
\end{Theorem}
\begin{proof} Formula (\ref{formule1}) follows by
\[
K_{t}^{\alpha}(\mathbf{x},\mathbf{n})=(e^{t Q_{\alpha}} \delta_{(\mathbf{0},\mathbf{a})})(\mathbf{x},\mathbf{n})=
\left(\mathcal{F}_{G/H}^{-1} \circ e^{ t \mathcal{F}_{G/H} \circ Q_{\alpha} \circ \mathcal{F}_{G/H}^{-1}} \circ \mathcal{F}_{G/H} \delta_{(\mathbf{0},\mathbf{a})} \right)(\mathbf{x},\mathbf{n}).
\]
Now $(\mathcal{F}_{G/H} \delta_{(\mathbf{0},\mathbf{a})})(\sigma^{p,s})=\textrm{1}_{\mathbb{L}_{2}(pS^2)}$ implies
$((\mathcal{F}_{G/H} \delta) (\sigma^{p,s})_{(\mathbf{0},\mathbf{a})})(\sigma^{p,s}))_{l,0,l',0}=\delta_{ll'}$ so that the result follows by setting $U=\delta_{(\mathbf{0},\mathbf{a})}$ (or more precisely, by taking $U$ a sequence that is a bounded approximation of the unity centered around $(\mathbf{0},\mathbf{a})$) in Theorem~\ref{th:decomposition}, where we recall the inversion formula from the first part of Theorem~\ref{corr:1}. 

Formula (\ref{formule2}) follows similarly by
\[
K_{t}^{\alpha}(\mathbf{x},\mathbf{n})=\left(e^{t Q_{\alpha}} \delta_{(\mathbf{0},\mathbf{a})}\right)(\mathbf{x},\mathbf{n})=
\left(\overline{\mathcal{F}}_{\R^3}^{-1} \circ e^{ t \overline{\mathcal{F}}_{\R^3} \circ Q_{\alpha} \circ \overline{\mathcal{F}}_{\R^3}^{-1}} \circ \overline{\mathcal{F}}_{\R^3} \delta_{(\mathbf{0},\mathbf{a})} \right)(\mathbf{x},\mathbf{n}).
\]
Now $\left(\overline{\mathcal{F}}_{\R^3} \delta_{(\mathbf{0},\mathbf{a})}\right)(\sigma^{p,s})=\frac{1}{(2\pi)^{\frac{3}{2}}} \delta_{\mathbf{a}}$
and the result follows from the second part of Theorem~\ref{corr:1} (again by taking $U$ a sequence that is a bounded approximation of the unity centered around $(\mathbf{0},\mathbf{a})$).
\end{proof}

\iftoggle{insertremarks}
{
\begin{Remark}\label{rem:analytic}
The coefficients $\left(\sigma^{p,s}_{(\mathbf{x},\mathbf{R}_{\mathbf{n}})}\right)_{l,0,l,0}$ relative to the basis functions $\{Y^{l,s}\}$, recall (\ref{matrixcoeff}), can be found in \cite[eq.10.35]{chirikjian_engineering_2000}
and by the basis transform (\ref{Philm}) we find
\[
\left[\sigma^{p,s}_{(\mathbf{x},\mathbf{R}_{\mathbf{n}})}\right]_{l,0,l,0}= \sum \limits_{l_1,l_2 \geq |s|} \frac{d_{l_1-s}^{l,0}(p) \; d_{l_2-s}^{l,0}(p)}{\|\mathbf{d}^{l,0}(p)\|^2} \left(\sigma^{p,s}_{(\mathbf{x},\mathbf{R}_{\mathbf{n}})}\right)_{l_1,0,l_2,0}\ .
\]
\end{Remark}
}{
}

\subsection{Monte-Carlo Approximations of the Kernels \label{ch:montecarlo}}

A stochastic approximation for the kernel $K_{t}^{\alpha}$ is computed by binning the endpoints of discrete random walks simulating $\alpha$-stable processes on the quotient $\R^{3}\rtimes S^2$ that we will explain next.
Let us first consider the case $\alpha=1$. For $M \in \mathbb{N}$ fixed, we have the discretization
\begin{equation}\label{eq:RandomWalkAlpha1}
\begin{cases}
\mathbf{X}_M = \mathbf{X}_0 + \sum\limits_{k=1}^{M} \sqrt{\frac{t D_{33}}{M}} \epsilon_k \mathbf{N}_{k-1},\\
\mathbf{N}_M = \left(\prod\limits_{k=1}^{M}  \mathbf{R}_{\mathbf{a}, \gamma_k} \mathbf{R}_{\mathbf{e}_y, \beta_k \sqrt{\frac{t D_{44}}{M}}} \right)  \mathbf{N}_0=
\left(\mathbf{R}_{\mathbf{a}, \gamma_M} \mathbf{R}_{\mathbf{e}_y, \beta_M \sqrt{\frac{t D_{44}}{M}}} \circ \ldots \circ \mathbf{R}_{\mathbf{a}, \gamma_1} \mathbf{R}_{\mathbf{e}_y, \beta_1 \sqrt{\frac{t D_{44}}{M}}}\right)\; \mathbf{N}_0\ ,
\end{cases}
\end{equation}
with $\epsilon_{k} \sim  G_{t=1}^{\R} \sim \mathcal{N}(0,\sigma = \sqrt{2})$ stochastically independent Gaussian distributed on $\R$ with $t=1$; with uniformly distributed $\gamma_{k}\sim
\textrm{Unif}\left(\mathbb{R}/(2\pi \mathbb{Z}) \equiv [-\pi,\pi)\right)$;
and $\beta_{k} \sim g$, where $g: \R \to \R^+$ equals $g(r)= \frac{|r|}{2}\; e^{-\frac{r^2}{4}}$ in view of the theory of isotropic stochastic processes on Riemannian manifolds by Pinsky \cite{pinsky_isotropic_1976}.
By the central limit theorem for independently distributed variables \emph{with finite variance} it is only the variances of the distributions for the random variables $g$ and $G_{t=1}^{\R}$ that matter. One may also take
\[
\begin{array}{l}
\epsilon_{k} \sim \sqrt{3}\; \textrm{Unif}\left[-\frac{1}{2},\frac{1}{2}\right] \textrm{ and }
\beta_{k} \sim \sqrt{6}\; \textrm{Unif}\left[-\frac{1}{2},\frac{1}{2}\right] \textrm{ or }
\epsilon_{k} \sim G_{t=1}^{\R} \textrm{ and }
\beta_{k} \sim G_{t=2}^{\R}.
\end{array}
\]
These processes are implemented recursively, for technical details and background see Appendix~\ref{ch:PT}.

\begin{Proposition}\label{prop:RW}
The discretization~(\ref{eq:RandomWalkAlpha1}) can be re-expressed, up to order $\frac{1}{M}$ for $M \gg 0$, as follows:
\begin{equation}\label{eq:62alt}
(\mathbf{X}_M,\mathbf{N}_M) \sim \mathbf{G}_{M} \odot (\mathbf{0},\mathbf{a}), \textrm{ with }
\mathbf{G}_{M}= \left(\prod \limits_{k=1}^M e^{\;\;\sum \limits_{i=3}^{5} \sqrt{\frac{t \, D_{ii}}{M}} \epsilon^{i}_{k} A_{i}}\right) \mathbf{G}_0,
\end{equation}
with $\epsilon^{i}_{k}\sim G_{t=1}^{\R}$ stochastically independent normally distributed variables with $t=\frac{1}{2}\sigma^2=1$, and $D_{44}=D_{55}$.
\end{Proposition}
\begin{proof}
In our construction, $\beta_k$ and $\gamma_k$ can be seen as the polar radius and the polar angle (on a periodic square $[-\pi, \pi] \times [-\pi,\pi])$ of a Gaussian process with $t=1$ on a plane spanned by rotational generators $A_4$ and $A_5$ . The key ingredient to obtain~(\ref{eq:62alt}) from~(\ref{eq:RandomWalkAlpha1}) is given by the following relation:
\begin{equation}\label{eq:firstosec}
\displaystyle e^{u \cos v A_5 - u \sin v A_4} = e^{v A_6} e^{u A_5} e^{-v A_6}, \textrm{ for all } u, v \in \R,
\end{equation}
which we use for $u = \beta_k \sqrt{\frac{t D_{44}}{M}}$ and $v = \gamma_k \sqrt{\frac{t D_{44}}{M}}$.

The second ingredient is given by the Campbell-Baker-Hausdorff-Dynkin formula:
$$\textrm{for all $a_i = O(\frac{1}{\sqrt{M}})$ and for $M$ large, we have } e^{a_3 A_3} e^{a_4 A_4} e^{a_5 A_5} = e^{(a_3 A_3 + a_4 A_4 + a_5 A_5)(1 + O(\frac{1}{M}))}, $$
that allows to decompose the stochastic process in $SE(3)$ into its spatial and angular parts.
\end{proof}

\begin{figure}[t!]
   \centering
   \includegraphics[width=\textwidth]{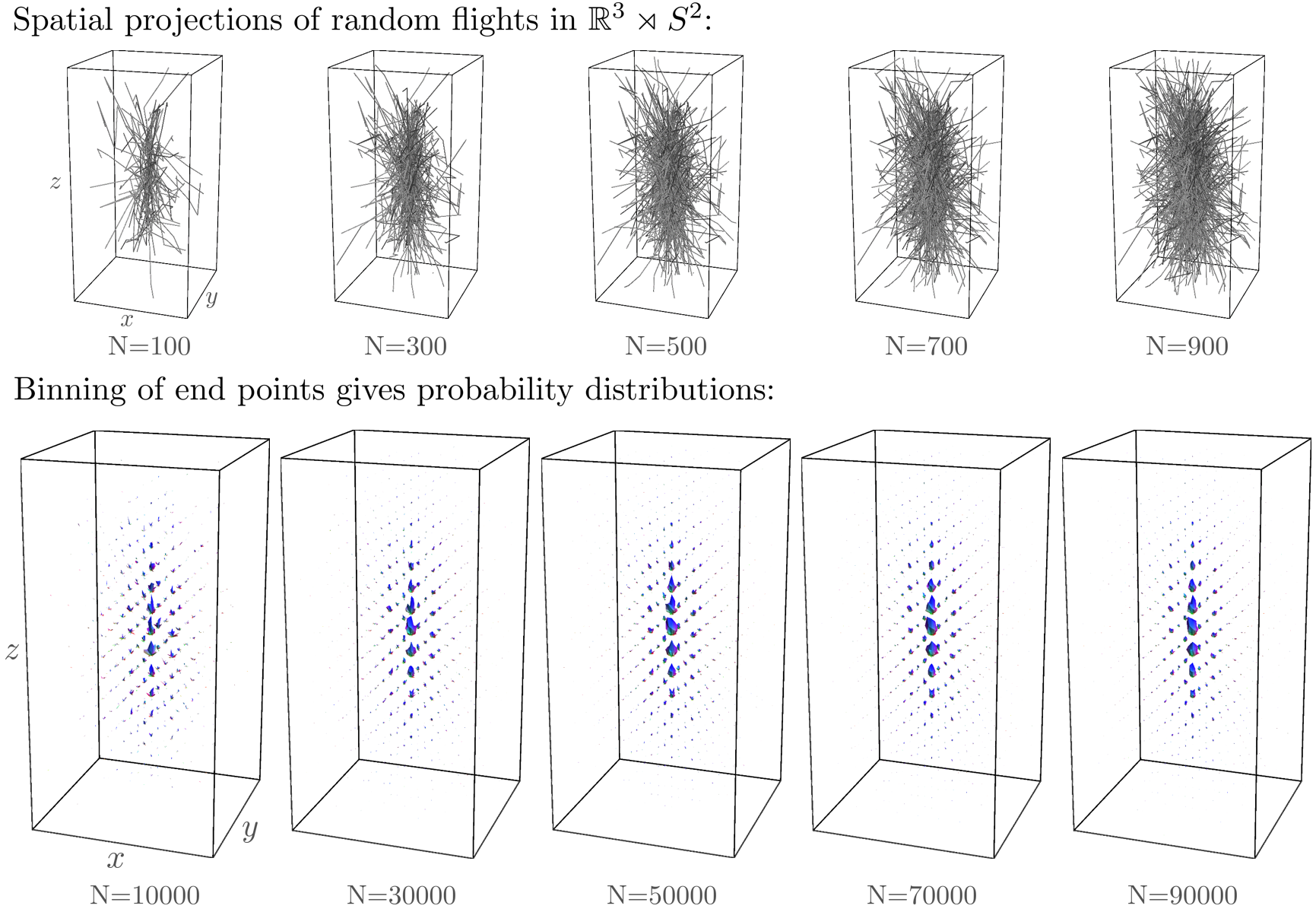}
 \caption{Top: Spatial projections in $\R^3$ of $N$ sample paths of the discrete random walks (or rather `drunk man's flights') in $\mathbb{R}^{3}\rtimes S^{2}$ for $\alpha=1$, given by (\ref{eq:RandomWalkAlpha1}), for increasing $N$ (with $\sigma=\frac{4\pi}{252}$, $\Delta s=1$, $M=40$).
 Bottom:  Convergence of the Monte-Carlo simulation kernel (\ref{MCS}) for $\alpha=1$ and $N\to \infty$.
 As $N$ increases the Monte-Carlo simulation converges towards the exact solution.
 For a comparison of the exact diffusion kernel (\ref{formule1}) and its Monte-Carlo approximation (\ref{MCS}), see Fig.~\!\ref{fig:comparison}. }\label{fig:figureMonteCarlo}
 \end{figure}

For the binning we divide $\mathbb{R}^3$ into cubes $c_{ijk}$, $i,j,k \in \mathbb{Z}$, of size $\Delta s \times \Delta s \times \Delta s$:
\begin{equation}\label{cubes}
c_{ijk} := \left[(i-\frac12)\Delta s,(i+\frac12)\Delta s\right] \times \left[(j-\frac12)\Delta s,(j+\frac12)\Delta s\right] \times \left[(k-\frac12)\Delta s,(k+\frac12)\Delta s\right].
\end{equation}
We divide $S^2$ into bins $B_l$, $l = \{1,\dots,b\}$ for $b \in \mathbb{N}$, with surface area $\sigma_{B_l}$ and maximal surface area $\sigma_{B}$. The number of random walks in a simulation with traveling time $t$ that have their end point $\bx_M \in c_{ijk}$ with their orientation $\bn_M \in B_l$ is denoted with $\#_t^{ijkl}$. Furthermore, we define the indicator function
\begin{equation*}
1_{c_{ijk},B_l}(\bx,\bn) := \begin{cases}
1 & \bx \in c_{ijk}, \bn \in B_l,\\
0 & \text{otherwise}.
\end{cases}
\end{equation*}
When the number of paths $N \rightarrow \infty$, the number of steps in each path $M \rightarrow \infty$ and the bin sizes tend to zero, the obtained distribution converges to the exact kernel:
\begin{equation} \label{MCS}
\boxed{
\begin{array}{l}
 \lim \limits_{N \rightarrow \infty}
\lim \limits_{\Delta s, \sigma_B \rightarrow 0}\lim \limits_{M \rightarrow \infty} p_t^{\Delta s, \sigma_B, N, M}(\bx,\bn) = K_t^{\alpha=1}(\bx,\bn), \\
\textrm{ with }p_t^{\Delta s, \sigma_B, N, M}(\bx,\bn) = \sum \limits_{l = 1}^b \sum \limits_{i,j,k \in \mathbb{Z}} 1_{c_{i,j,k},B_l}(\bx,\bn) \frac{\#_t^{ijkl}}{M (\Delta s)^3 \sigma_{B_l}}.
\end{array}
}
\end{equation}
The convergence is illustrated in Figure~\ref{fig:figureMonteCarlo}.
\subsubsection{Monte-Carlo simulation for $\alpha \in (0,1]$.}
\begin{figure}[t!]
   \centering
   \includegraphics[width=0.82\textwidth]{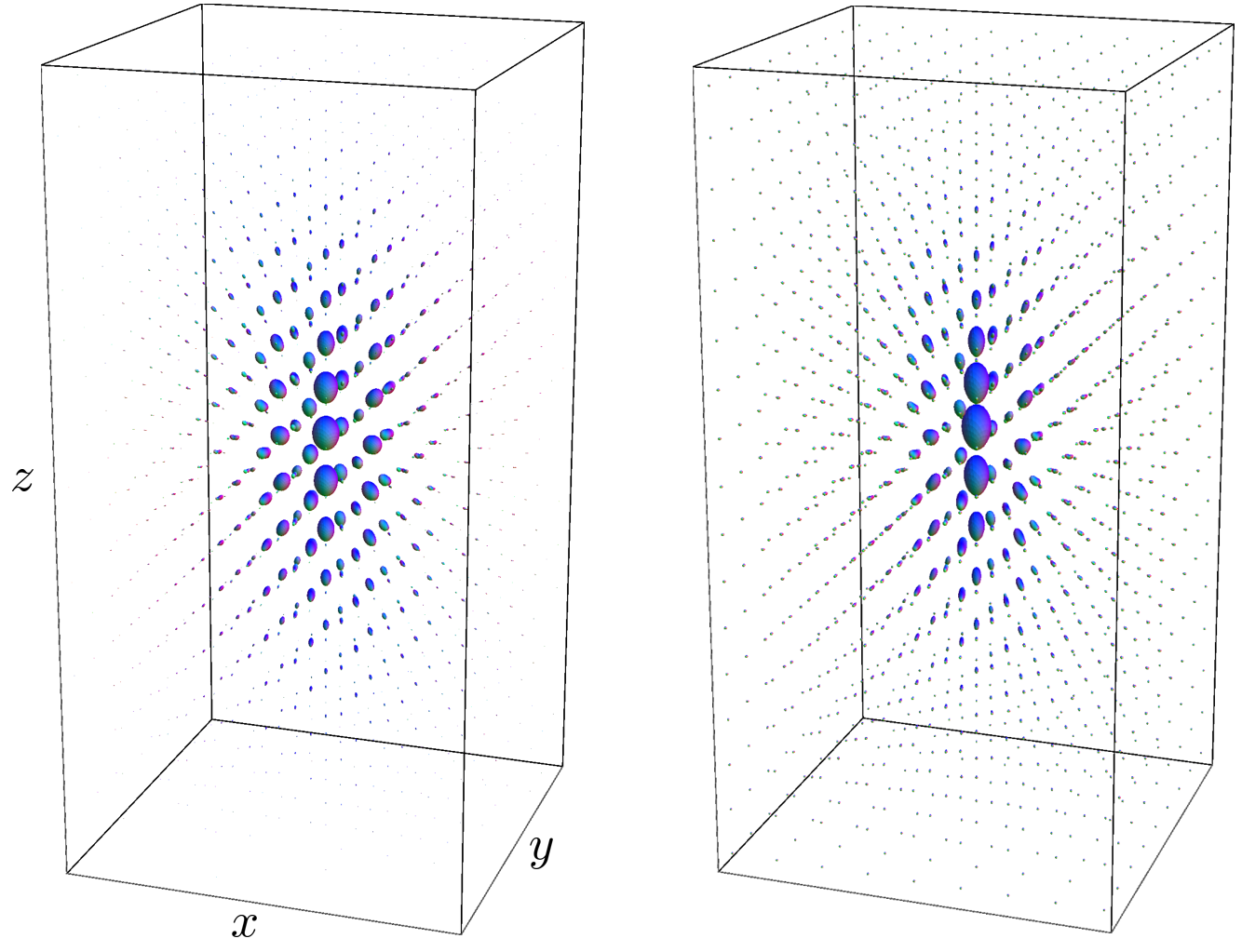}
 \caption{Left: The hypo-elliptic diffusion kernel (Eq.~\!(\ref{formule1}) for $\alpha=1$ and $t=2$). Right:
 the hypo-elliptic Poisson kernel (Eq.~\!(\ref{formule1}) for $\alpha=\frac{1}{2}$ and $t=3.5$). Parameters settings: $D_{44}=0.2,\! D_{33}=1,\!D_{11}=0$.
 \label{fig:alpha}}
 \end{figure}

\noindent Let $q_{t,\alpha}:\mathbb{R}^{+} \to \mathbb{R}^+$ be the temporal probability density given by the inverse Laplace transform
\begin{equation}\label{eq:qker}
\begin{array}{l}
q_{t,\alpha}(\tau) = \mathcal{L}^{-1}\left(\lambda \to e^{-t \lambda^\alpha}\right)(\tau), \textrm{ with in particular: } \\[6pt]
\qquad \textrm{for }\alpha=\frac{1}{2} \textrm{ it is } q_{t,\frac12}(\tau) = \frac{t}{2 \tau \sqrt{\pi \tau}} e^{-\frac{t^2}{4 \tau}}, \\
\qquad \textrm{for }\alpha \uparrow 1 \textrm{ we find } q_{t,\alpha}(\cdot) \to \delta_{t} \textrm{ in distributional sense .}
\end{array}
\end{equation}
For explicit formulas in the general case $\alpha \in (0,1]$ see~\cite{yosida_functional_1980}.
Then one can deduce from Theorem~\ref{th:three} that
\begin{equation} \label{kernelrel}
K_t^{\alpha}(\mathbf{x},\mathbf{n}) = \int\limits_0^\infty q_{t,\alpha}(\tau) \; K_{\tau}^{\alpha=1}(\mathbf{x},\mathbf{n}) \; {\rm d} \tau.
\end{equation}
This allows us to directly use the Monte-Carlo simulations for the diffusion kernel $\alpha=1$ for several time instances to compute a Monte-Carlo simulation of the $\alpha$-stable L\'{e}vy kernels for $\alpha \in (0,1]$. To this end we replace the Monte Carlo approximation (\ref{MCS}) for $\alpha=1$ in the above
formula (\ref{kernelrel}). See Figure~\ref{fig:alpha}, where we compare the diffusion kernel $K_{t}^{\alpha=1}$ to the Poisson kernel $K_{t}^{\alpha=\frac{1}{2}}$. See also Appendix~\ref{app:A21}.

\subsection{Comparison of Monte-Carlo Approximations of the Kernels to the Exact Solutions}

In this section we compute the probability density kernels $K_t^\alpha$ via the analytic approach of Subsection \ref{sec:Expl} (Eq.~\ref{formule2}, Thm.~\ref{th:three}) and via the Monte-Carlo approximation of Subsection~\ref{ch:montecarlo}. The kernels are computed on a regular grid with each $(x_i,y_j,z_k)$ at the center of the cubes $c_{ijk}$ of (\ref{cubes}) with $i, j = -3,\dots,3$, $k = -5,\dots,5$, and $\Delta s = 0.5$.
The Monte-Carlo simulations also require spherical sampling which we did by a geodesic polyhedron that sub-divides each mesh triangle of an icosahedron into $n^2$ new triangles and projects the vertex points to the sphere. We set $n=4$ to obtain $252$ (almost) uniformly sampled points on $S^2$.

The exact solution is computed using (truncated) spherical harmonics with $l \leq 12$. To obtain the kernel we first solve the solution in the spatial Fourier domain and then do an inverse spatial Fast Fourier Transform. The resulting kernel $K_t^\alpha$ (where we literally follow~(\ref{formule2})) is only spatially sampled and provides for each $(x_i,y_j,z_k)$ an analytic spherical distribution expressed in spherical harmonics. 

For the Monte-Carlo approximation we follow the procedure as described in Subsection~\ref{ch:montecarlo}. The kernel $K_t^\alpha$ is obtained by binning the end points of random paths on the quotient $\mathbb{R}^3 \rtimes S^2$ (cf. Eq.~(\ref{eq:RandomWalkAlpha1})) and thereby approximate the limit in Eq.~(\ref{MCS}). Each path is discretized with $M=40$ steps and in total $N=10^{10}$ random paths were generated. The sphere $S^2$ is divided into $252$ bins with an average surface area of $\sigma_{B_l} \approx \frac{4 \pi}{252}$.

\begin{figure}
\centerline{
\includegraphics[width=0.8\textwidth]{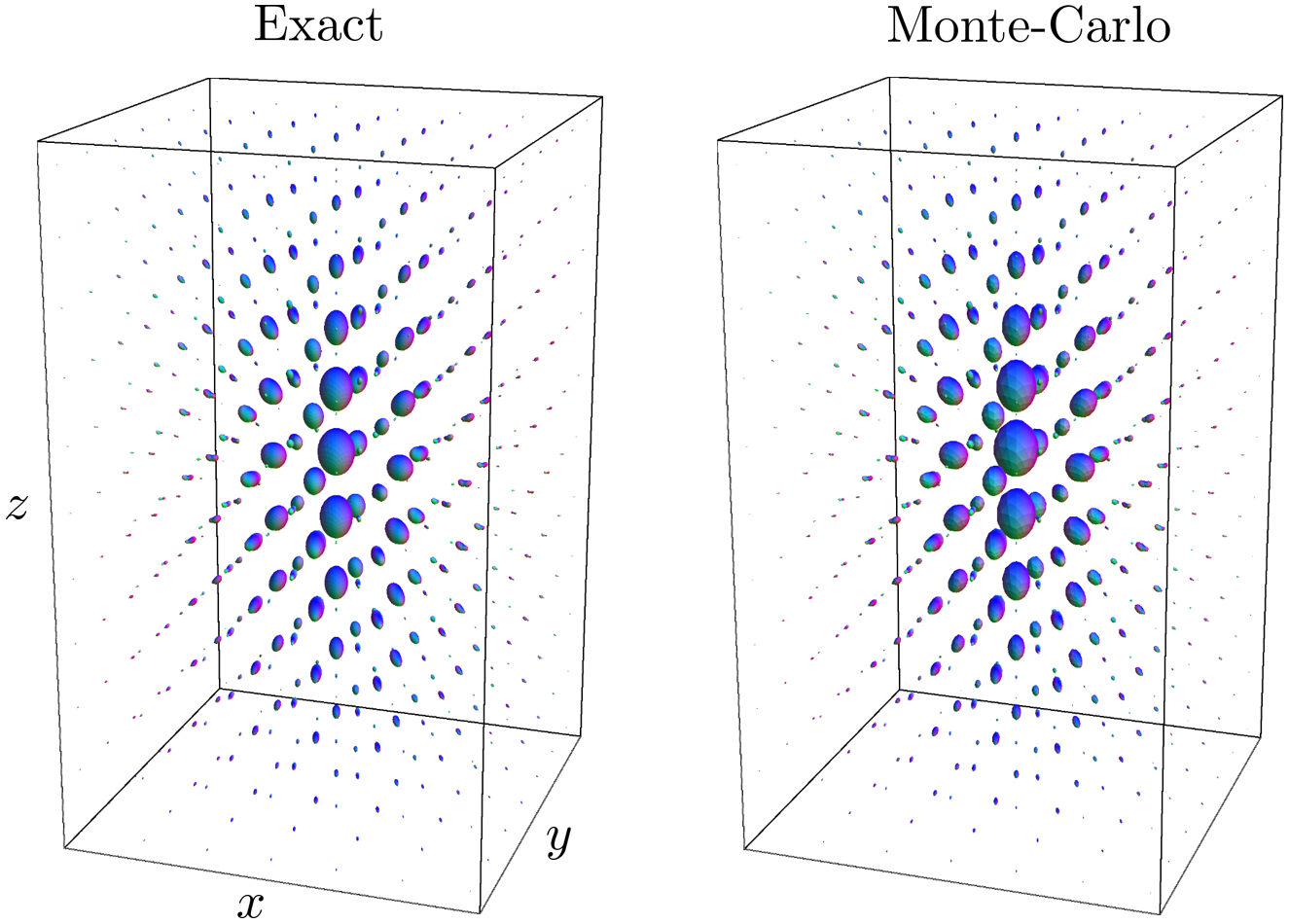}
}
\caption{The exact kernel $K_t^\alpha$ and its Monte-Carlo approximation for $t=2$, $\alpha=1$, $D_{33}=1$, $D_{44}=0.2$. \label{fig:comparison}}
\end{figure}

In Figures \ref{fig:sample-paths}, \ref{fig:figureMonteCarlo}, \ref{fig:alpha} and \ref{fig:comparison} we set $D_{33}=1$, $D_{44}=0.2$. In the comparison between the kernels $K_t^{\alpha=1}$ with $K_t^{\alpha=0.5}$ we set $t=2$ and $t=3.5$ respectively in order to match the full width at half maximum value of respectively the distributions. In Figures \ref{fig:sample-paths}, \ref{fig:figureMonteCarlo} and \ref{fig:comparison} we set $\alpha=1$ and $t=2$.
In Figures \ref{fig:sample-paths}, \ref{fig:figureMonteCarlo}, \ref{fig:alpha} we sample the grid~(\ref{cubes}) with $|i|,|j| \leq 4$, $|k|\leq 8$.

Figure~\ref{fig:comparison} shows that the Monte-Carlo kernel closely approximates the exact solution and since the exact solutions can be computed at arbitrary spherical resolution, it provides a reliable way to validate numerical methods for $\alpha$-stable L\'{e}vy processes on $\mathbb{R}^3 \rtimes S^2$.

\newpage
\section{Conclusions}\label{ch:coclusion}

We have set up a Fourier transform $\mathcal{F}_{G/H}$ on the homogeneous space of positions and orientations.
The considered Fourier transform acts on functions that are bi-invariant with respect to the action of subgroup $H$.
We provide explicit formulas (relative to a basis of modified spherical harmonics) for the transform,
its inverse, and its Plancherel formula, in Theorem~\ref{corr:1}.

Then we used this Fourier transform to derive new exact solutions to the probability kernels of $\alpha$-stable L\'{e}vy processes on $G/H$, including the diffusion PDE for Wiener processes which is the special case $\alpha=1$. They are obtained by spectral decomposition of the evolution operator in Theorem~\ref{th:decomposition}. 

New formulas for the probability kernels are presented in Theorem~\ref{th:three}.
There the general case $0<\alpha<1$ followed from the case $\alpha=1$ by taking the fractional power of the eigenvalues.
In comparison to previous formulas in~\cite{portegies_new_2017} for the special case $\alpha=1$ obtained via a spatial Fourier transform, we have more concise formulas with a more structured evolution operator in the Fourier domain of $G/H$, where we rely on ordinary spherical harmonics, and where we reduced the dimension of the manifold over which is integrated from 3 to 1 (as can be seen in Theorem~\ref{th:three}).

We introduced stochastic differential equations (or rather stochastic integral equations) for the $\alpha$-stable L\'{e}vy processes in Appendix~\ref{app:A1}, and we provide
simple discrete approximations where we rely on matrix exponentials in the Lie group $SE(3)$ in Proposition~\ref{prop:RW}.

We verified the exact solutions and the stochastic process formulations, by Monte-Carlo simulations that confirm to give the same kernels, as is shown in Figure~\ref{fig:comparison}.
We also observed the expected behavior that the probability kernels for $0<\alpha<1$ have heavier tails, as shown in Figure~\ref{fig:alpha}.

The PDEs and the probability kernels have a wide variety of applications in image analysis (crossing-preserving, contextual enhancement of diffusion-weighted MRI, cf.~\cite{duits_left-invariant_2010-3,reisert_fiber_2011,meesters_fast_2016,momayyez-siahkal_3d_2009,prckovska_contextual_2015,portegies_improving_2015,meesters_s.p.l._cleaning_2016}
or in crossing-preserving diffusions in 3D scalar images \cite{Janssen2018}),
and in robotics \cite{chirikjian_engineering_2000,chirikjian_stochastic_2011,Barbaresco} and in probability theory \cite{Feller,Benoist}.
The generalizations to $\alpha \in (0,1]$ allow for longer range interactions between local orientations (due to the heavy tails). This is also of  interest in machine learning, where convolutional neural networks on the homogeneous space of positions and orientations \cite{BekkersMICCAI,Wellink} can be extended to 3D~\!\cite{worrall2018cubenet,winkels20183d}, which may benefit from the PDE descriptors and the Fourier transform presented here. \\
\
\\





\acknowledgments{ \ \\
We gratefully acknowledge former PhD student J.M.~Portegies (ASML, The Netherlands) for providing us with the \emph{Mathematica} code for the exact solutions and Monte-Carlo simulations for the diffusion case $\alpha=1$ that we simplified and generalized to the general case $\alpha \in (0,1]$.

The research leading to the results of this paper has received funding from the European Research Council under the European Community's Seventh Framework Programme (FP7/2007-2013) / ERC grant \emph{Lie Analysis}, agr.~nr.~335555.
}

\authorcontributions{
R.~Duits led the project, and wrote the main body/general theory of this article. This was done in a close and fruitful collaboration with A.~Mashtakov (resulting in the final theoretical formulations and the final structure of this article) and with E.J.~Bekkers
(resulting in the experiments, simulations and discrete stochastic process formulations in the article).

}
\conflictsofinterest{``The authors declare no conflict of interest.'' }

\newpage
\abbreviations{The following abbreviations and symbols are used in this manuscript:\\
\noindent
\begin{tabular}{@{}lll}
UIR & Unitary Irreducible Representation & \ \; \\
$G$ & The rigid body motions group $SE(3)$ & Eq.~(\ref{product}) \\
$\ba$ & The reference axis $\ba = \be_z = (0,0,1)^T$ & Eq.~(\ref{quot})\\
$H$ & The subgroup that stabilizes $(\mathbf{0}, \ba)$ & Eq.~(\ref{quot}) \\
$G/H$ & The homogeneous space of positions and orientations $\R^{3}\rtimes S^2$ & Eq.~(\ref{quot}) \\
$\overline{U}$ & The spatial Fourier transform of $U$ & Eq.~(\ref{eq:spatFT0}) \\
$\hat{U}$ & The Fourier transform $\hat{U} = \mathcal{F}_{G/H}U$ & Eq.~(\ref{FThom}) \\
$\alpha$ & Parameter of the $\alpha$-Stable processes (indexing fractional power of the generator) & Eq.~\!(\ref{PDEgroup})\\
$\overline{\alpha}$ & Rotation angle around reference axis $\mathbf{a}=\mathbf{e}_{z}=(0,0,1)$ & Remark~\ref{rem:notation} \\
$\sigma^{p,s}$ & UIR of $G = SE(3)$ & Eq.~(\ref{inter}) \\
$\overline{\sigma}^{p,s}$  & the action on the quotient corresponding to $\sigma^{p,s}$ & Definition~\ref{def:quotrep} \\
$\tilde{K}_t^{\alpha}$ & The probability kernel on $G$ & Eq.~(\ref{solconv})\\
$K_t^{\alpha}$ & The probability kernel on $G/H$ & Eq.~(\ref{eq:16b})\\
$\tilde{W}_{\alpha}$ & Solution of the PDE on $G$ & Eq.~(\ref{PDEgroup}) \\
$W_{\alpha}$  & Solution of the PDE on $G/H$ & Eq.~(\ref{eq:differentialQ}) \\
$\tilde{Q}_{\alpha}$  & Evolution generator of the PDE on $G$ & Eq.~(\ref{eq:Qtildealpha}) \\
$Q_{\alpha}$  & Evolution generator of the PDE on $G/H$ & Eq.~(\ref{Qalphafund}) \\
$\bR_{\bn}$ & Any rotation that maps $\ba$ onto $\bn$ & Remark~\ref{rem:terminology} \\
$\bR_{\mathbf{v},\phi}$  & A counter-clockwise rotation about axis $\mathbf{v}$ with angle $\phi$& Remark~\ref{rem:terminology} \\
$\bP_t$ & L\'{e}vy Processes on $G/H$ & Definition~\ref{def:levyprocquot}\\
$\overline{\bP}_t$  & L\'{e}vy Processes on $\mathbb{R}^3 \times \mathbb{R}^3$ & Eq.~(\ref{eq:SDEembd}) \\
$q_{t,\alpha}$ & The kernel relating $K_t^{\alpha}$ and $K_t^{1}$ & Eq.~(\ref{eq:qker}) \\
$Y^{l,m}$  & The ordinary spherical harmonics & Proposition~\ref{prop} \\
$Y^{l,m}_s$ & The modified spherical harmonics according to~\cite{chirikjian_engineering_2000} & Proposition~\ref{prop} \\
$Y^{l,m}_{\bomega}$  & The generalized spherical harmonics according to~\cite{portegies_new_2017} & Definition~\ref{def:1} \\
$\Phi^{l,m}_{\bomega}$ & The spheroidal wave basis function for $\mathbb{L}_{2}(S^2)$ & Definition~\ref{def:2} \\
$(\overline{\alpha}, \beta , \gamma)$ & ZYZ Euler angles. & (\ref{Eulerangles1}) 
\end{tabular}
}
\ \\
\appendixtitles{no} 
\appendixsections{multiple} 
\appendix

\section{Probability theory \label{ch:PT}}
%
%
%
%
%
%


\subsection{L\'{e}vy Processes on $\R^3 \rtimes S^2$}\label{app:A1}

In the next definition, we define L\'{e}vy processes on our manifold of interest $G/H=\quot$. Recall, that the action of $G = SE(3)$ on $G/H$ is given by~(\ref{actiononset}).
As a prerequisite, we define the `difference' of 2 random variables $\bP_{1} = (\bX_1, \bN_1)$ and $\bP_{2} = (\bX_2, \bN_2)$ in $\quot$:
\begin{equation}\label{eq:randdifquot}
\bG_2^{-1} \odot \bP_1 =
(\bX_{2},\bR_{\bN_{2}})^{-1} \odot (\bX_{1},\bN_{1})=
(\bR_{\bN_{2}}^{T}(\bX_{1}-\bX_{2}),\bR_{\bN_{2}}^{T} \bN_{1}),
\end{equation}
where we relate random variables on $G/H$ and in $G$ via $\bP = \bG \odot (\mathbf{0},\mathbf{a})$, according to ~(\ref{eq:qoutfromgroup}).

We will assume that $\mathbf{P}_{1}$ and $\mathbf{P}_{2}$ are chosen such that the distribution of $\bG_2^{-1} \odot \bP_1$ is invariant under the choice of rotation variable $\mathbf{R}_{\mathbf{N}_2} \in SO(3)$ which maps reference axis $\mathbf{a}$ onto $\mathbf{N}_2$.
This is done in view of the homogeneous space structure (\ref{quot}) and the fact that L\'{e}vy processes on Lie groups such as $G=SE(3)$ require Lie group inversion in their definition (see, e.g.~\cite{MIngLiao}).

\begin{Definition}\label{def:levyprocquot}
A stochastic process $\{\mathbf{P}_t \, : \, t \geq 0\}$ on $G/H$ is a L\'{e}vy process if the following conditions hold:
\begin{enumerate}
\item
For any $n\geq 1$ and $0\!\leq\! t_0<\!t_1\!<\!\ldots\!<\!t_n$, the variables $\bP_{t_0}$, $\bG_{t_0}^{-1} \odot \bP_{t_1}$, $\ldots$, $\bG_{t_{n-1}}^{-1} \odot \bP_{t_n}$ are independent.
\item The distribution of $\bG_{s}^{-1} \odot \bP_{s+t}$ does not depend on $s\geq0$.
\item $\mathbf{P}_0 = (\mathbf{0},\ba)$ almost surely.
\item It is stochastically continuous, i.e.
$\lim_{s\to t} P[d(\bP_s, \bP_t) > \eps] = 0$, $\forall \eps >0$.
\\
Here $d((\bx_1,\bn_1), (\bx_2,\bn_2)) = |\bx_1-\bx_2|^2 + \arccos^2(\bn_1 \cdot \bn_2)$. 
\end{enumerate}
\end{Definition}
Let us consider the solutions
\[
W_{\alpha}(\mathbf{x},\mathbf{n},t)=(K_{t}^{\alpha}* U)(\mathbf{x},\mathbf{n})
\]
of our linear PDEs of interest (\ref{eq:differentialQ}) for $\alpha \in (0,1]$ fixed.
 Let us consider the case where $U \sim \delta_{(\mathbf{0},\mathbf{a})}$, so that the solutions are the probability kernels $K_{t}^{\alpha}$ themselves. We consider the random variables $\mathbf{P}_t^{\alpha}$ such that their probability densities are given by
\begin{equation} \label{Pt}
P(\mathbf{P}_{t}^{\alpha}=(\mathbf{x},\mathbf{n}))=K_{t}^{\alpha}(\mathbf{x},\mathbf{n}) \textrm{ for all }t\geq 0, (\mathbf{x},\mathbf{n})\in \quot.
\end{equation}
\begin{Proposition}\label{prop:A1}
The stochastic process $\{\mathbf{P}_t^{\alpha}\, : \, t \geq 0\}$ is a  L\'{e}vy processes on $\quot$.
\end{Proposition}
\begin{proof}
First, address items 1, 2.  On $G = SE(3)$, one has for two stochastically independent variables:
$$P(\bG_1 \bG_2 = g) = \int\limits_{G} P(\bG_2 = h^{-1} g) P(\bG_1 = h) \; {\rm d} h.$$
In particular, for $\bG_1 = \bG_t \sim \tilde{K}_t^{\alpha}$ and $\bG_2 = \bG_s \sim \tilde{K}_s^{\alpha}$ we have
$$\bG_s \bG_t \sim \tilde{K}_t^{\alpha} * \tilde{K}_s^{\alpha} = \tilde{K}_{t+s}^{\alpha} \quad \textrm{   and   } \quad \bG_s^{-1} \bG_{t+s}  = \bG_t \sim \tilde{K}_t^{\alpha},$$
which is due to $e^{t \tilde{Q}_\alpha}\circ e^{s \tilde{Q}_\alpha} = e^{(t+s) \tilde{Q}_\alpha}$, recall~(\ref{eq:etQ}).
Similarly, on the quotient $G/H$ we have
$$ \bG_s^{-1} \odot \bP_{s+t} = \bP_t \sim K_t^{\alpha}.$$
Furthermore, the choice of $\bG_s$ such that $\bG_s \odot (\mathbf{0}, \ba) = (\mathbf{0}, \ba)$ does not matter, since
$$P((\mathbf{0},\bR_{\ba,\bar{\alpha}})^{-1} \bG_s^{-1} \odot \bP_{s+t} = (\bx,\bn)) = K_t^{\alpha} ((\mathbf{0}, \bR_{\ba, \bar{\alpha}}) \odot (\bx,\bn)) = K_t^{\alpha}(\bx,\bn),$$
recall Eq.~(\ref{nice}).
Item 3 is obvious since we have $\bP_0 = \delta_{(\mathbf{0}, \ba)}$.
Item 4 follows by strong continuity of the semigroup operators~\cite{yosida_functional_1980},\!~\cite[Thm. 2]{duits_axioms_2004}.
\end{proof}
\begin{Lemma}\label{lem:A1}
The kernels $K_{t}^{\alpha}$ are infinitely divisible, i.e.
\[
K_{t}^{\alpha} * K_{s}^{\alpha}=K_{t+s}^{\alpha} \quad \textrm{  for all }s,t \geq 0.
\]
\end{Lemma}
\begin{proof}
The infinite divisibility directly follows from Corollary~\ref{corr:new} and $\mathcal{F}_{G/H}(K_{t}^{\alpha} * K_{s}^{\alpha})=
\mathcal{F}_{G/H}(K_{t}^{\alpha}) \circ  \mathcal{F}_{G/H}(K_{t}^{\alpha})= \mathcal{F}_{G/H}(K_{t+s}^{\alpha})$ which is clear due to~(\ref{formule1}).
\end{proof}

\iftoggle{insertremarks}{
\begin{Remark}
In our convention~(\ref{eq:randdifquot}), the ordering in general matters. Also on Lie group $SE(3)$  we have stationary and independent left increments or right increments~\cite{Nielsen}. In our case, the kernels carry inverse symmetry~(\ref{sym2}) and~(\ref{sym22}), thus this process has both stationary independent left increments and right increments.
\end{Remark}
}{%
 }

\begin{Remark}
Recall that on $\R^n$ a L\'{e}vy process $\mathbf{X}_{t}$ is called $\alpha$-stable if \begin{equation} \label{astab}
a^{-\frac{1}{2\alpha}}\mathbf{X}_{at} \sim \mathbf{X}_t \quad \textrm{ for all } a>0.
\end{equation}
This convention and property applies to all $n \in \mathbb{N}$, cf.~\cite{Feller}. Next we will come to a generalization of $\alpha$-stability but then for the processes $\mathbf{P}_t$. Here an embedding of $\R^{3}\rtimes S^{2}$ into $\R^6=\R^{3} \times \R^3$ will be required in order to give a meaning to $\alpha$-stability and a scaling relation on $\mathbf{P}_t=(\mathbf{X}_t,\mathbf{N}_t)$ that is similar to (\ref{astab}).
\end{Remark}

\subsection{SDE formulation of $\alpha$-Stable L\'{e}vy Processes on $\R^3 \rtimes S^2$ \label{ch:SDE}}

Consider the L\'{e}vy processes $\{\mathbf{P}_t \, : \, t \geq 0\}$  on $\R^{3}\rtimes S^2$ given by (\ref{Pt}). They give rise to the Forward Kolmogorov PDEs
(\ref{eq:differentialQ}) in terms of their stochastic differential equation (SDE) according to the book of Hsu on Stochastic Analysis on Manifolds \cite{Hsu}.

We will apply \cite[Prop.1.2.4]{Hsu} on the embedding map $\Phi: \R^{3} \times \R^3 \to \R^3 \rtimes S^2$ given by
\[
\Phi: (\mathbf{x}, \overline{\mathbf{n}}) \mapsto  \Phi(\mathbf{x}, \overline{\mathbf{n}}) =\left(\mathbf{x}, \frac{\overline{\mathbf{n}}}{\|\overline{\mathbf{n}}\|}\right)  = (\mathbf{x}, \mathbf{n}).
\]
Note that $\Phi_* = \mathcal{D} \Phi  = \left(I, \frac{1}{\|\overline{\mathbf{n}}\|} \left(I - \frac{\overline{\mathbf{n}}}{\|\overline{\mathbf{n}}\|} \otimes \frac{\overline{\mathbf{n}}}{\|\overline{\mathbf{n}}\|} \right) \right)$. Here $I$ denotes the identity map on $\R^3$.

Let us first concentrate on $\alpha=1$. In this case, our PDE~(\ref{eq:differentialQ}) becomes a diffusion PDE that is the forward Kolmogorov equation of a Wiener process $\mathbf{P}_t = (\bX_t, \bN_t)$ on $\R^3 \rtimes S^2$. Next we relate this Wiener process to a Wiener process $(\mathbf{W}_t^{(1)}, \overline{\mathbf{W}}_t^{(2)})$ in the embedding space $\R^3 \times \R^3$. We will write down the stochastic differential equation (SDE) and show that (\ref{eq:RandomWalkAlpha1}) boils down to discretization of the stochastic integral (in   \^{I}to sense) solving the SDE.

Next, we define $\overline{\mathbf{P}}_t = (\mathbf{X}_t, \overline{\mathbf{N}}_t)$ by the SDE in the embedding space:
\begin{equation}\label{eq:SDEembd}
{\rm d} \overline{\mathbf{P}}_t = \overline{\gothic{s}}|_{\overline{\mathbf{P}}_t} \circ {\rm d} \mathbf{W}_t,
\end{equation}
where $\mathbf{W}_t = (\mathbf{W}_t^{(1)}, \overline{\mathbf{W}}_t^{(2)})$, with $\mathbf{W}_t^{(1)}$ and $\overline{\mathbf{W}}_t^{(2)}$ being Wiener processes in $\R^3$; and where
$$
\overline{\gothic{s}}|_{\overline{\mathbf{P}}} ({\rm d}\mathbf{x},\rm d \overline{\mathbf{n}}) =
\left(\begin{array}{c}
\overline{\gothic{s}}^{(1)}|_{\overline{\mathbf{P}}}\left({\rm d}\mathbf{x},\rm d \overline{\mathbf{n}}\right) \\
\overline{\gothic{s}}^{(2)}|_{\overline{\mathbf{P}}}\left({\rm d}\mathbf{x},\rm d \overline{\mathbf{n}}\right)
 \end{array}\right)
=
 \left(\begin{array}{c}
\sqrt{D_{33}} \;
 \frac{\overline{\mathbf{N}}}{\|\overline{\mathbf{N}}\|}\,
\; \left(\frac{\overline{\mathbf{N}}}{\|\overline{\mathbf{N}}\|} \cdot {\rm d}\mathbf{x}\right) \\
\sqrt{D_{44}}\; {\rm d}\overline{\mathbf{n}}
 \end{array}\right).
$$
Here index (1) stands for the spatial part and (2) stands for the angular part.

Now we define a corresponding  process on $\quot$:
\[
\mathbf{P}_t = \Phi(\overline{\mathbf{P}}_t).
\]
Then the SDE for $\mathbf{P}_t = (\mathbf{X}_t, \mathbf{N}_t)$ becomes (see \cite[Prop.1.2.4]{Hsu})
$$
{\rm d} \mathbf{P}_t  = {\rm d}\left(\Phi \circ \overline{\mathbf{P}}_t \right) \Leftrightarrow
\begin{cases}
{\rm d}\mathbf{X}_t= \left.\overline{\gothic{s}}^{(1)}\right|_{\overline{\mathbf{P}}_t} \circ {\rm d}\mathbf{W}^{(1)}_t, \\
{\rm d}\mathbf{N}_t= \mathbb{P}_{\langle \mathbf{N}_t \rangle^\perp} \left.\overline{\gothic{s}}^{(2)}\right|_{\overline{\mathbf{P}}_t} \circ {\rm d}\overline{\mathbf{W}}^{(2)}_t,
\end{cases}
$$
where $\mathbf{N}_t = \frac{\overline{\mathbf{N}}_t}{\|\overline{\mathbf{N}}_t\|}$; and where $\mathbb{P}_{\langle \mathbf{N}_t \rangle^\perp} =\left(I - \mathbf{N}_t \otimes \mathbf{N}_t \right)$ denotes the orthogonal projection to the tangent plane perpendicular to  $\mathbf{N}_t$.

Therefore, we have the following SDE on $\quot$:
\begin{equation} \label{ene}
\boxed{
\begin{cases}
{\rm d}\mathbf{X}_t= \sqrt{D_{33}}\, \mathbf{N}_t  (\mathbf{N}_t  \cdot {\rm d}\mathbf{W}^{(1)}_t), \\
{\rm d}\mathbf{N}_t= \sqrt{D_{44}}\, \mathbb{P}_{\langle \mathbf{N}_t \rangle^\perp} {\rm d}\overline{\mathbf{W}}^{(2)}_t
\end{cases}
}
\end{equation}
So, integrating the SDE, we obtain the following stochastic integral    (in \^{I}to form):
\begin{equation} \label{SI}
\begin{cases}
\mathbf{X}_t = \mathbf{X}_0 + \sqrt{D_{33}}\int\limits_0^t \mathbf{N}_s (\mathbf{N}_s  \cdot {\rm d}\mathbf{W}^{(1)}_s) = \mathbf{X}_0 + \sqrt{D_{33}} \mslim\limits_{M \to \infty} \sum\limits_{k=1}^{M} \mathbf{N}_{t_{k-1}} \left(\mathbf{N}_{t_{k-1}} \cdot \left(\mathbf{W}^{(1)}_{t_k} - \mathbf{W}^{(1)}_{t_{k-1}} \right) \right),\\
\mathbf{N}_t =  \mslim\limits_{M\to \infty} \prod\limits_{k=1}^{M} \exp_{S^2} \left(\sqrt{D_{44}} \left(I - \mathbf{N}_{t_{k-1}} \otimes \mathbf{N}_{t_{k-1}}\right) \left(\overline{\mathbf{W}}^{(2)}_{t_k} - \overline{\mathbf{W}}^{(2)}_{t_{k-1}}\right)\right) \mathbf{N}_0. 
\end{cases}
\end{equation}
Here $\exp_{S^2}(V) \mathbf{n}_0$ denotes the exponential map on a sphere, i.e. its value is the end point (for $t=1$) of a geodesic starting from $\mathbf{n}_0 \in S^2$ with the tangent vector $V \in T_{\mathbf{n}_0}S^2$. Note that in the formula above, the symbol $\prod$ denotes the composition $$ \prod\limits_{k=1}^M \exp_{S^2}(V_k) \mathbf{n}_0 = \left(\exp_{S^2}(V_M) \circ \ldots \circ \exp_{S^2}(V_1)\right) \mathbf{n}_0.$$
Note that $\sqrt{D_{33}} \left(\bW_{t_k}^{(1)} -  \bW_{t_{k-1}}^{(1)}\right) = \sqrt{D_{33}} \bW_{t_k-t_{k-1}}^{(1)} = \sqrt{\frac{t D_{33}}{M}} \epsilon_k$, where $\epsilon_k \sim \bW^{(1)}_{1}$, i.e. $\epsilon_k \sim G_{t=1}$.

For $M \in \mathbb{N}$ fixed, we propose a discrete approximation for the stochastic integrals in (\ref{SI}):
\begin{equation}\label{eq:DescritR3S2}
\boxed{
\begin{cases}
\mathbf{X}_M = \mathbf{X}_0 + \sum\limits_{k=1}^{M} \sqrt{\frac{t D_{33}}{M}} \epsilon_k \mathbf{N}_{k-1},\\
\mathbf{N}_M = \left(\prod\limits_{k=1}^{M} \mathbf{R}_{\mathbf{a}, \gamma_k} \mathbf{R}_{\mathbf{e}_y, \beta_k \sqrt{\frac{t D_{44}}{M}}} \right)  \mathbf{N}_0,
\end{cases}
}
\end{equation}
with $\epsilon_{k} \sim  G_{t=1}^{\R} \sim \mathcal{N}(0,\sigma = \sqrt{2})$ stochastically independent Gaussian distributed on $\R$ with $t=1$; with uniformly distributed $\gamma_{k}\sim \textrm{Unif}\left( \mathbb{R}/(2\pi \mathbb{Z})\equiv [-\pi,\pi) \right)$;
and with $\beta_{k} \sim g$, where $g: \R \to \R^+$ equals $g(\beta)= \frac{|\beta|}{2}\; e^{-\frac{\beta^2}{4}}$. The choice of $g$ is done by application of the theory of isotropic stochastic processes on Riemannian manifolds by Pinsky \cite{pinsky_isotropic_1976}, where we note that
$$G_{t}^{\R^2}(\beta \cos\gamma, \beta \sin\gamma) = g(\beta) \; \textrm{Unif}\left([-\pi,\pi)\right)(\gamma),\quad \beta \in \R,\, \gamma \in [-\pi,\pi).$$
Now, in the numerical simulation we can replace $g$ by $G_{t=2}^{\R}$ due to the central limit theorem on $\R$ and
$$\textrm{Var}(\beta) = \int\limits_{-\infty}^{\infty} \beta^2 g(\beta) {\rm d} \beta = 2 \int\limits_{0}^{\infty} \beta^2 g(\beta) {\rm d} \beta = 2.$$

\subsubsection{\textbf{From the diffusion case $\alpha=1$ to the general case $\alpha \in (0,1]$}}\label{app:A21}

For the case $\alpha \in (0,1]$ we define the (fractional) random processes by their probability densities
\begin{equation} \label{Palpha}
\begin{array}{l}
P(\mathbf{P}^{\alpha}_t=(\mathbf{x},\mathbf{n}))= \int \limits_{0}^{\infty} q_{t,\alpha}(\tau) \; P(\mathbf{P}_\tau=(\mathbf{x},\mathbf{n}))\, {\rm d}\tau, \\
P(\overline{\mathbf{P}}^{\alpha}_t=(\mathbf{x},\overline{\mathbf{n}}))= \int \limits_{0}^{\infty} q_{t,\alpha}(\tau) \;  P(\overline{\mathbf{P}}_\tau=(\mathbf{x},\overline{\mathbf{n}}))\, {\rm d}\tau.
\end{array}
\end{equation}
recall that
the kernel $q_{t,\alpha}(\tau)$  was given by (\ref{eq:qker}).
For Monte-Carlo simulations one can use (\ref{kernelrel}), or alternatively 
use
$ \bP_{t_M}^{\alpha} \approx \prod\limits_{i=1}^{M} \bG_{T_i}\odot\bP_{0}, \textrm{ for } M \gg 0$, where $\bP_0$ is almost surely $(\mathbf{0},\mathbf{a})$, with $T_i$ a temporal random variable with $
P(T_{i}=\tau)=
q_{t_i,\alpha}(\tau)$, with $t_i = \frac{i}{M} t$ and $\mathbf{G}_{t_i}$ given by (\ref{eq:62alt}).

%

%
%

\subsubsection{\textbf{$\alpha$-Stability of the L\'{e}vy Process}}

Due to the absence of suitable dilations on $G/H$, we resort
to the embedding space where $\alpha$-stability is defined.
The L\'{e}vy process
$\{\overline{\mathbf{P}}_{t}^{\alpha}=(\mathbf{X}_{t}^{\alpha},\overline{\mathbf{N}}_{t}^{\alpha})\; | \; t \geq 0\}$ associated to the
L\'{e}vy process
$\{\mathbf{P}_{t}^{\alpha}=(\mathbf{X}_{t}^{\alpha},\mathbf{N}_{t}^{\alpha})\; |\; t \geq 0\}$ in $\R^{3}\rtimes S^{2}$
is $\alpha$-stable, i.e. for all $a,t>0$ we have (by (\ref{ene}) and (\ref{kernelrel}))
\[
\begin{array}{l}
a^{-\frac{1}{2\alpha}}\mathbf{X}_{at}^{\alpha} \sim \mathbf{X}_{t}^{\alpha} \textrm{ and } \
a^{-\frac{1}{2\alpha}}\overline{\mathbf{N}}_{at}^{\alpha} \sim \overline{\mathbf{N}}^{\alpha}_{t}.
\end{array}
\]

\section{Left-Invariant Vector fields  on SE(3) via 2 Charts \label{ch:LI}}

We need two charts to cover $SO(3)$.
When using the following coordinates (ZYZ-Euler angles) for $SE(3)=\R^{3}\rtimes SO(3)$ for the first chart:
\begin{equation} \label{Eulerangles1}
g=(x,y,z,\mathbf{R}_{\mathbf{e}_{z},\gamma} \mathbf{R}_{\mathbf{e}_{y},\beta} \mathbf{R}_{\mathbf{e}_{z},\overline{\alpha}}), \textrm{ with }\beta \in (0,\pi), \overline{\alpha},\gamma \in [0,2\pi),
\end{equation}
formula (\ref{li}) yields the following formulas for the left-invariant vector fields:
\begin{equation}
\begin{array}{l}
\!\!\mathcal{A}_1|_g = (\cos \overline{\alpha} \cos \beta \cos \gamma - \sin \overline{\alpha} \sin \gamma) \partial_x
+ (\sin \overline{\alpha} \cos \gamma + \cos \overline{\alpha} \cos \beta \sin \gamma) \partial_y - \cos \overline{\alpha} \sin \beta \, \partial_z \\
\!\!\mathcal{A}_2|_g = (- \sin \overline{\alpha} \cos \beta \cos \gamma - \cos \overline{\alpha} \sin \gamma) \partial_x
+ (\cos \overline{\alpha} \cos \gamma - \sin \overline{\alpha} \cos \beta \sin \gamma) \partial_y + \sin \overline{\alpha} \sin \beta \, \partial_z \\
\!\!\mathcal{A}_3|_g = \sin \beta \cos \gamma \,\partial_x + \sin \beta \sin \gamma \,\partial_y + \cos \beta \, \partial_z,
\\[6pt]
\!\!\mathcal{A}_4|_g =\; \cos \overline{\alpha} \textrm{cot} \beta \, \partial_{\overline{\alpha}} + \sin \overline{\alpha}\, \partial_\beta - \frac{\cos \overline{\alpha}}{\sin \beta}\, \partial_\gamma,\\
\!\!\mathcal{A}_5|_g = -\sin \overline{\alpha} \textrm{cot} \beta \,\partial_{\overline{\alpha}} + \cos \overline{\alpha}\, \partial_\beta + \frac{\sin \overline{\alpha}}{\sin \beta}\, \partial_\gamma,\\
\!\!\mathcal{A}_6|_g = \partial_{\overline{\alpha}}.
\end{array}
\end{equation}
We observe that
\begin{equation} \label{Za}
\underline{\mathcal{A}}_{g h_{\overline{\alpha}}} \equiv
(\bR_{\mathbf{e}_z,\overline{\alpha}} \oplus \bR_{\mathbf{e}_z,\overline{\alpha}})^T \underline{\mathcal{A}}_{g}, \quad \textrm{ where } \underline{\mathcal{A}}_{g} = \left(\mathcal{A}_1|_g, \ldots, \mathcal{A}_6|_g\right).
\end{equation}
\iftoggle{insertremarks}
{
The reason for the equivalence sign is that the tangent vectors are attached to different points in the group and for the equivalence we apply ordinary parallel transport
from
$\left.(\partial_{\overline{\alpha}},\partial_{\beta},\partial_{\gamma})\right|_{g h_{\overline{\alpha}}}$ to $\left.(\partial_{\overline{\alpha}},\partial_{\beta},\partial_{\gamma})\right|_{g}$
to identify the two tangent vectors.
}{}

The above formula's do not hold for $\beta=\pi$ or $\beta=0$. So we even lack expressions for our left-invariant vector fields at the unity element $(\mathbf{0},\mathbf{I}) \in SE(3)$ when using the standard ZYZ-Euler angles.
Therefore, one formally needs a second chart, for example the XYZ-coordinates in
\cite{portegies_new_2015,Duits2011,Duits2016JDCS}:
\begin{equation} \label{Eulerangles1}
g=(x,y,z,\mathbf{R}_{\mathbf{e}_{x},\tilde{\gamma}} \mathbf{R}_{\mathbf{e}_{y},\tilde{\beta}} \mathbf{R}_{\mathbf{e}_{z},\overline{\alpha}}), \  \textrm{ with }\tilde{\beta} \in [-\pi,\pi), \overline{\alpha} \in [0,2\pi),\ \tilde{\gamma} \in (-\pi/2,\pi/2),
\end{equation}
formula (\ref{li}) yields the following formulas for the left-invariant vector fields (only for $|\tilde{\beta}| \neq \frac{\pi}{2}$):
\begin{equation}
\begin{array}{l}
\!\!\mathcal{A}_1|_g = \cos \overline{\alpha} \cos \tilde{\beta} \, \partial_x
+ (\cos \tilde{\gamma} \sin \overline{\alpha} + \cos \overline{\alpha} \sin \tilde{\beta} \sin \tilde{\gamma})\, \partial_y +
(\sin \overline{\alpha} \sin \tilde{\gamma}- \cos \overline{\alpha} \sin \tilde{\beta} \cos{\tilde{\gamma}})\, \partial_z
\\
\!\!\mathcal{A}_2|_g = - \sin \overline{\alpha} \cos \tilde{\beta}\, \partial_x
+ (\cos \overline{\alpha} \cos \tilde{\gamma} - \sin \overline{\alpha} \sin \tilde{ \beta} \sin \tilde{\gamma}) \partial_y +
(\sin \overline{\alpha} \sin \tilde{\beta} \cos \tilde{\gamma} +\cos \overline{\alpha} \sin \tilde{\gamma})\,
\partial_z \\
\!\!\mathcal{A}_3|_g = \sin \tilde{\beta}\, \partial_x -
\cos \tilde{\beta} \sin \tilde{\gamma}\, \partial_y + \cos \tilde{\beta} \cos \tilde{\gamma} \, \partial_z,
\\[6pt]
\!\!\mathcal{A}_4|_g = -\cos \overline{\alpha} \textrm{tan} \tilde{\beta} \, \partial_{\overline{\alpha}} + \sin \overline{\alpha} \, \partial_{\tilde{\beta}} + \frac{\cos \overline{\alpha}}{\cos \tilde{\beta}}\partial_{\tilde{\gamma}},\\
\!\!\mathcal{A}_5|_g = \sin \overline{\alpha} \textrm{tan} \tilde{\beta}\,  \partial_{\overline{\alpha}} + \cos \overline{\alpha} \partial_{\tilde{\beta}} - \frac{\sin \overline{\alpha}}{\cos \beta}\, \partial_{\tilde{\gamma}},\\
\!\!\mathcal{A}_6|_g = \partial_{\overline{\alpha}}.
\end{array}
\end{equation}


%

\reftitle{References}
\bibliography{references}

\begin{thebibliography}{-------}
\providecommand{\natexlab}[1]{#1}

\bibitem[Zettl(2005)]{Zettl}
Zettl, A.
\newblock {\em Sturm-Liouville Theory}; Vol. 121, {\em Mathematical Surveys and
  Monographs}, American Mathematical Society,  2005.

\bibitem[Kato(1976)]{kato_operators_1976}
Kato, T.
\newblock Operators in {Hilbert} spaces. In {\em Perturbation {Theory} for
  {Linear} {Operators}}; Classics in {Mathematics}, Springer Berlin Heidelberg,
   1976; pp. 251--308.

\bibitem[Rudin(1991)]{Rudin}
Rudin, W.
\newblock {\em Functional Analysis}, 2 ed.; McGraw-Hill, Inc.,  1991.

\bibitem[Chirikjian and Kyatkin(2000)]{chirikjian_engineering_2000}
Chirikjian, G.S.; Kyatkin, A.B.
\newblock {\em Engineering {Applications} of {Noncommutative} {Harmonic}
  {Analysis}: {With} {Emphasis} on {Rotation} and {Motion} {Groups}}; CRC
  Press,  2000.

\bibitem[Chirikjian(2011)]{chirikjian_stochastic_2011}
Chirikjian, G.S.
\newblock {\em Stochastic {Models}, {Information} {Theory}, and {Lie} {Groups}:
  {Analytic} {Methods} and {Modern} {Applications}}; Vol.~2, Springer Science
  \& Business Media,  2011.

\bibitem[Saccon \em{et~al.}(2012)Saccon, Aguiar, Hausler, Hauser, and
  Pascoal]{Saccon}
Saccon, A.; Aguiar, A.P.; Hausler, A.J.; Hauser, J.; Pascoal, A.M.
\newblock Constrained motion planning for multiple vehicles on SE(3).
\newblock  2012 IEEE 51st IEEE Conference on Decision and Control (CDC),  2012,
  pp. 5637--5642.

\bibitem[Henk~Nijmeijer(1990)]{Nijmijer}
Henk~Nijmeijer, A.v.d.S.
\newblock {\em Nonlinear Dynamical Control Systems};  1990; p. 426.

\bibitem[Ali \em{et~al.}(1999)Ali, Antoine, and Gazeau]{Alibook}
Ali, S.; Antoine, J.; Gazeau, J.
\newblock {\em Coherent States, Wavelets and Their Generalizations}; Springer
  Verlag: New York, Berlin, Heidelberg,  1999.

\bibitem[Bekkers \em{et~al.}(2018)Bekkers, Lafarge, Veta, Eppenhof, Pluim, and
  Duits]{BekkersMICCAI}
Bekkers, E.; Lafarge, M.; Veta, M.; Eppenhof, K.; Pluim, J.; Duits, R.
\newblock Roto-Translation Covariant Convolutional Networks for Medical Image
  Analysis.
\newblock  Medical Image Computing and Computer Assisted Intervention -- MICCAI
  2018; et~al., F., Ed.; Springer International Publishing: Cham,  2018; pp.
  440--448.

\bibitem[Bekkers \em{et~al.}()Bekkers, Loog, ter Haar~Romeny, and
  Duits]{bekkers_template_2018}
Bekkers, E.; Loog, M.; ter Haar~Romeny, B.; Duits, R.
\newblock Template matching via densities on the roto-translation group.
\newblock {\em 40},~452--466.

\bibitem[Cohen \em{et~al.}(2018)Cohen, Geiger, and
  Weiler]{cohen2018intertwiners}
Cohen, T.S.; Geiger, M.; Weiler, M.
\newblock Intertwiners between Induced Representations (with Applications to
  the Theory of Equivariant Neural Networks).
\newblock {\em arXiv preprint arXiv:1803.10743} {\bf 2018}.

\bibitem[Cohen and Welling()]{Wellink}
Cohen, T.; Welling, M.
\newblock Group equivariant convolutional networks.
\newblock  Int. Conf. on Machine Learning, pp. 2990--2999.

\bibitem[Sifre and Mallat(2013)]{Mallat}
Sifre, L.; Mallat, S.
\newblock Rotation, scaling and deformation invariant scattering for texture
  discrimination.
\newblock  CVPR. IEEE,  2013, pp. 1233--1240.

\bibitem[Duits \em{et~al.}(2006)Duits, Felsberg, Granlund, and ter
  Haar~Romeny]{duits_image_2006}
Duits, R.; Felsberg, M.; Granlund, G.; ter Haar~Romeny, B.
\newblock Image {Analysis} and {Reconstruction} using a {Wavelet} {Transform}
  {Constructed} from a {Reducible} {Representation} of the {Euclidean} {Motion}
  {Group}.
\newblock {\em Int J Comput Vision} {\bf 2006}, {\em 72},~79--102.

\bibitem[Citti and Sarti(2006)]{citti_cortical_2006}
Citti, G.; Sarti, A.
\newblock A {Cortical} {Based} {Model} of {Perceptual} {Completion} in the
  {Roto}-{Translation} {Space}.
\newblock {\em J Math Imaging Vis} {\bf 2006}, {\em 24},~307--326.

\bibitem[Duits \em{et~al.}(2013)Duits, Fuehr, Janssen, Florack, and van
  Assen]{DuitsACHA}
Duits, R.; Fuehr, H.; Janssen, B.; Florack, L.; van Assen, H.
\newblock Evolution equations on Gabor transforms and their applications.
\newblock {\em ACHA} {\bf 2013}, {\em 35},~483--526.

\bibitem[Prandi and Gauthier(2018)]{prandigauthierbook}
Prandi, D.; Gauthier, J.P.
\newblock {\em A Semidiscrete Version of the Citti-Petitot-Sarti Model as a
  Plausible Model for Anthropomorphic Image Reconstruction and Pattern
  Recognition};  2018; p. 113.

\bibitem[Janssen \em{et~al.}(2018)Janssen, Janssen, Bekkers, Besc{\'o}s, and
  Duits]{Janssen2018}
Janssen, M.H.J.; Janssen, A.J.E.M.; Bekkers, E.J.; Besc{\'o}s, J.O.; Duits, R.
\newblock Design and Processing of Invertible Orientation Scores of 3D Images.
\newblock {\em Journal of Mathematical Imaging and Vision} {\bf 2018}.

\bibitem[Boscain \em{et~al.}(2012)Boscain, Duplaix, Gauthier, and
  Rossi]{boscain_anthropomorphic_2012}
Boscain, U.; Duplaix, J.; Gauthier, J.; Rossi, F.
\newblock Anthropomorphic {Image} {Reconstruction} via {Hypoelliptic}
  {Diffusion}.
\newblock {\em SIAM J. Control Optim.} {\bf 2012}, {\em 50},~1309--1336.

\bibitem[Schur(1968)]{Schur}
Schur, I.
\newblock {\em {V}orlesungen {\"{u}}ber {I}nvariantentheorie}; P. Noordhoff:
  Groningen,  1968.

\bibitem[Dieudonn{\'{e}}(1977)]{Dieudonne}
Dieudonn{\'{e}}, J.
\newblock {\em Treatise on Analysis}; Vol.~V,  1977.

\bibitem[Folland(1994)]{folland_course_1994}
Folland, G.B.
\newblock {\em A {Course} in {Abstract} {Harmonic} {Analysis}}; CRC Press,
  1994.

\bibitem[Agrachev \em{et~al.}(2009)Agrachev, Boscain, Gauthier, and
  Rossi]{agrachev_intrinsic_2009}
Agrachev, A.; Boscain, U.; Gauthier, J.P.; Rossi, F.
\newblock The intrinsic hypoelliptic {Laplacian} and its heat kernel on
  unimodular {Lie} groups.
\newblock {\em J Funct Anal} {\bf 2009}, {\em 256},~2621--2655.

\bibitem[F{\"u}hr(2005)]{fuhr_abstract_2005}
F{\"u}hr, H.
\newblock {\em Abstract {Harmonic} {Analysis} of {Continuous} {Wavelet}
  {Transforms}}; Springer Science \& Business Media,  2005.

\bibitem[Mackey(1949)]{mackey_imprimitivity_1949}
Mackey, G.W.
\newblock Imprimitivity for {Representations} of {Locally} {Compact} {Groups}
  {I}.
\newblock {\em Proc Natl Acad Sci U S A} {\bf 1949}, {\em 35},~537--545.

\bibitem[Sugiura(1990)]{sugiura_unitary_1990}
Sugiura, M.
\newblock {\em Unitary {Representations} and {Harmonic} {Analysis}: {An}
  {Introduction}}; Elsevier,  1990.

\bibitem[Dixmier(1981)]{Dixmier}
Dixmier, J.
\newblock {\em $C^*$-algebras}; North Holland,  1981.

\bibitem[Gaveau(1977)]{Gaveau1977}
Gaveau, B.
\newblock Principe de moindre action, propagation de la chaleur et estimees
  sous elliptiques sur certains groupes nilpotents.
\newblock {\em Acta Math.} {\bf 1977}, {\em 139},~95--153.

\bibitem[Duits and van Almsick(2008)]{duits_explicit_2008-1}
Duits, R.; van Almsick, M.
\newblock The explicit solutions of linear left-invariant second order
  stochastic evolution equations on the 2D {Euclidean} motion group.
\newblock {\em Quarterly of Applied Mathematics} {\bf 2008}, {\em 66},~27--67.

\bibitem[Duits and Franken(2009)]{duits_line_2009}
Duits, R.; Franken, E.
\newblock Line {Enhancement} and {Completion} via {Linear} {Left} {Invariant}
  {Scale} {Spaces} on {SE}(2).
\newblock  {SSVM}; Springer-Verlag: Berlin, Heidelberg,  2009; pp. 795--807.

\bibitem[Duits and van Almsick(2005)]{DuitsCASA2005}
Duits, R.; van Almsick.
\newblock The explicit solutions of linear left-invariant second order
  stochastic evolution equations on the 2D-{E}uclidean motion group.
\newblock Technical Report CASA-report, nr.43, 37 pages., Eindhoven University
  of Technology Dep. of mathematics and computer science,  2005.
\newblock \url{http://www.win.tue.nl/analysis/reports/rana05-43.pdf}.

\bibitem[Duits and Franken(2010)]{duits_left-invariant_2010-1}
Duits, R.; Franken, E.
\newblock Left-invariant parabolic evolutions on {SE}(2) and contour
  enhancement via invertible orientation scores {Part} {II}: {Nonlinear}
  left-invariant diffusions on invertible orientation scores.
\newblock {\em Quart. Appl. Math.} {\bf 2010}, {\em 68},~293--331.

\bibitem[Zhang \em{et~al.}(2016)Zhang, Duits, Sanguinetti, and ter
  Haar~Romeny]{zhang_numerical_2016}
Zhang, J.; Duits, R.; Sanguinetti, G.; ter Haar~Romeny, B.M.
\newblock Numerical {Approaches} for {Linear} {Left}-invariant {Diffusions} on
  {SE}(2), their {Comparison} to {Exact} {Solutions}, and their {Applications}
  in {Retinal} {Imaging}.
\newblock {\em Numerical Methods Theory and Applications (NM-TMA)} {\bf 2016},
  {\em 9},~1--50.

\bibitem[Mumford(1994)]{mumford_elastica_1994-1}
Mumford, D.
\newblock Elastica and {Computer} {Vision}. In {\em Algebraic {Geometry} and
  its {Applications}}; Springer New York,  1994; pp. 491--506.

\bibitem[Petitot(2003)]{petitot_neurogeometry_2003}
Petitot, J.
\newblock The neurogeometry of pinwheels as a sub-{Riemannian} contact
  structure.
\newblock {\em Journal of Physiology-Paris} {\bf 2003}, {\em 97},~265--309.

\bibitem[Agrachev \em{et~al.}(2009)Agrachev, Boscain, Gauthier, and
  Rossi]{agrachev_intrinsic_2009-1}
Agrachev, A.; Boscain, U.; Gauthier, J.P.; Rossi, F.
\newblock The intrinsic hypoelliptic {Laplacian} and its heat kernel on
  unimodular {Lie} groups.
\newblock {\em Journal of Functional Analysis} {\bf 2009}, {\em 256},~2621 --
  2655.

\bibitem[Portegies and Duits(2017)]{portegies_new_2017}
Portegies, J.M.; Duits, R.
\newblock New exact and numerical solutions of the
  (convection{\textendash})diffusion kernels on {SE}(3).
\newblock {\em Differential Geometry and its Applications} {\bf 2017}, {\em
  53},~182--219.

\bibitem[Portegies \em{et~al.}(2015)Portegies, Fick, Sanguinetti, Meesters,
  Girard, and Duits]{portegies_improving_2015}
Portegies, J.M.; Fick, R.H.J.; Sanguinetti, G.R.; Meesters, S.P.L.; Girard, G.;
  Duits, R.
\newblock Improving {Fiber} {Alignment} in {HARDI} by {Combining} {Contextual}
  {PDE} {Flow} with {Constrained} {Spherical} {Deconvolution}.
\newblock {\em PLoS ONE} {\bf 2015}, {\em 10},~e0138122.

\bibitem[Momayyez-Siahkal and Siddiqi(2009)]{momayyez-siahkal_3d_2009}
Momayyez-Siahkal, P.; Siddiqi, K.
\newblock 3D {Stochastic} {Completion} {Fields} for {Fiber} {Tractography}.
\newblock  Proc {IEEE} {Comput} {Soc} {Conf} {Comput} {Vis} {Pattern}
  {Recognit},  2009, pp. 178--185.

\bibitem[Skibbe and Reisert(2017)]{Skibbe}
Skibbe, H.; Reisert, M.
\newblock Spherical Tensor Algebra: A Toolkit for 3D Image Processing.
\newblock {\em Journal of Mathematical Imaging and Vision} {\bf 2017}, {\em
  58},~349--381.

\bibitem[Meesters \em{et~al.}(2017)Meesters, Ossenblok, Wagner, Schijns, Boon,
  Florack, Vilanova, and Duits]{MEESTERS2017}
Meesters, S.; Ossenblok, P.; Wagner, L.; Schijns, O.; Boon, P.; Florack, L.;
  Vilanova, A.; Duits, R.
\newblock Stability metrics for optic radiation tractography: Towards damage
  prediction after resective surgery.
\newblock {\em Journal of Neuroscience Methods} {\bf 2017}, {\em 288},~34 --
  44.

\bibitem[Reisert and Kiselev(2011)]{reisert_fiber_2011}
Reisert, M.; Kiselev, V.G.
\newblock Fiber {Continuity}: {An} {Anisotropic} {Prior} for {ODF}
  {Estimation}.
\newblock {\em IEEE Transactions on Medical Imaging} {\bf 2011}, {\em
  30},~1274--1283.

\bibitem[Pr{\v c}kovska \em{et~al.}(2010)Pr{\v c}kovska, Rodrigues, Duits,
  Haar~Romenij, and Vilanova]{prckovska_extrapolating_2010}
Pr{\v c}kovska, V.; Rodrigues, P.; Duits, R.; Haar~Romenij, B.t.; Vilanova, A.
\newblock Extrapolating fiber crossings from {DTI} data: can we infer similar
  fiber crossings as in {HARDI}?
\newblock  {MICCAI},  2010, Workshop on {Computational} {Diffusion} {MRI}.

\bibitem[Iijima(1959)]{Iiji59a}
Iijima, T.
\newblock Basic theory of pattern observation.
\newblock {\em Technical Group on Automata and Automatic Control, IECE, Japan}
  {\bf 1959}.

\bibitem[Koenderink(1984)]{Koenderink}
Koenderink, J.J.
\newblock The structure of images.
\newblock {\em Biological Cybernetics} {\bf 1984}, {\em 50},~363--370.

\bibitem[ter Haar~Romeny(2003)]{ter_haar_romeny_front-end_2003}
ter Haar~Romeny, B.M.
\newblock {\em Front-end {Vision} and {Multi}-{Scale} {Image} {Analysis}:
  {Multi}-{Scale} {Computer} {Vision} {Theory} and {Applications}, {Written} in
  {Mathematica}}; Kluwer Academic Publishers,  2003.

\bibitem[Weickert(1998)]{weickert_anisotropic_1998}
Weickert, J.
\newblock {\em Anisotropic {Diffusion} in {Image} {Processing}};  1998.

\bibitem[Duits and Burgeth(2007)]{duits_scale_2007}
Duits, R.; Burgeth, B.
\newblock Scale {Spaces} on {Lie} {Groups}.
\newblock  {SSVM}. Springer,  2007, Vol. 4485, {\em Lecture {Notes} in
  {Computer} {Science}}, pp. 300--312.

\bibitem[Benoist and Quint(2016)]{Benoist}
Benoist, Y.; Quint, J.F.
\newblock Central limit theorem for linear groups.
\newblock {\em Annals of Probability} {\bf 2016}, {\em 44},~1306--1340.

\bibitem[Pilt\'{e} \em{et~al.}(2018)Pilt\'{e}, Bonnabel, and
  Barbaresco]{Barbaresco}
Pilt\'{e}, M.; Bonnabel, S.; Barbaresco, F.
\newblock Maneuver Detector for Active Tracking Update Rate Adaptation.
\newblock  2018 19th International Radar Symposium (IRS),  2018, pp. 1--10.

\bibitem[Berger \em{et~al.}(2015)Berger, Neufeld, Becker, Lenzen, and
  Schn{\"o}rr]{Berger}
Berger, J.; Neufeld, A.; Becker, F.; Lenzen, F.; Schn{\"o}rr, C.
\newblock Second Order Minimum Energy Filtering on {SE(3)} with Nonlinear
  Measurement Equations.
\newblock  Scale Space and Variational Methods in Computer Vision; Aujol, J.F.;
  Nikolova, M.; Papadakis, N., Eds.; Springer: Cham,  2015; pp. 397--409.

\bibitem[Oksendal(1998)]{Oksendal}
Oksendal, B.
\newblock {\em Stochastic differential equations}; Universitext, Springer,
  1998.

\bibitem[Hsu(2002)]{Hsu}
Hsu, E.
\newblock {\em Stochastic Analysis on Manifolds}; Contemporary Mathematics,
  American Mathematical Society,  2002.

\bibitem[Feller(1966)]{Feller}
Feller, W.
\newblock {\em An Introduction to Probability Theory and Its Applications};
  Vol.~II, Wiley Series in Probability and Mathematical Statistics,  1966.

\bibitem[Felsberg \em{et~al.}(2003)Felsberg, Duits, and Florack]{Fels2003}
Felsberg, M.; Duits, R.; Florack, L.
\newblock The Monogenic Scale Space on a Bounded Domain and its Applications.
\newblock {\em Proceedings Scale Space Conference, volume 2695 of Lecture Notes
  of Computer Science, Springer, Isle of Skye, UK.} {\bf 2003}, pp. 209--224.

\bibitem[Duits \em{et~al.}(2003)Duits, Felsberg, and Florack]{Duits2003a}
Duits, R.; Felsberg, M.; Florack, L.M.J.
\newblock $\alpha$ Scale Spaces on a Bounded Domain.
\newblock {\em Proceedings Scale Space Conference, Isle of Skye,UK.} {\bf
  2003}, pp. 494--510.

\bibitem[Duits \em{et~al.}(2004)Duits, Florack, Graaf, and
  Romeny]{duits_axioms_2004}
Duits, R.; Florack, L.; Graaf, J.d.; Romeny, B.t.H.
\newblock On the {Axioms} of {Scale} {Space} {Theory}.
\newblock {\em Journal of Mathematical Imaging and Vision} {\bf 2004}, {\em
  20},~267--298.

\bibitem[Pedersen \em{et~al.}(2005)Pedersen, Duits, and Nielsen]{Pedersen}
Pedersen, K.S.; Duits, R.; Nielsen, M.
\newblock On $\alpha$ Kernels, L{\'e}vy Processes, and Natural Image
  Statistics.
\newblock  Scale Space and PDE Methods in Computer Vision; Kimmel, R.; Sochen,
  N.A.; Weickert, J., Eds.; Springer Berlin Heidelberg: Berlin, Heidelberg,
  2005; pp. 468--479.

\bibitem[Yosida(1980)]{yosida_functional_1980}
Yosida, K.
\newblock {\em Functional {Analysis}}; Springer Berlin Heidelberg: Berlin,
  Heidelberg,  1980.

\bibitem[Winkels and Cohen(2018)]{winkels20183d}
Winkels, M.; Cohen, T.S.
\newblock 3D G-CNNs for Pulmonary Nodule Detection.
\newblock {\em arXiv preprint arXiv:1804.04656} {\bf 2018}.

\bibitem[Worrall and Brostow(2018)]{worrall2018cubenet}
Worrall, D.; Brostow, G.
\newblock CubeNet: Equivariance to 3D Rotation and Translation.
\newblock {\em arXiv preprint arXiv:1804.04458} {\bf 2018}.

\bibitem[Weiler \em{et~al.}(2018)Weiler, Geiger, Welling, Boomsma, and
  Cohen]{weiler20183d}
Weiler, M.; Geiger, M.; Welling, M.; Boomsma, W.; Cohen, T.
\newblock 3D Steerable CNNs: Learning Rotationally Equivariant Features in
  Volumetric Data,  2018,
  \href{http://xxx.lanl.gov/abs/1807.02547}{{\normalfont
  [arXiv:cs.LG/1807.02547]}}.

\bibitem[Montobbio \em{et~al.}(2018)Montobbio, Sarti, and Citti]{CittiX}
Montobbio, N.; Sarti, A.; Citti, G.
\newblock A metric model for the functional architecture of the visual cortex
  {\bf 2018}.
\newblock  \href{http://xxx.lanl.gov/abs/1807.02479}{{\normalfont
  [arXiv:math.MG/1807.02479]}}.

\bibitem[Oyallon \em{et~al.}(2013)Oyallon, Mallat, and
  Sifre]{oyallon2013generic}
Oyallon, E.; Mallat, S.; Sifre, L.
\newblock Generic deep networks with wavelet scattering.
\newblock {\em arXiv preprint arXiv:1312.5940} {\bf 2013}.

\bibitem[Kanti V.~Mardia(1999)]{MardiaJuppBook1999}
Kanti V.~Mardia, P.E.J.
\newblock {\em Directional Statistics};  1999; p. 350.

\bibitem[Wu(2007)]{Wu}
Wu, L.
\newblock Chapter 3 Modeling Financial Security Returns Using Lévy Processes.
  In {\em Financial Engineering}; Birge, J.R.; Linetsky, V., Eds.; Elsevier,
  2007; Vol.~15, {\em Handbooks in Operations Research and Management Science},
  pp. 117 -- 162.

\bibitem[Belkic and Belkic(2010)]{BelkicandBelkic2010CRCPress}
Belkic, D.D.; Belkic, K.
\newblock {\em Signal processing in magnetic resonance spectroscopy with
  biomedical applications}; Boca Raton : CRC Press,  2010.

\bibitem[Duits and Franken(2010)]{duits_left-invariant_2010-2}
Duits, R.; Franken, E.
\newblock Left-invariant parabolic evolutions on {SE}(2) and contour
  enhancement via invertible orientation scores {Part} {I}: {Linear}
  left-invariant diffusion equations on {SE}(2).
\newblock {\em Quart. Appl. Math.} {\bf 2010}, {\em 68},~255--292.

\bibitem[Citti and Sarti(2006)]{citti_cortical_2006-1}
Citti, G.; Sarti, A.
\newblock A {Cortical} {Based} {Model} of {Perceptual} {Completion} in the
  {Roto}-{Translation} {Space}.
\newblock {\em Journal of Mathematical Imaging and Vision} {\bf 2006}, {\em
  24},~307--326.

\bibitem[H{\"o}rmander(1967)]{hormander_hypoelliptic_1967}
H{\"o}rmander, L.
\newblock Hypoelliptic second order differential equations.
\newblock {\em Acta Mathematica} {\bf 1967}, {\em 119},~147--171.

\bibitem[Misiorek and Weron(2012)]{Misiorek2012}
Misiorek, A.; Weron, R., Heavy-Tailed Distributions in VaR Calculations.
\newblock In {\em Handbook of Computational Statistics: Concepts and Methods};
  Gentle, J.E.; H{\"a}rdle, W.K.; Mori, Y., Eds.; Springer Berlin Heidelberg:
  Berlin, Heidelberg,  2012; pp. 1025--1059.

\bibitem[Felsberg and Sommer(2004)]{felsberg_monogenic_2004}
Felsberg, M.; Sommer, G.
\newblock The {Monogenic} {Scale}-{Space}: {A} {Unifying} {Approach} to
  {Phase}-{Based} {Image} {Processing} in {Scale}-{Space}.
\newblock {\em Journal of Mathematical Imaging and Vision} {\bf 2004}, {\em
  21},~5--26.

\bibitem[Kanters \em{et~al.}(2007)Kanters, Florack, Duits, Platel, and ter
  Haar~Romeny]{Kanters2007}
Kanters, F.; Florack, L.; Duits, R.; Platel, B.; ter Haar~Romeny, B.
\newblock ScaleSpaceViz: $\alpha$-Scale spaces in practice.
\newblock {\em Pattern Recognition and Image Analysis} {\bf 2007}, {\em
  17},~106--116.

\bibitem[Schmidt and Weickert(2016)]{schmidt_morphological_2016}
Schmidt, M.; Weickert, J.
\newblock Morphological {Counterparts} of {Linear} {Shift}-{Invariant}
  {Scale}-{Spaces}.
\newblock {\em J Math Imaging Vis} {\bf 2016}, {\em 56},~352--366.

\bibitem[Duits and Franken(2011)]{Duits2011}
Duits, R.; Franken, E.
\newblock Left-Invariant Diffusions on the Space of Positions and Orientations
  and their Application to Crossing-Preserving Smoothing of HARDI images.
\newblock {\em International Journal of Computer Vision} {\bf 2011}, {\em
  92},~231--264.

\bibitem[Duits \em{et~al.}(2012)Duits, Dela~Haije, Creusen, and
  Ghosh]{duits_morphological_2012}
Duits, R.; Dela~Haije, T.; Creusen, E.; Ghosh, A.
\newblock Morphological and {Linear} {Scale} {Spaces} for {Fiber} {Enhancement}
  in {DW}-{MRI}.
\newblock {\em J Math Imaging Vis} {\bf 2012}, {\em 46},~326--368.

\bibitem[Portegies \em{et~al.}(2015)Portegies, Sanguinetti, Meesters, and
  Duits]{portegies_new_2015}
Portegies, J.; Sanguinetti, G.; Meesters, S.; Duits, R.
\newblock New {Approximation} of a {Scale} {Space} {Kernel} on {SE}(3) and
  {Applications} in {Neuroimaging}. In {\em {SSVM}}; Aujol, J.F.; Nikolova, M.;
  Papadakis, N., Eds.; Number 9087 in {LNCS}, Springer International
  Publishing,  2015; pp. 40--52.

\bibitem[Bukhvalov and W.(1994)]{Bukhvalov}
Bukhvalov, A.; W., A.
\newblock Integral representation of resolvent and semigroups.
\newblock {\em Forum Math. 6} {\bf 1994}, {\em 6},~111--137.

\bibitem[Liao(2004)]{MIngLiao}
Liao, M.
\newblock {\em Lévy Processes in Lie Groups}; Cambridge Tracts in Mathematics,
  Cambridge University Press,  2004.

\bibitem[Griffiths(1994)]{griffiths}
Griffiths, D.
\newblock {\em Introduction to Quantum Mechanics}; Prentice-Hall: USA,  1994.

\bibitem[Wigner(1931)]{Wigner}
Wigner, E.
\newblock Gruppentheorie und ihre Anwendungen auf die Quantenmechanik der
  Atomspektren.
\newblock {\em Braunschweig: Vieweg Verlag.} {\bf 1931}.
\newblock translated into English by Griffin, J. J. (1959). Group Theory and
  its Application to the Quantum Mechanics of Atomic Spectra.

\bibitem[Margenau and Murphy(1956)]{margenau_mathematics_1956}
Margenau, H.; Murphy, G.M.
\newblock {\em The mathematics of physics and chemistry}; Van Nostrand,  1956.

\bibitem[Pinsky(1976)]{pinsky_isotropic_1976}
Pinsky, M.A.
\newblock Isotropic transport process on a {Riemannian} manifold.
\newblock {\em Trans. Amer. Math. Soc.} {\bf 1976}, {\em 218},~353--360.

\bibitem[Duits and Franken(2010)]{duits_left-invariant_2010-3}
Duits, R.; Franken, E.
\newblock Left-invariant diffusions on the space of positions and orientations
  and their application to crossing-preserving smoothing of {HARDI} images.
\newblock {\em International Journal of Computer Vision} {\bf 2010}, {\em
  92},~231--264.

\bibitem[Meesters \em{et~al.}(2016)Meesters, Sanguinetti, Garyfallidis,
  Portegies, and Duits]{meesters_fast_2016}
Meesters, S.P.L.; Sanguinetti, G.R.; Garyfallidis, E.; Portegies, J.M.; Duits,
  R.
\newblock Fast implementations of contextual {PDE}'s for {HARDI} data
  processing in {DIPY}.
\newblock  2016.

\bibitem[Pr{\v c}kovska \em{et~al.}(2015)Pr{\v c}kovska, Andorr{\`a},
  Villoslada, Martinez-Heras, Duits, Fortin, Rodrigues, and
  Descoteaux]{prckovska_contextual_2015}
Pr{\v c}kovska, V.; Andorr{\`a}, M.; Villoslada, P.; Martinez-Heras, E.; Duits,
  R.; Fortin, D.; Rodrigues, P.; Descoteaux, M.
\newblock Contextual {Diffusion} {Image} {Post}-processing {Aids} {Clinical}
  {Applications}. In {\em Visualization and {Processing} of {Higher} {Order}
  {Descriptors} for {Multi}-{Valued} {Data}}; Hotz, I.; Schultz, T., Eds.;
  Mathematics and {Visualization}, Springer International Publishing,  2015;
  pp. 353--377.
\newblock DOI: 10.1007/978-3-319-15090-1\_18.

\bibitem[{Meesters, S.P.L.} \em{et~al.}(2016){Meesters, S.P.L.}, Sanguinetti,
  Garyfallidis, Portegies, Ossenblok, and Duits]{meesters_s.p.l._cleaning_2016}
{Meesters, S.P.L.}.; Sanguinetti, G.R.; Garyfallidis, E.; Portegies, J.M.;
  Ossenblok, P.; Duits, R.
\newblock Cleaning output of tractography via fiber to bundle coherence, a new
  open source implementation.
\newblock  2016.

\bibitem[E.~Barndorff-Nielsen \em{et~al.}(2001)E.~Barndorff-Nielsen, Mikosch,
  and I.~Resnick]{Nielsen}
E.~Barndorff-Nielsen, O.; Mikosch, T.; I.~Resnick, S.
\newblock {\em Lévy processes. Theory and applications}; Vol.~48,  2001; p.
  415.

\bibitem[Duits \em{et~al.}(2016)Duits, Ghosh, Dela~Haije, and
  Mashtakov]{Duits2016JDCS}
Duits, R.; Ghosh, A.; Dela~Haije, T.C.J.; Mashtakov, A.
\newblock On Sub-Riemannian Geodesics in SE(3) Whose Spatial Projections do not
  Have Cusps.
\newblock {\em Journal of Dynamical and Control Systems} {\bf 2016}, {\em
  22},~771--805.

\end{thebibliography}
\end{document}